\documentclass[11pt,preprint]{imsart}
\RequirePackage[OT1]{fontenc}
\RequirePackage{amsthm,amsmath}
\RequirePackage[numbers]{natbib}
\RequirePackage[colorlinks,citecolor=blue,urlcolor=blue]{hyperref}

\usepackage{graphicx,amsmath,amsthm,amsfonts,amssymb}
\usepackage{booktabs, rotating, float}
\usepackage[hmargin=3cm,vmargin=2.7cm]{geometry}
\usepackage{mathtools}
\usepackage{bbm}
\mathtoolsset{showonlyrefs}
\usepackage{enumerate}
\usepackage{enumitem}
\usepackage{color}
\usepackage{listings}
\usepackage{epstopdf}
\usepackage{makecell}
\usepackage{tikz}
\usetikzlibrary{decorations.pathreplacing}
\usepackage{array}
\usepackage{caption}
\usepackage{subcaption}
\usepackage{soul}
\usepackage{hyperref}
\usepackage{verbatim}
\usepackage{siunitx}
\usetikzlibrary{calc} \tikzset{>=latex}
\usetikzlibrary{calc,decorations.markings}
\usetikzlibrary{positioning,shapes,decorations.text}
\theoremstyle{plain}
\newtheorem{thm}{Theorem}[section]
\newtheorem{cor}[thm]{Corollary}

\newtheorem{prop}[thm]{Proposition}
\newtheorem{lemma}[thm]{Lemma}
\newtheorem{assumption}[thm]{Assumption}
\theoremstyle{definition}
\newtheorem{defn}[thm]{Definition}

\newtheorem{eg}[thm]{Example}
\newtheorem{fact}[thm]{Fact}
\newtheorem{observe}[thm]{Observation}

\numberwithin{equation}{section}

\newcommand{\rpm}{\sbox0{$1$}\sbox2{$\scriptstyle\pm$}
  \raise\dimexpr(\ht0-\ht2)/2\relax\box2 }

\newcommand{\R}{\mathbb{R}}

\newcommand{\ind}[1]{\mathbf{1}_{#1}}

\tikzstyle{nd} = [anchor=base, inner sep=0pt]
\tikzstyle{ndpic} = [remember picture, baseline, every node/.style={nd}]

\newcommand{\chia}{\chi^{\mathrm{Ai}}}
\def\beq{\begin{equation}}
\def\eeq{\end{equation}}
\def\ba{\begin{enumerate}[(a)]}
\def\bei{\begin{enumerate}[(i)]}
\def\be{\begin{enumerate}[(1)]}
\def\ee{\end{enumerate}}
\def\bi{\begin{itemize}}
\def\ei{\end{itemize}}
\def\beg{\begin{eg}}
\def\eeg{\end{eg}}
\def\bd{\begin{defn}}
\def\ed{\end{defn}}
\def\bt{\begin{thm}}
\def\et{\end{thm}}
\def\bl{\begin{lemma}}
\def\el{\end{lemma}}
\def\bfac{\begin{fact}}
\def\efac{\end{fact}}

\def\bc{\begin{cor}}
\def\ec{\end{cor}}
\def\bp{\begin{prop}}
\def\ep{\end{prop}}
\def\bo{\begin{observe}}
\def\eo{\end{observe}}
\def\bas{\begin{assumption}}
\def\eas{\end{assumption}}
\def\RR{\mathbb{R}}
\def\CC{\mathbb{C}}
\def\EE{\mathbb{E}}
\def\ZZ{\mathbb{Z}}

\def\PP{\mathbb{P}}

\def\beg{\begin{eg}}
\def\eeg{\end{eg}}

\def\fluc{{\Upsilon}}

\def\Ai{{\textrm{Ai}}}
\makeindex

\numberwithin{equation}{section}
\numberwithin{table}{section}

\startlocaldefs
\endlocaldefs

\begin{document}

\begin{frontmatter}

\title{Lower tail of the KPZ equation}
\runtitle{Lower tail of the KPZ equation}
\runauthor{Corwin \& Ghosal}

\begin{aug}
  \author{\fnms{Ivan}  \snm{Corwin}\textsuperscript{1}\corref{}\ead[label=e1]{corwin@math.columbia.edu}}
  \author{\fnms{Promit} \snm{Ghosal}\textsuperscript{2}\corref{}\ead[label=e2]{pg2475@columbia.edu}\vspace{0.1in}\\{Columbia University}}

\address{\textsuperscript{1}Department of Mathematics, 2990 Broadway, New York, NY 10027, USA\\ \printead{e1}}
\address{\textsuperscript{2}Department of Statistics, 1255 Amsterdam, New York, NY 10027, USA\\ \printead{e2}}




\end{aug}

\begin{abstract}
We provide the first tight bounds on the lower tail probability of the one point distribution of the KPZ equation with narrow wedge initial data. Our bounds hold for all sufficiently large times $T$ and demonstrates a crossover between super-exponential decay with exponent $5/2$ (and leading pre-factor $\frac{4}{15\pi} T^{1/3}$) for tail depth greater than $T^{2/3}$, and exponent $3$ (with leading pre-factor $\frac{1}{12}$) for tail depth less than $T^{2/3}$. 

\end{abstract}




\end{frontmatter}









\setcounter{tocdepth}{1}
\tableofcontents

\section{Introduction}\label{introduction}
The $(1+1)$d stochastic heat equation (SHE) with multiplicative space time white noise $\xi$ is
\begin{align}\label{eq:SHEDef}
\partial_T \mathcal{Z}(T,X) = \frac{1}{2}\partial^2_X \mathcal{Z}(T,X) +\mathcal{Z}(T,X) \xi(T,X),
\end{align}
where $T\geq 0$ and $X\in \R$. The solution theory for this stochastic PDE is classical \cite{Walsh, Corwin12, Quastel12}, based on  It\^{o} stochastic integrals or martingale problems.
The SHE is ubiquitous, modeling the density of particles diffusing in space-time random environments (with random killing / branching \cite{Mol96, Davar} or random drifts \cite{BarraquandCorwin15,CG17}). Via Feynman-Kac, it is the partition function for the continuum directed polymer model \cite{Amir11, Comets16, HuseHenley}. Taking logarithms formally leads to the Kardar-Parisi-Zhang (KPZ) equation
\begin{align}\label{eq:KPZDef}
\partial_T \mathcal{H}(T,X)= \frac{1}{2} \partial^2_X \mathcal{H}(T,X)+ \frac{1}{2} \big(\partial_X \mathcal{H}(T,X)\big)^2+ \xi(T,X),
\end{align}
which is a paradigm for random interface growth \cite{KPZ86} and a testing ground for the study of non-linear stochastic PDEs \cite{Hairer13, GJ14, GIP15, HQ15, GP17}. The KPZ equation's spatial derivative formally solves the stochastic Burgers equation -- a continuum model for turbulence \cite{FNS,BCJ94}, interacting particle systems and driven lattice gases \cite{BKS85}.

The {\it Cole-Hopf solution to the KPZ equation} with   {\it narrow wedge initial data} is given by
\begin{equation}\label{eq:flucs}
\mathcal{H}(T,X):=\log\mathcal{Z}(T,X),\qquad \textrm{with}\quad \mathcal{Z}(0,X)=\delta_{X=0}.
\end{equation}
The well-definedness of $\log\mathcal{Z}$ for all $T>0$ and $X\in \R$ relies upon the almost-sure strict positivity of the $\mathcal{Z}$ proved in \cite{CM91} to hold for a wide class of initial data (including the delta function). This is the physically relevant notion of solution and has been shown to arise quite generally from various regularization or discretization schemes for the equation and noise \cite{Bertini1995, Bertini1997,CL15,CS16,Hairer13,HS15,HQ15,GIP15,GP17,GJ14}. The Cole-Hopf solution also coincides with the solutions constructed from regularity structures \cite{Hairer13}, paracontrolled distributions \cite{GIP15} and energy solution methods \cite{GP17}.

\medskip

This paper establishes tight bounds on the {\it lower tail} probability that $\mathcal{Z}(T,X)$ is close to zero, or equivalently that $\mathcal{H}(T,X)$ is very negative\footnote{To avoid confusion, let us distinguish our present investigation from earlier work of \cite{EKMS,EV99,EV00} which studied the stochastic Burgers equation (the spatial derivative of KPZ) but with a noise which is smooth in space and white in time. In that case, which has no direct relationship to our work, the tail of the local slope has $-\frac{7}{2}$ power-law (not exponential) decay. A proxy for the question we consider here, \cite{AE95} studied the tail behavior of the invicid Burgers equation with white-noise initial data, showing cubic exponential decay. That result, however, also has no direct bearing on our present work.}. The first result in this direction was the aforementioned almost-sure positivity of $\mathcal{Z}$ established in \cite{CM91} via large deviation bounds and a comparison principal. Using Malliavin calculus, \cite{CN08} proved a quantitative upper bound on the decay of the lower tail probability. Working with the SHE on an interval with Dirichlet boundary conditions and constant initial data, they show that for any $\delta>0$ there are constants $c_1,c_2>0$ so that $\mathbb{P}\big(\mathcal{H}(T,X)\leq -s\big) \leq c_1 \exp\big(- c_2 s^{\frac{3}{2}-\delta}\big)$. Using Talagrand's concentration of measure methods, \cite{GMF14} improved the exponent. In particular, \cite{GMF14} considered the full-line SHE with $\mathcal{Z}(0,X)=\delta_{X=0}$ initial data (this is the setting we address in this paper) and proved a similar bound to \cite{CN08} but with the $3/2-\delta$ exponent replaced by the Gaussian exponent $2$. Quite recently, using Malliavin calculus \cite{HuLe} extended these sort of results to noises with more general covariance structure. There is some work in progress \cite{Davar} which seeks to use stochastic analytic methods to prove a lower bound with exponent $5/2$ on this tail probability. As we prove here, the exponents accessed in earlier work are not optimal and, moreover, these previous results are (in a sense we now describe) not well adapted to study the long (or intermediate) time solution tail.

When time increases, the KPZ equation shows an overall decay at linear rate $-T/24$ with fluctuations which grows like $T^{1/3}$. \cite{Amir11} proved (see also \cite{SS10} for a less rigorous treatment done in parallel, and \cite{CDR10,Dot10} for physics results) that when $\mathcal{Z}(0,X) = \delta_{X=0}$,
\begin{equation}\label{eq:eqnchit}
\lim_{T\to \infty}\mathbb{P}\big(\fluc_T \leq s\big) = F_{{\rm GUE}}(s),\qquad \textrm{where}\quad \fluc_T:=\frac{\mathcal{H}(2T,0) + \frac{T}{12}}{T^{\frac{1}{3}}}.
\end{equation}
The $T^{1/3}$ scaling is a characteristic of models in the KPZ universality class, as is the limiting GUE Tracy-Widom distribution $F_{{\rm GUE}}(s)$ \cite{Corwin12}. We consider $\mathcal{H}$ at time $2T$ to simplify some factors of 2 in formulas. Reinterpreting the tail bounds of \cite{CN08,GMF14} in terms of the lower tail of $\fluc_T$, one sees that their effectiveness degrades as $T$ grows (i.e. they do not reflect the centering or scaling associated with the long-time fluctuations).

While the distributional limit in \eqref{eq:eqnchit} does not provide control over the tails of $\fluc_T$ for finite $T$, it does suggest a natural conjecture. For $s$ large, $F_{{\rm GUE}}(-s) \approx e^{-\frac{1}{12}s^3}$ (see Proposition~\ref{TracyWidom} herein, or \cite{TracyWidom94,BBD08,RRV11}). Thus one might expect a similar lower tail bound for $\fluc_T$, at least for large enough $T$. As we prove, this is only half true. In fact, there are two types of decay regimes for the lower tail $\mathbb{P}\big(\fluc_T<-s\big)$: for $T^{2/3}\gg s\gg 0$ the cubic exponent controls the tail decay whereas for $s\gg T^{2/3}$ the tail exponent becomes $5/2$ (and the leading constant in the exponential is $\frac{4}{15\pi}T^{1/3}$ instead of $\frac{1}{12}$ in the first regime).

\medskip

We now state the main result of this paper\footnote{There is also forthcoming work (done independently and in parallel) of \cite{KlD} which probes (non-rigorously and in some cases numerically) some elements of the crossover behavior of the lower tail also using the formula in Proposition \ref{MomentMatch} as a starting point. The lower tail for the half-space problem is also discussed therein. Results from that, as well as from our present work will be briefly summarized in a forthcoming physic letter \cite{Letter}.}.
\bt\label{MainResult}
Let $\fluc_{T}$ denote the centered and scaled KPZ solution with narrow wedge initial data as in \eqref{eq:flucs}. Fix $\epsilon, \delta\in(0,\frac{1}{3})$ and $T_0>0$. Then, there exist $S = S(\epsilon, \delta, T_0)$, $C=C(T_0)>0$, $K_1=K_1(\epsilon,\delta, T_0)>0$ and $K_2=K_2(T_0)>0$ such that for all $s\geq S$ and $T\geq T_0$,
\begin{align}\label{eq:MainUpperBound}
\mathbb{P}\big(\fluc_T\leq -s\big)\leq  e^{-\frac{4(1-C\epsilon) }{15\pi}T^{\frac{1}{3}} s^{5/2}} + e^{-K_1s^{3-\delta}-\epsilon T^{1/3}s} + e^{-\frac{(1-C\epsilon)}{12}s^3},
\end{align}
and
\begin{equation}\label{eq:MainLowerBound}
\mathbb{P}\big(\fluc_T\leq -s\big) \geq e^{-\frac{4(1+C\epsilon)}{15\pi}T^{\frac{1}{3}}s^{5/2}}+  e^{-K_2 s^{3}}.
\end{equation}
\et

We prove this in Section~\ref{sec:Proof}. Note that the right side of \eqref{eq:MainUpperBound} is a sum of three terms. The first dominates the other two when $s\gg T^{\frac{2}{3}}$. In the region $T^{\frac{2}{3}}\gg s\gg 0$, the second and third terms dominate, and when $T\to \infty$, the third dominates the second and recovers the $\tfrac{1}{12}s^3$ tail behavior of the GUE Tracy-Widom distribution. There is a similar interplay between the two terms in \eqref{eq:MainLowerBound}, though in this lower bound we do not recover\footnote{We expect this is just a limitation of our result and would follow from a finer analysis.} the $\tfrac{1}{12}$ constant as $T\to \infty$.

The KPZ equation is believed to be the unique heteroclinic orbit between the Edwards-Wilkinson (i.e. weak coupling) and KPZ (i.e. strong coupling) fixed points \cite{KPZfixed}. Theorem \ref{MainResult} essentially demonstrates how the tail behavior crosses over between neighborhoods of these two fixed points.
The $\frac{5}{2}$ exponent in Theorem \ref{MainResult} corresponds with the short-time lower tail exponent\footnote{Note, however, that the Edwards-Wilkinson equation itself has Gaussian tails with exponent $2$.} while the $3$ exponent is that of the long-time lower tail. This can also be interpreted in terms of large deviations for the KPZ equation -- see Section \ref{secLDP}. While the fluctuations for models which renormalize to the KPZ fixed point should be universal (e.g. GUE Tracy-Widom in this case), the large deviation rate function should vary from model to model. Theorem \ref{MainResult} is, to our knowledge, the first large deviation result for a non-determinantal model in the KPZ universality class -- see Section \ref{secother} for further discussion.


\medskip

\noindent We now briefly explain the three steps in our proof, though to simplify the exposition we will leave off the $\epsilon$ and $\delta$'s which are present in the statement and proof.

\smallskip
\noindent {\bf Step 1:} Our starting point is the KPZ equation one-point formula \cite{Amir11, SS10, CDR10, Dot10}. \cite{BorGor16} reformulated that result as an identity between the Laplace transform of the SHE and the expectation of a specific multiplicative functional of the Airy point process (see Proposition \ref{MomentMatch}). Armed with this, our first deduction is that the large parameter ($u$ in \eqref{eq:NarrowWedgeToAirtyRelation}) asymptotics of the SHE Laplace transform translate into lower tail asymptotics for the KPZ equation. This reduces Theorem \ref{MainResult} to Proposition \ref{MostImportantTheorem} (see the proof of Theorem \ref{MainResult} for details). Using Proposition \ref{MomentMatch} we reduce Proposition \ref{MostImportantTheorem} to Proposition \ref{ppn:ImportantLemma} which studies Airy point process asymptotics and whose proof is the main technical feat of this paper.

\smallskip
\noindent {\bf Step 2:} The proof of Proposition \ref{ppn:ImportantLemma} relies upon three results (Theorems \ref{thm:LRigidityBound}, \ref{UpperTailLemma} and \ref{cor:AiryTails}) about large deviations of the number of Airy points in large intervals and their rigidity around typical locations. Theorems \ref{thm:LRigidityBound} and  \ref{UpperTailLemma} respectively probe the lower and upper large deviation tails for the fluctuations of the number of Airy points in a large interval $[-s,\infty)$. The mean\footnote{The variance grows like $\log(s)$ and the fluctuations satisfy a central limit theorem in this scale \cite{Soshnikov2000}.} number of points grows (Proposition \ref{ppn:ExpVarOfLS}) like $\frac{2}{3\pi}s^{\frac{3}{2}}$ and these theorems probe the probability of finding a different constant than $\frac{2}{3\pi}$. On the lower tail, Theorem \ref{thm:LRigidityBound} shows that the exponential decay power law has exponent $3$, while Theorem \ref{UpperTailLemma} shows that the corresponding upper tail exponent is $\frac{3}{2}$. To our knowledge, such large deviation result are new for the Airy point process (cf. Sections \ref{secrigid} and \ref{secLDP} for further discussion). Theorem \ref{cor:AiryTails} controls the maximum (over the entire Airy point process) deviation of points outside bands around their typical locations. We do not expect this result is nearly as tight Theorems \ref{thm:LRigidityBound} and  \ref{UpperTailLemma}, but it suffices for our purposes. Using these three theorems we can establish control the probabilities of various scenarios for the Airy point process and hence establish precise upper and lower bounds on the expectation value needed to prove Proposition \ref{ppn:ImportantLemma}.

\smallskip
\noindent {\bf Step 3:} The proofs of Theorems \ref{thm:LRigidityBound}, \ref{UpperTailLemma} and \ref{cor:AiryTails} are each rather different. The first two rely on the determinantal structure of the Airy point process (Section \ref{secdet}), while the third uses its relation to the stochastic Airy operator (Section \ref{SAO}). The proof of Theorem \ref{thm:LRigidityBound} is technically the most challenging. Via Markov's inequality, it reduces to a bound on the cumulant generating function for the number of Airy points in the interval $[-s,\infty)$, when the parameter $v$ of the generating function is of order $s^{\frac{3}{2}}$ (see Section \ref{secablowitz}). Theorem \ref{thm:CGFexpansion} relates (via standard determinantal methods) this generating function $F(x;v)$ to the Ablowitz-Segur solution to the Painlev\'e II equation, and then proves the needed decay bound on the generating function using a delicate analysis of an asymptotic formal (given in recent work in \cite{Bothner15} in terms of oscillatory Jacobi elliptic functions) for this solution to Painlev\'e II. The proof Theorem \ref{UpperTailLemma} is considerably simpler. It uses the fact that the number of Airy points in an interval equals (in law) the sum of independent Bernoulli random variables (with parameters related to the eigenvalues of the Airy kernel projected onto the interval). The theorem follows by combining Bennett's concentration inequality on such sums, along with estimates on mean and variance given in Proposition \ref{ppn:ExpVarOfLS}. Theorem \ref{cor:AiryTails} uses the identity in law (Proposition \ref{ppn:connection}) between the Airy point process and the spectrum of the stochastic Airy operator. The typical locations of the Airy points are given by the zeros of the Airy function, and the estimate on uniform deviations from bands around those typical locations can be reduced (through operator manipulations such used in \cite{RRV11,Balint14}) to an exponential tail estimate (proved in Lemma \ref{lm:ControlRegularity}) on of the maximum oscillation of Brownian motion\footnote{The Brownian motion is the driving noise for the stochastic Airy operator.}.

\medskip

The rest of this introduction records the main results (summarized above) which go into our proof of Theorem \ref{MainResult}. Section \ref{seclaplace} provides the key identity relating the Laplace transform of the SHE and the expectation of a multiplicative functional of the Airy point process. Section \ref{secrigid} records the Airy point process large deviation and rigidity estimates that we rely upon. Section \ref{secablowitz} records the precise asymptotics of the Ablowitz-Segur solution of the Painlev\'{e} II equation needed in the proof of Theorem \ref{thm:LRigidityBound}.

\subsection{Laplace transform formula}\label{seclaplace}

The starting point for our study is the exact formula characterizing the one-point distribution of the SHE with delta initial data. This was simultaneously and independently computed in \cite{Amir11, SS10, CDR10, Dot10} (rigorous proof provided in \cite{Amir11}). That formula can, by straight-forward manipulations, be reformulated (Proposition \ref{MomentMatch}) in terms the expectation of a multiplicative functional of the Airy point process (Section \ref{sec:Airy}). This was done in \cite[Theorem~2.2]{BorGor16}, and the resulting formula offers a major benefit since it enables one to bring to bare on the KPZ equation the vast range of tools and understanding developed for the Airy point process. In fact, prior to our present work, it was not clear how to prove directly that the formula in \cite{Amir11, SS10, CDR10, Dot10} defines a probability distribution\footnote{The hard part is to prove that the lower tail probability decays to 0.}. Armed with Proposition \ref{MomentMatch} such a result is immediate, and the lower tail decay becomes  more tractable.

Proposition \ref{MomentMatch} is a special limit  case of a general matching between stochastic vertex models and Macdonald measures in \cite[Corollary 4.4]{BorodinMoments}. In special cases, the Macdonald measures reduce to determinantal Schur measures and hence are analyzable in the spirit of this paper (see \cite{BO16,BBW16,BBCW17} or Sections \ref{sechalfspace}  and \ref{secother} for further discussion).

\bp[Theorem~2.2 of \cite{BorGor16}]\label{MomentMatch}
Let $\mathcal{Z}(T,X)$ be the unique solution to the SHE \eqref{eq:SHEDef} with $\mathcal{Z}(0,X)=\delta_{X=0}$. Denote the ordered  points of the Airy point process (Section~\ref{sec:Airy}) by $\mathbf{a}_1>  \mathbf{a}_2>\ldots $. Then, for any $T,u>0$, we have\footnote{A similar result holds for any $X$ up to multiplying $\mathcal{Z}$ by a Gaussian factor -- see \cite[Proposition 1.4]{Amir11}.}
\begin{align}\label{eq:NarrowWedgeToAirtyRelation}
\mathbb{E}_{\mathrm{SHE}}\Bigg[\exp\bigg(-u \mathcal{Z}(2T,0)\exp\Big(\frac{T}{12}\Big)\bigg)\Bigg] = \mathbb{E}_{\mathrm{Airy}}\Bigg[\prod_{k=1}^{\infty} \frac{1}{1+u \exp\big(T^{\frac{1}{3}}\mathbf{a}_k\big)}\Bigg]
\end{align}
\ep

Setting $u=\exp\big(T^{\frac{1}{3}}s\big)$ and rewriting the above result in terms of $\fluc_T$ from \eqref{eq:eqnchit}, we find
\begin{align}\label{eq:SHEAiryeqn}
\mathbb{E}_{\mathrm{SHE}}\Bigg[\exp\bigg(-\exp\Big(\big(T^{\frac{1}{3}}\big(\fluc_T+s\big)\Big)\bigg)\Bigg] = \mathbb{E}_{\mathrm{Airy}}\Bigg[\prod_{k=1}^{\infty} \frac{1}{1+\exp\big(T^{\frac{1}{3}}(s+\mathbf{a}_k)\big)}\Bigg]
\end{align}
The function $\exp\big(-\exp(x)\big)$ is an approximate version of $\ind{x<0}$ and thus when $s$ is large, the expectation on the left-side of \eqref{eq:SHEAiryeqn} is approximately $\PP\big(\fluc_T+s<0\big)$ which is exactly the tail we are looking to control. Now consider the right-hand side of \eqref{eq:SHEAiryeqn}. If $s+\mathbf{a}_k\gg 0$ then the corresponding term in the product will be exponentially small, whereas if $s+\mathbf{a}_k\ll 0$ then the term will be very close to 1. Thus, the tail decay on the left-hand side is linked with the number of exponentially small terms (and their exponential factors) on the right-hand side.

Typically, the Airy point process is close to the zeros of the Airy function (Proposition \ref{ppn:SandwichResult}), and hence $\mathbf{a}_k \sim -\big(\frac{3\pi}{2}k\big)^{\frac{2}{3}}$ (Proposition \ref{MeanPosition}). Plugging in this estimate readily yields decay like $\exp\big(-\frac{4}{15\pi}T^{\frac{1}{3}} s^{\frac{5}{2}}\big)$. The Airy points may, however, differ from these typical locations. For instance, $\mathbf{a}_1$ (which is GUE Tracy-Widom distributed) may dip below $-s$ in which case the product in the expectation on the right-hand side of \eqref{eq:SHEAiryeqn} becomes very close to 1. The probability of such a drastic dip behaves like $\exp\big(-\frac{1}{12}s^3\big)$. Of course, there are many other scenarios in which the Airy points deviate from their typical locations in less drastic ways, and the contributions of those to the overall expectation need to be controlled and ultimately contribute to the other terms in our bounds in Theorem \ref{MainResult}. We give a brief overview of this in Section \ref{secrigid}. Proposition \ref{MostImportantTheorem} (which follows directly from Proposition \ref{ppn:ImportantLemma}) contains precise statements of the bounds that we prove on the behavior of the right-hand side of \eqref{eq:SHEAiryeqn}.

\subsection{Rigidity bounds for Airy point process}\label{secrigid}

The Airy point process $\mathbf{a}_1>  \mathbf{a}_2>\ldots $ (see Section~\ref{sec:Airy}) is a determinantal point process on the real line introduced by Tracy and Widom \cite{TracyWidom94} as the scaling limit of the edge of the spectrum of the Gaussian unitary ensemble (GUE). \cite{TracyWidom94} found that $F(s):=\PP(\mathbf{a}_1<s)$ can be written in terms of the Hastings-McLeod (HM) solution of Painlev\'{e} II:
\begin{equation}\label{eq:HMSol}
F(s)= \exp\left(-\int_{x}^{\infty}(y-x)u^2_{\mathrm{HM}}(y)dy\right)
\end{equation}
where $u_{\mathrm{HM}}(y)$ is the solution of the Painlev\'{e} II equation (introduced in \cite{P00,P02} -- see also the review \cite{FIKN}) with specific boundary behavior as $x\to \infty$:
\begin{align}
u^{\prime\prime}_{\mathrm{HM}} &= xu_{\mathrm{HM}} +2 u^3_{\mathrm{HM}}, \quad (') = \frac{d}{dx},\label{eq:Painleve}\\
u_{\mathrm{HM}}(x) &= \frac{x^{-\frac{1}{4}}}{2\sqrt{\pi}} e^{-\frac{2}{3}x^{\frac{3}{2}}} (1+ o(1)) , \text{ as }x\to \infty.\label{eq:HMBoundary}
\end{align}
This solution was introduced in \cite{HM} wherein they determined an asymptotic formula for $u_{HM}(x)$ as $x\to -\infty$ (this is called solving the {\it connection problem}). Plugging this into \eqref{eq:HMSol}, allowed \cite{TracyWidom94} to demonstrate that $F(-s)$ decays like $\exp(-\frac{s^3}{12})$. Later, using the nonlinear steepest descent technique pioneered by \cite{DZ92}, \cite{DZ,DIK,BBD08} determined smaller order terms in the asymptotic expansion of $F(-s)$. Similar results have been established for other determinant point processes related to KPZ class models, e.g. \cite{BDJ,BKMM}.

In order to make rigorous the heuristic described in the last section we need to establish sufficiently precise control over the deviations of a large number of Airy points from their typical locations. Controlling deviations of eigenvalues from their typical locations is a central theme in some random matrix universality works (see, for example \cite{ESY,EYY} and subsequent works). We require very precise upper and lower bounds on large deviations than do not seem to be present in the existing literature. In fact, we must ultimately rely upon the Ablowitz-Segur solution of Painlev\'{e} II to establish suitably precise bounds.

Our rigidity bound are established in terms of counting Airy points in intervals. Define
\begin{align}\label{eq:CharacteristicFunction}
\chia: \mathcal{B}(\RR)\to \ZZ_{\geq 0}, \quad \chia(B) := \# \{k: \mathbf{a}_k\in B\}, \quad \forall B\in \mathcal{B}(\RR)
\end{align}
where $\mathcal{B}(\RR)$ denotes the Borel $\sigma$-algebra of $\RR$. The cumulants of $\chia(B)$ have been studied in \cite{Soshnikov2000} when the Borel set $B$ is a semi-infinite interval of the form $[-s,\infty)$ or a finite interval of the form $[-ks, -(k-1)s)$.
Following \cite[Theorem~1]{Soshnikov2000} we can record the following formula for the expectation and the variance of the random variable $\chia(B)$.
\bp\label{ppn:ExpVarOfLS}
Define intervals $\mathfrak{B}_k(s):= [-ks,-(k-1)s)$ for $k\in \ZZ_{>1}$ and $\mathfrak{B}_1(s):= [-s,\infty)$. For any $s>0$,
\begin{align}\label{eq:}
\mathbb{E}_{\mathrm{Airy}}\Big[\chia\big(\mathfrak{B}_1(s)\big)\Big] &= \frac{2}{3\pi}s^{\frac{3}{2}} + \mathfrak{D}_1(s), \\
\mathrm{Var}_{\mathrm{Airy}}\Big(\chia\big(\mathfrak{B}_k(s)\big)\Big) &= \frac{11}{12\pi^2}\log(s) + \mathfrak{D}^{(k)}_2(s),\quad \forall k\geq \ZZ_{>0}
\end{align}
where  $\sup_{s\geq 0}\big|\mathfrak{D}_1(s)\big|<\infty$, and $\sup_{s\geq 0}\big|\mathfrak{D}^{(k)}_2(s)\big|<\infty$ for all $k\in \ZZ_{\geq 1}$.
\ep

These estimates can be used to prove a central limit theory for linear statistics (including the number of particles in large intervals) for the Airy point process (see, e.g. \cite{Soshnikov2000}).
By studying higher order cumulants,  \cite[Theorem 5.2]{DE13} derives a moderate deviation result for $\chia(B_n)$ where $B_n:=[-n,n]$. However, their result does not probe far enough into the tails of the distribution (it is still effectively Gaussian) to be of use in our desired application.

The theorems which we now state effectively show that the deviations of $\chia([-s,\infty))$ have the same exponential order (up to some small correct terms) tail behavior as the deviations of $\mathbf{a}_1$. In other terms, the probability of having far too few or far too many points in a large interval is similar to the probability of having the first point far to the left or to right.

\bt\label{thm:LRigidityBound}
For any $\delta>0$, there exist $s_0=s_0(\delta)>0$ and $K= K(\delta)>0$ such that for all $s\geq s_0$ and $c>0$
\begin{align}\label{eq:ChiL}
\mathbb{P}\Big(\chia\big([-s,\infty)\big) - \mathbb{E}\big[\chia\big([-s,\infty)\big)\big]\leq - cs^{\frac{3}{2}}\Big)\leq \exp\Big(-cs^{3-\delta}\big(1-Ks^{-\frac{4\delta}{15}}\big)\Big).&&&
\end{align}
\et

\bt\label{UpperTailLemma}
Recall $\mathfrak{B}_k(s)$ from Proposition \ref{ppn:ExpVarOfLS}.
Fix any $k\in \ZZ_{\geq 1}$, $c>0$ and $\epsilon\in (0,1)$. Then, there exists $s=s_0(k,\epsilon)$ such that for all $s\geq s_0$,
 \begin{align}\label{eq:ChiD}
\mathbb{P}\Big(\chia\big(\mathfrak{B}_k(s)\big)- \mathbb{E}\big[\chia(\mathfrak{B}_k(s))\big]\geq c s^{\frac{3}{2}}\Big) \leq \exp\Big(- cs^{\frac{3}{2}}\big(\log(cs^{\frac{3}{2}})-(1+\epsilon)\log(\log (s))\big)\Big).\qquad
\end{align}
    \et
Theorems \ref{thm:LRigidityBound} and \ref{UpperTailLemma} are respectively proved in Sections \ref{sec:LRgidity} and \ref{sec:UpperTail}. The proof of Theorem~\ref{thm:LRigidityBound} is based on a connection between the cumulant generating function of $\chia\big([-s,\infty)\big)$ and the Ablowitz-Segur solution of Painlev\'{e} II (Section~\ref{secablowitz}). The proof of Theorem \ref{UpperTailLemma} is simpler, relying on estimates in Proposition \ref{ppn:ExpVarOfLS} along with Bennett's concentration inequality.

In addition to controlling the number of Airy points in large intervals, we require some uniform bound on the distance between the points and their typical locations. Let $\lambda_1<\lambda_2<\cdots$ denote the eigenvalues of the Airy operator (see Section \ref{SAO}). As shown in Proposition \ref{MeanPosition}, $\lambda_n\approx \big(\frac{3\pi}{2}n\big)^{\frac{2}{3}}$. The following result follows directly from combining Proposition \ref{ppn:SandwichResult} with $\beta=2$, and Proposition \ref{ppn:connection}. Proposition \ref{ppn:SandwichResult} is a similar bound for the Airy$_\beta$ point process, and its proof relies on studying the spectrum of the stochastic Airy operator. To control the deviations of that random operator's spectrum, we prove a result (Lemma \ref{lm:ControlRegularity}) which precisely controls the oscillations of Brownian motion. We do not claim that the next rigidity result is optimal and it may be possible to prove similar (or better) results about uniform rigidity of Airy points via other methods, e.g. \cite[Theorem 3.1]{Borg}.

\bt\label{cor:AiryTails}
For $\epsilon\in (0,1)$, let $C^{\mathrm{Ai}}_{\epsilon}$ be the smallest real number such that for all $k\geq 1$
    \begin{align}\label{eq:A-ASandwitch}
    (1-\epsilon) \mathbf{\lambda}_k -C^{\mathrm{Ai}}_\epsilon \leq -\mathbf{a}_k \leq (1+\epsilon)\mathbf{\lambda}_k +C^{\mathrm{Ai}}_\epsilon
    \end{align}
Then, for all $\epsilon,\delta\in (0,1)$ there exist $s_0= s_0(\epsilon, \delta)$, and $\kappa= \kappa(\epsilon,\delta)$ such that for~$s\geq s_0$,
\begin{equation}\label{eq:SupDevTail}
\mathbb{P}\big(C^{\mathrm{Ai}}_{\epsilon}\geq s\big)\leq \kappa\exp\big(-\kappa s^{1-\delta}\big).
\end{equation}
\et





\subsection{Asymptotics of Ablowitz-Segur solution of Painlev\'{e} II}\label{secablowitz}


The proof of Theorem \ref{thm:LRigidityBound} relies on Markov's inequality which shows that for any $v>0$,
\begin{align}\label{eq:LowerTailDemo}
 \mathbb{P}\left(\chia\big([-s,\infty)\big)- \mathbb{E}[\chia\big([-s,\infty)\big)]\leq -cs^{\frac{3}{2}}\right)\leq e^{-cvs^{\frac{3}{2}}+v\mathbb{E}[\chia([-s,\infty)]} F(-s;v)
\end{align}
where
\begin{align}\label{eq:fxv}
F(x;v):=\mathbb{E}\Big[\exp\big(-v\chia\big([x,\infty)\big)\big)\Big].
\end{align}

In \eqref{eq:LowerTailDemo} we choose $v=s^{\frac{3}{2}-\delta}$. In order to extract asymptotics of $F(x;v)$ (see Theorem \ref{thm:CGFexpansion}), we rely on a connection to the Ablowitz-Segur (AS) solution to the Painlev\'{e} II equation.

The Ablowitz-Segur (AS) solution $u_{\mathrm{AS}}(\cdot;\gamma)$ of the Painlev\'{e} II equation is an one parameter family of solutions to \eqref{eq:Painleve} characterized by the following boundary condition
\begin{align}\label{eq:ASBoundary}
     u_{\mathrm{AS}}(x;\gamma) = \sqrt{\gamma}\frac{x^{-\frac{1}{4}}}{2\sqrt{\pi}} e^{- \frac{2}{3}x^{\frac{3}{2}}}\big(1+o(1)\big) \text{ as } x\to \infty.
\end{align}
(Here $o(1)$ means any function which goes to 0 as $x\to \infty$.) For fixed $\gamma\in (0,1)$, \cite{AS76, AS77} solved the connection problem (behavior as $x\to -\infty$). The case $\gamma=1$ is the Hastings-McLeod solution analyzed in \cite{HM}, and the case when $\gamma>1$ was subsequently studied in \cite{K92}.

\bt\label{thm:CGFexpansion}
For $K^{\mathrm{Ai}}$ the Airy point process correlation kernel (Section~\ref{sec:Airy}) and $\gamma~=~1~-~e^{-v}$,
\begin{align}\label{eq:ASsol}
F(x;v)= \mathrm{det}\big(I- \gamma K^{\mathrm{Ai}}\big)_{L^2([x, \infty))}= \exp\Big(- \int^{\infty}_{x} (y-x)u^2_{\mathrm{AS}}(y;\gamma)dy\Big).
\end{align}
Fix any $\delta\in (0, \frac{2}{3})$ and set $v= s^{\frac{3}{2}-\delta}$. Then, as $s$ goes to $\infty$,
 \begin{align}\label{eq:F_2}
\log F(-s;v) \leq  -\frac{2}{3\pi}vs^{\frac{3}{2}}  + \mathcal{O}(s^{3-\frac{19\delta}{15}}).
\end{align}
\et

The first part of this result, \eqref{eq:ASsol}, contains two equalities. The first follows from general theory relating multiplicative functions of determinant point processes to Fredholm determinants (see \cite[Section~3.4]{AGZ10} for background on Fredholm determinants): For a determinantal point process $X$ with state space $\mathcal{X}$ and correlation kernel $K^X$, and a function $\phi:\mathcal{X}\to \CC$,
\begin{align}\label{eq:multfunct}
\EE\Big[\prod_{x\in X} \phi(x)\Big] = \det\big(1-(1-\phi)K^X\big)_{L^2(\mathcal{X})}.
\end{align}
This identity requires $(1-\phi)K^X$ to be trace-class (see \cite{BorDet} for more details). The second equality in \eqref{eq:ASsol} relies on the integrable structure of the Airy kernel \cite[Section~1.C]{TracyWidom94}.

Proving the second part of the theorem, namely \eqref{eq:F_2}, requires a close analysis of the AS solution to Painlev\'{e} II, as is provided in Section \ref{App:ASSol}.

\medskip

The AS solution has received some attention recently in \cite{BohigasCarvalhoPato,BB2017} due to the fact that $\gamma K^{\mathrm{Ai}}$ represents the kernel for a {\it thinned} version of the Airy point process -- each particle is removed with probability $1-\gamma$. This thinning represents one way to achieve a crossover between the GUE Tracy-Widom distribution and more classical extreme value statistics. The study of positive temperature free-Fermions in Section \ref{secferm} represents  another such mechanism.

\cite{AS76,AS77} solved the connection problem for the AS solution for $\gamma\in (0,1)$ fixed. For our application, $\gamma$ (or equivalently $v$) fixed would only yield an exponent of $s^{\frac{3}{2}-\delta}$ in Theorem \ref{thm:LRigidityBound} (not the desired $s^{3-\delta}$). Recently, utilizing Riemann-Hilbert steepest descent, \cite{Bothner15} computed the asymptotic form of the AS solution $u_{\mathrm{AS}}(x;\gamma)$ as $x\to -\infty$ for a more general range of $\gamma$. The formulas are written in terms of Jacobi elliptic theta functions and take different forms depending on the values of $\gamma$. In particular, setting $\tau:= -\frac{1}{(-x)^{3/2}}\log(1-\gamma)$, \cite{Bothner15} computes asymptotic formulas in three different ranges of parameters:
(a) $\tau \in \big(0,  (-x)^{-\delta}\big]$; (b)  $\tau\in \big(0, \tfrac{2}{3}\sqrt{2}-\eta\big]$; (c) $\tau \in \big(\tfrac{2}{3}\sqrt{2}-\aleph\tfrac{\log (-x)^{\frac{3}{2}}}{(-x)^{\frac{3}{2}}}, \infty\big)$. Here $\delta, \eta>0$ are arbitrary small numbers and $\aleph\in \big(-\infty, \frac{7}{6}\big]$. For $\tau\in\big(0,\tfrac{2}{3}\sqrt{2}\big)$ the resulting asymptotic form of $u_{\mathrm{AS}}(x;\gamma)$ as $x\to -\infty$ is pseudoperiodic, thus making it rather challenging to compute the integral in the exponential in \eqref{eq:ASsol} (as necessary to recover asymptotics for $F(x;v)$). As $\tau$ approaches $0$ and $\tfrac{2}{3}\sqrt{2}$ the oscillations die out, though due to different mechanisms in each case.

\cite{Bothner15,Bothner16} managed to translate his asymptotic result for  $u_{\mathrm{AS}}$ into a corresponding result for $F$ only in the (c) region\footnote{\cite{Bothner15} achieved this for $\tau>\tfrac{2}{3}\sqrt{2}$ based on the lack of oscillations in $u_{\mathrm{AS}}$ for such $\tau$, and \cite{Bothner16} provided an extension to the full region (c) (and slightly beyond).}.
For region $(a)$, \cite{Bothner15} demonstrated a simplified form of $u_{\mathrm{AS}}(x;\gamma(x))$ for $\tau \in \big(0, (-x)^{-\delta}\big)$ for any fixed $\delta>0$. However, this simplified form still retains its oscillatory nature which is one of the difficulties in getting a full expansion for $F(-s;1-e^{-s^{3/2-\eta}})$. Recently,  \cite{BB2017} showed that
for any $0<\epsilon<\frac{1}{2}$, there exist constants $s_0=s_0(\epsilon)$ and $c^{\prime}_j=c^{\prime}_j(\epsilon)$ for $j=1,2$ so that for $s\geq s_0$ and $0\leq v= -\log(1-\gamma)<s^{\frac{1}{2}-\epsilon}$,
\begin{equation}\label{eq:BoBu}
\log F(-s;v) = -\frac{2v}{3\pi}s^{\frac{3}{2}}+ \frac{v^2}{4\pi^2}\log(8s^{\frac{3}{2}})+ \log\bigg(G\Big(1+\frac{\mathbf{i }v}{2\pi}\Big)G\Big(1- \frac{\mathbf{i}v}{2\pi}\Big)\bigg)+r(s,v).
\end{equation}
Here $G(x)$ is the Barnes G-function and $|r(s,v)|\leq c^{\prime}_1\frac{v^3}{s^{\frac{3}{2}}} + c^{\prime}_2\frac{v}{s}$ for all $s\geq s_0, 0\leq v\leq s^{\frac{1}{2}-\epsilon}$.

Since \eqref{eq:BoBu} gives the full expansion of $\log F(-s;s^{\frac{3}{2}-\delta})$ only when $\delta>\frac{2}{3}$, plugging it into the right side of \eqref{eq:LowerTailDemo} only yields a leading term (in the upper bound of the lower tail probability of $\chia([-s,\infty))$) like $\exp(-cs^{2-})$. However, Theorem~\ref{thm:LRigidityBound} asks that the upper bound is like $\exp(-cs^{3-})$. In Section \ref{App:ASSol} we demonstrate how we can work with $\delta$ close to $0$. Presently we cannot justify a full expansion of $F(-s;v)$ in Theorem~\ref{thm:CGFexpansion} like that of \eqref{eq:BoBu}. However, the weaker result in Theorem~\ref{thm:CGFexpansion} suffices for our present needs.

\subsection*{Outline}
Rest of this paper is organized as follows. Section \ref{sec.con} includes a brief discussion of how our results and methods connect to other problems and may be extended in other directions. Section \ref{sec:Proof} reduces the proof of our main result (Theorem~\ref{MainResult}) to a result (Proposition~\ref{MostImportantTheorem}) for a cumulant generating function. Proposition \ref{MostImportantTheorem} is subsequently proved in Section \ref{sec:MostImpTheo} by reducing it to a result (Proposition \ref{ppn:ImportantLemma}) about the Airy point process. The rest of Section \ref{sec:Airy} develops and proves various properties about the Airy point process, including the key rigidity estimates stated in the introduction as Theorems \ref{thm:LRigidityBound}, \ref{UpperTailLemma} and \ref{cor:AiryTails}. Proposition \ref{ppn:ImportantLemma} is proved in Section \ref{sec:propproof}. Finally, Section \ref{App:ASSol} contains a discussion on asymptotics of the Ablowitz-Segur solution to Painlev\'{e} II and a proof of Theorem \ref{thm:CGFexpansion}, stated earlier in the introduction.
\subsection*{Acknowledgements}
I.C. and P.G. wish to thank A. Aggarwal, J. Baik, A. Borodin, P. Bourgade, T. Bothner, P. Deift, V. Gorin, T. Halpin-Healy, A. Krajenbrink, P. Le Doussal, K. Liechty, B. Meerson, H. Spohn, B. Virag and O. Zeitouni for discussions and comments related to this project.
I.C. and P.G. initiated this project during the 2017 Park City Mathematics Institute, funded in part by NSF grant DMS:1441467.
I.C. was partially funded by NSF grant DMS:1664650 and the Packard Foundation through a Packard Fellowship for Science and Engineering.

\section{Connections and extensions}\label{sec.con}

We discuss various applications and extensions of our results and methods. Section \ref{secintegro} describes the relationship between our analysis and an inverse-scattering problem generalizing the Painlev\'{e} II equation.  Section \ref{secferm} explains how our results relate to the lower-tail decay for positive temperature free-Fermions. Section \ref{secLDP} discusses extending our analysis to study the KPZ equation large deviation rate function, as well as relates our work to recent physics literature.
Sections \ref{secuppertail}, \ref{secotherinitialdata}, \ref{sechalfspace} and \ref{secother} touch upon extensions of our methods and results to (respectively) the KPZ equation upper-tail decay, general initial data, half-space geometry, and certain discretizations of the KPZ equation like ASEP or the stochastic six vertex model.

\subsection{An integro-differential generalization of Painlev\'{e} II}\label{secintegro}

Using the explicit form of the Airy kernel and the fact \eqref{eq:multfunct} that expectations of multiplicative functions of determinant point processes can be written as Fredholm determinants we can rewrite the equality in Proposition \ref{MomentMatch} (actually \eqref{eq:SHEAiryeqn}) as
\begin{align}\label{eq:SHEdet}
\mathbb{E}_{\mathrm{SHE}}\Bigg[\exp\bigg(-\exp\Big(\big(T^{\frac{1}{3}}\big(\fluc_T+s\big)\Big)\bigg)\Bigg] = \det(I-K)_{L^2(s,\infty)}=: Q(s)
\end{align}
where $K$ is the Airy kernel deformed by a Fermi-factor:
\begin{align}\label{eq:Kcross}
K(x,x') = \int_{-\infty}^{\infty} dr \sigma(r) \Ai(x+r)\Ai(x'+r),\quad \textrm{with} \quad\sigma(r)= \frac{1}{1+e^{-T^{\frac{1}{3}} r}}.
\end{align}

It was proved in \cite[Section 5.2]{Amir11} (following \cite{TracyWidom02}) that for any choice of $\sigma(r)$ (which is smooth except at a finite number of points at which it has bounded jumps, and which approaches $0$ at $-\infty$ and $1$ at $+\infty$ exponentially fast) and the resulting $Q(s)$ satisfies
\begin{align}
\frac{d^2}{ds^2} \log Q(s)  &= \int_{-\infty}^{\infty} dr\sigma'(r) q_r^2(s),\\
Q(s)  &= \exp\Big(-\int_s^{\infty}dx (x-s) \int_{-\infty}^{\infty} dr\sigma'(r) q_r^2(x)\Big),
\end{align}
where $q_r(s)$ solves the following integro-differential generalization of Painlev\'{e} II:
\begin{equation}\label{eq:piigen}
\frac{d^2}{ds^2} q_r(s) = \Big(s+r + 2 \int_{-\infty}^{\infty}dr'\sigma'(r') q_{r'}^2(s)\Big) q_r(s),\qquad \textrm{with } q_r(s) \sim \Ai(r+s)\textrm{ as }s\to +\infty.
\end{equation}
If $\sigma(r) = \ind{r\geq 0}$ then the above equation recovers the Hastings-McLeod solution to Painlev\'{e} II. The derivation of the above result in \cite[Section 5.2]{Amir11} came from an attempt to directly study the lower tail for the KPZ equation\footnote{Due to the complexity of this equation, \cite{Amir11} was unable to even show that the lower tail decays to zero  and resorted to a more indirect route via the results of \cite{CM91} -- see Section \ref{secLDP}, however, for mention of some recent non-rigorous physics attempts at doing this.}. We may reverse the direction of inference and try to use our methods for studying the KPZ tail to deduce results for the solution to \eqref{eq:piigen}.

The connection problem for \eqref{eq:piigen} asks how the Airy behavior as $s\to \infty$ propagates through as  $s\to -\infty$. This problem also falls under the realm of inverse scattering on the line \cite{DeiftTrubowitz, BDT88}. For the Hastings-McLeod solution of the Painlev\'{e} II equation, this problem has been resolved to a great level of detail using the steepest descent method for an associated $2\times 2$ Riemann-Hilbert problem \cite{DZ92,DZ,DIK,BBD08, BB2017,FIKN}.

For a general choice of $\sigma(r)$, the kernel $K$ may be rewritten as
$$
K(x,x') = \int_{-\infty}^{\infty} dr \sigma'(r) \frac{\Ai(x+r)\Ai'(x'+r)-\Ai'(x+r)\Ai(x'+r)}{x-x'}
$$
and hence takes the form of an {\it integrable integral operator}. As shown in \cite{IIKS}, the associated $Q(s)$ can be written in terms of an operator valued\footnote{When $\sigma'(r)$ is a sum of $N$ delta functions, the resulting Riemann-Hilbert problem is $2N\times 2N$ dimensional.} Riemann-Hilbert problem. The analysis of such problems is considerably more involved than in the finite dimensional (namely $2\times 2$) matrix setting (cf. \cite{ItsKoz1,ItsKoz2} for some recent advances in this direction).

The approach developed in this present paper may offer an alternative to studying the operator valued Riemann-Hilbert problem. In our analysis there is nothing particularly special about the choice of $\sigma(r)$ (which translates into the choice of multiplicative functional). For another $\sigma(r)$ we could just as well similarly derive asymptotics for $Q(s)$. Turning this into a solution to the connection problem in \eqref{eq:piigen} may still be a challenge. Should this work, the study of the operator valued Riemann-Hilbert problem would be reduced to the study of the $2\times 2$ matrix problem associated with the Hastings-McLeod and Ablowitz-Segur solutions. We do not pursue this idea further in the present text and leave it for further investigation.

\subsection{Positive temperature free-Fermions}\label{secferm}

Positive temperature free-Fermions and the equivalent MNS matrix model have recently been studied in \cite{DDMS,LW} (and earlier in \cite{Joh} in a grand-canonical form). These ensembles are indexed by an inverse temperature $\beta$. When $\beta\to \infty$ this recovers the Gaussian Unitary Ensemble. \cite{Joh,DDMS,LW} consider taking the number of Fermions (or matrix dimension) $N\to \infty$.
When $\beta$ is fixed, the distribution of the rightmost Fermion converges to the GUE Tracy-Widom distribution (see \cite[Theorem 2(a)]{LW}); when $\beta$ tends to 0 sufficiently fast relative to $N$ going to infinity, the rightmost Fermion converges to a Gumbel distribution; and when $\beta$ tends to 0 and $N$ tend to infinity in a critical manner, there is a crossover between the GUE Tracy-Widom and Gumbel distribution. The limit of the correlation kernel for Fermion point process at the edge converges under this critical scaling to the Fermi-factor deformation of the Airy kernel given in \eqref{eq:Kcross}. As such the Fredholm determinant in \eqref{eq:SHEdet} gives the probability that the right-most limiting Fermion is located to the left of $s$, and Proposition \ref{MostImportantTheorem} provides the lower tail probability decay of that distribution.

\subsection{Large deviation rate function}\label{secLDP}

Theorem \ref{MainResult} shows that there is a crossover between two types of tail decay which occurs when $s$ is of order $T^{2/3}$. This can be understood in terms of large deviations. For $z\leq 0$ let
\begin{equation}\label{eq.phi}
\Phi_-(z) = -\lim_{T\to \infty} T^{-2} \log\bigg( \PP\Big(\mathcal{H}(2T,0) + \frac{T}{12} \leq z T\Big)\bigg).
\end{equation}
The existence of the above limit has not, to our knowledge, been proved\footnote{\cite{Sly} has an approach to proving the existence of such rate functions for first and directed last passage percolation. Whether this approach lifts to positive temperature models like KPZ remains to be seen.}.

In terms of $\Phi_{-}$, Theorem \ref{MainResult} shows that $\Phi_{-}(z)\approx \frac{1}{12}(-z)^3$ for $z$ near 0 and $\Phi_{-}(z) \approx \frac{4}{15\pi} (-z)^{5/2}$ for $z$ near $-\infty$. In order for a large deviation principle to subsume Theorem \ref{MainResult}, the existence of the limit in \eqref{eq.phi} would need to be uniform in $z$ as $T\to \infty$.

In the physics literature, the crossover between the exponents $\frac{5}{2}$ to $3$ seems to have been first predicted via weak noise theory\footnote{Weak noise theory (WNT), sometimes also called `optimal fluctuation theory' studies the large deviations of the noise necessary to produce a given space-time trajectory of the KPZ equation (or more general systems). It is a valid method only under `weak coupling' or when there is an exceedingly small parameter in front of the noise term. In many instances, this approach is only valid for short times (when the noise is, through rescaling, effectively weak). However, for the KPZ equation it seems that it remains valid for longer times, if one probes deep enough into the tail. WNT has a long and rich history within physics dating back to the 1960s in condensed matter physics \cite{HL66,ZL66,Lif67} and was introduced into the study of the noisy Burgers equation by Fogedby in the late 90s \cite{Fog98}. It also goes under names such as the `instanton method' in turbulence, `macroscopic fluctuation theory' in lattice gases \cite{BDGJL}, and `WKB method' in reaction-diffusion systems (see \cite{MS17} for a more extensive history). Within mathematics, the WNT for diffusions goes under the name Fredilin-Wentzell theory. For field valued / infinite dimensional diffusion processes \cite{BDM08} and for certain non-linear stochastic PDEs  \cite{HW17,CD16}, it has recently received some rigorous treatment. WNT alone does not provide the $\frac{5}{2}$ exponent or associated prefactor. Once the large deviations for the sample path (e.g. evolution of the KPZ equation) is determined, one still needs to solve a Hamiltonian variational problem to figure out the most likely trajectory among all those which achieve a given one-point large deviation. In the physics literature, \cite{KK07,KK09,MKV16a} worked through this calculation for KPZ with flat initial data and predicted the $\frac{5}{2}$ exponent along with a prefactor of $\frac{8}{15\pi}$. \cite{MKV16b} worked with parabolic initial data (which interpolates between flat and narrow wedge) and predicted that the prefactor becomes $\frac{4}{15\pi}$ in the narrow wedge limit. These short-time predictions have been confirmed through exact formulas in physics works such as \cite{DMRS16, KLExact}.}
by \cite{KK07} in the context of directed polymers, and quite recently by \cite{MKV16a} in the context of the KPZ equation.
Recently, this crossover has been studied via analysis of the integro-differential equation discussed in Section~\ref{secintegro}. \cite{PSG16} performed a rough (non-rigorous) analysis of the equation and predicted the existence of a LDP with speed $T^2$ and cubic behavior for small $z$. However, their analysis missed the behavior of $\Phi_{-}(z)$ for $z\ll 0$ and hence did not predict that the $\frac{5}{2}$ exponent remains for long time. Via non-rigorous WKB approximation analysis, \cite{SMP17} predicted not only that the $\frac{5}{2}$ to $3$ crossover should hold for all times sufficiently large, but also predicted a formula for the large deviation rate function $\Phi_{-}(z)$ from \eqref{eq.phi}. The \cite{SMP17} prediction
\begin{equation}
\Phi_-(z)= \frac{4}{15\pi^6}(1-\pi^2 z)^{5/2} - \frac{4}{15 \pi^6} + \frac{2}{3\pi^4}z -\frac{1}{2\pi^2}z^2 \label{eqMeerson}
\end{equation}
indeed formula recovers the desired small and large $z$ asymptotics (see Section \ref{secLDP}). Quite recently, \cite{HLMRS} has performed high-precision simulations (via methods of importance sampling and parallel computation) which numerically confirm the $5/2$ exponent for short and moderate values of time. The cubic exponent is harder to access numerically, though there seems to be some convergence towards that exponent.

We now explain heuristically how our present work could be extended to prove a formula for $\Phi_{-}(z)$. The core challenge  is that there is no proved large deviation theory for the empirical density of the Airy point process (such as done for the GUE point process in \cite{BaG97} -- see also \cite{LS17} and references therein). Since there are infinitely many points in the Airy point process, one cannot naively apply the  Coulumb-gas / electrostatics approach to formulate a large deviation principle. We leave this challenge to future work.

In light of \eqref{eq:NarrowWedgeToAirtyRelation} and the argument used to prove Theorem \ref{MainResult}, $\Phi_{-}(z)$ should be given by
$$
\Phi_{-}(z)=\lim_{T\to \infty} \frac{1}{T^2} \log \mathbb{E}\bigg[\exp\Big(-\sum_{i=1}^\infty \varphi_{T,-zT^{2/3}}(\mathbf{a}_i )\Big)\bigg],
$$
where the $\mathbf{a}_i$ are the Airy point process, and $\varphi_{t,s}(a) := \log\big(1+\exp\big(T^{1/3}(a+s)\big)\big)$.
For large $T$, $\varphi_{T,-zT^{2/3}}(\mathbf{a}_i )\approx T(a-z)_+$ (where $(\cdot)_+ := \max(\cdot,0)$). Letting $\mu_{T}(\cdot) = T^{-1} \sum_{i\geq 1} \delta_{\mathbf{a}_i T^{-2/3}}(\cdot)$ denote the scaled empirical Airy point process measure,
$$
\Phi_{-}(z)=\lim_{T\to \infty} \frac{1}{T^2} \log \mathbb{E}\bigg[\exp\Big(-T^2 \int_{\mathbb{R}}\mathrm{d}a \mu_T(a)(a-z)_+ \Big)\bigg].
$$

Now, assume\footnote{This is where things become quite heuristic and non-rigorous.} that for a suitable class of functions $\mu$, the empirical measure $\mu_T$ satisfies $\mathbb{P}(\mu_T \approx \mu)\approx \exp\big(-T^{-2} I(\mu)\big)$ for a rate functional $I$. Then, we would expect that
\begin{equation}\label{eqImu}
\Phi_{-}(z)=\min_{\mu} \bigg(\int_{\mathbb{R}}\mathrm{d}a\, \mu(a)(a-z)_{+} +I(\mu)\bigg)
\end{equation}
where the minimum is over the class of functions upon which $I$ is finite.

Assuming \eqref{eqImu}, we can derive upper bounds on $\Phi_{-}$. For instance $I(\mu)$ should be minimized and equal to 0 for the limiting density\footnote{This can be calculated, for instance, by taking the trace of the Airy kernel.} of the Airy point process $\mu_{*}(a) = \pi^{-1} \sqrt{-a}\mathbf{1}_{a\leq 0}$. Plugging this choice into \eqref{eqImu} and evaluating the integral gives $\Phi_{-}(z) \leq \frac{4}{15\pi}(-z)^{\frac{5}{2}}$. On the other hand, consider the limiting density of the Airy point process conditioned on $\mathbf{a}_1\leq zT^{2/3}$ (after the scaling discussed above). Since that density will be supported strictly on $(-\infty,z]$, the integral in \eqref{eqImu} will be zero. For that density, $I(\mu)=\frac{(-z)^3}{12}$, as can be determined by the known large deviations for $\mathbf{a}_1$ in Proposition \ref{TracyWidom}. Thus we fine that $\Phi_{-}(z) \leq \frac{(-z)^3}{12}$.

\begin{figure}[ht]
\includegraphics[width=.6\linewidth]{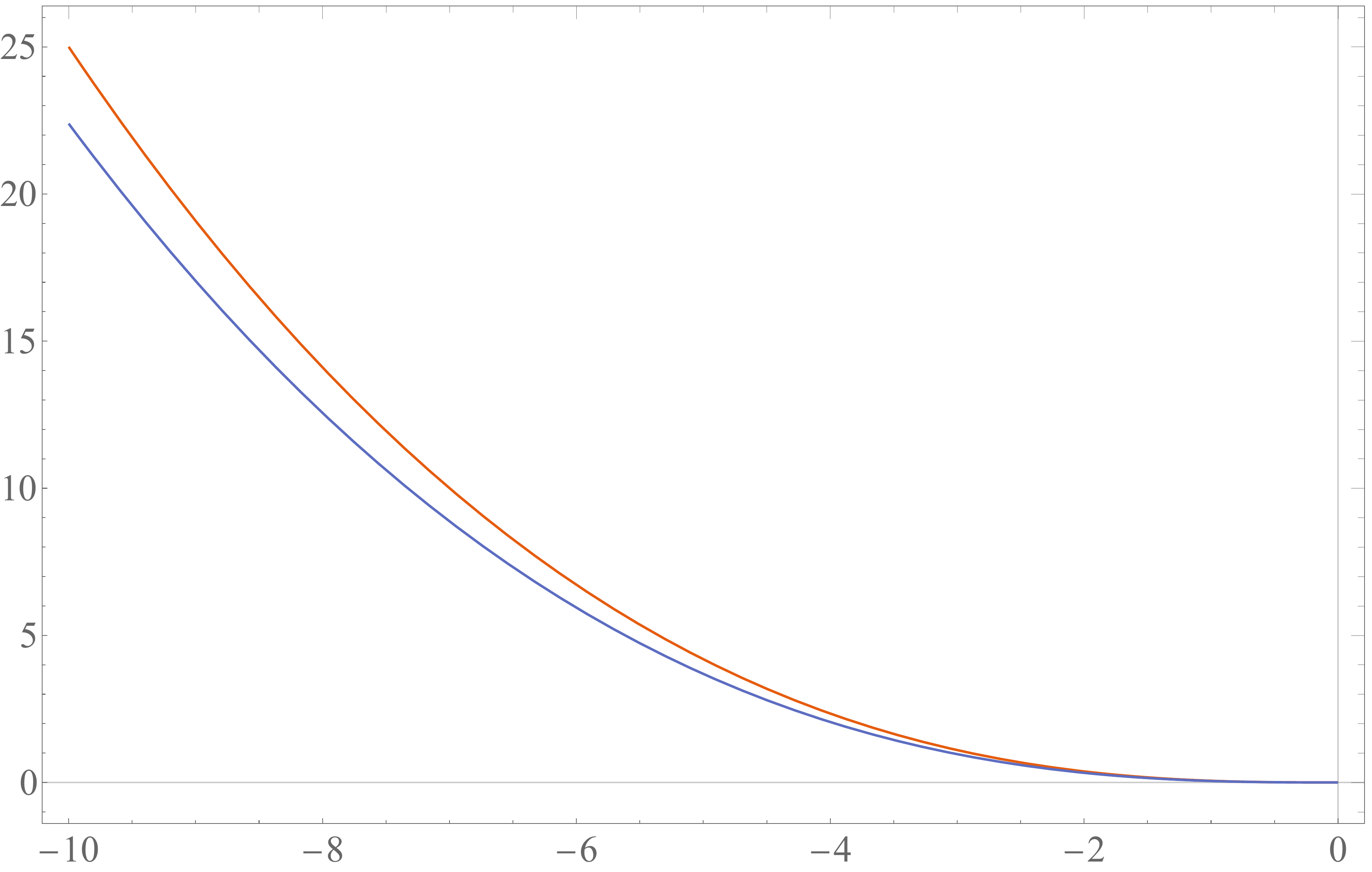}
\caption{The rate function $\Phi_{-}(z)$ from (\ref{eqMeerson}) (bottom curve) compared to the upper-bound rate function $\tilde{\Phi}_{-}(z)$ from (\ref{equpperbd}) (top curve) for $z\in [-10,0]$. Notice that $\Phi_{-}(z)\leq \tilde{\Phi}_{-}(z)$. The ratio $\tilde{\Phi}_{-}(z)/\Phi_{-}(z)$ stays bounded by 1.15 as $z$ varies and tends to 1 as $z\to 0$ or $z\to -\infty$.}
\label{Fig_ratio}
\end{figure}

Without knowing the full rate function $I(\mu)$, we can still improve our upper bound on $\Phi_{-}(z)$ by considering the effect of conditioning on $\big\{\mathbf{a}_1\leq rT^{2/3}\big\}$ for various $r\in [z,0]$. As above, the rate function cost of such a conditioning is $\frac{(-r)^3}{12}$. We should determine the limiting density of the conditioned point process. Non-rigorously, this can be extract by taking a suitable edge limit of physics result contained in \cite{dean2006large,dean2008extreme} for the limit shape of the GUE ensemble under a similar conditioning\footnote{Conditioning the Airy point process on events like $\big\{\mathbf{a}_1\leq rT^{2/3}\big\}$ result in a new determinantal point process whose kernel is modified by the inclusion of a resolvent (see \cite{BorDet,Buf}). An analysis of the trace of this kernel should also reveal the formula for $\mu_{*}^r$.}. Calling $\mu_{*}^r$ the conditional limit density, we find that
$\mu_{*}^r(a) = \frac{r-2a}{2\pi \sqrt{r-a}} \mathbf{1}_{a\leq r}$.
Note that as $r\to 0$, this recovers $\mu_*$. Using $\mu_{*}^r$ we find that
$$
\Phi_{-}(z)\leq \min_{z\leq r\leq 0} \bigg(\int_{\mathbb{R}}\mathrm{d}a\, \mu_{*}^r(a)(a-z)_{+} +\frac{(-r)^3}{12}\bigg).
$$
The argument above is minimized at $r= 4 \pi^{-2} (2-Z)$ with $Z=\sqrt{4-z \pi^2}$ which yields
\begin{align}\label{equpperbd}
\Phi_{-}(z)\leq \tilde{\Phi}_{-}(z):=\frac{2}{15\pi^6} \Big(40(-2+Z)^3 + 2(8-z\pi^2 -4Z)^{3/2}\big(-z\pi^2 +6(-2+Z)\big)\Big).
\end{align}
The two expressions agree in the limits $z\to 0$ and $z\to -\infty$. See Figure \ref{Fig_ratio} for a comparison of $\Phi_{-}$ and $\tilde{\Phi}_{-}$ for intermediate values of $z$. Numerically it is clear that $\Phi_{-}(z)\leq \tilde{\Phi}_{-}(z)$ as desired. A proof the formula for $\Phi_{-}$ will require a better understanding of $I(\mu)$.

\subsection{Upper tail}\label{secuppertail}



Unlike for the lower tail, the upper tail probability $\mathbb{P}(\Upsilon_T> s)$ can be studied via Fredholm determinants \cite[Proposition 10]{CJ13}. The large deviation rate should be $T$ (instead of $T^2$ for the lower tail) and it is predicted in  \cite{PSG16, SMP17} that the rate function is $\frac{4}{3} s^{\frac{3}{2}}$. We leave for future work the problem of proving this via the methods of this paper.
%
%

\subsection{General initial data}\label{secotherinitialdata}


\cite{CorHam16} introduced a method (based on the KPZ line ensemble Gibbs property) to extend tail probabilities for the narrow wedge initial data KPZ equation to corresponding results for quite general initial data. In \cite[Theorem~13]{CorHam16}, the inputs came from \cite{Amir11} and \cite{GMF14} and were far from optimal. In future work we plan to employ our newly proved tight tail bounds from Theorem~\ref{MainResult} to try to derive similar results for general initial data.
%

\subsection{Half-space KPZ}\label{sechalfspace}

The $(1+1)$-dimensional SHE $\mathcal{Z}^{\mathrm{hs}}(T,X)$ in the half space $\RR_{+}$ with delta initial data at the origin is uniquely defined (see \cite{CS16}) by the SPDE in \eqref{eq:SHEDef} and the Robin boundary condition (parametrized by $A\in \RR$) which is formally given as
$
\partial_X \mathcal{Z}^{\mathrm{hs}}(T,X)\big|_{X=0} = A \mathcal{Z}^{\mathrm{hs}}(T,0).
$
for all $T\geq 0$. 
The above half space SHE/KPZ equation has been recently studied in \cite{CS16, SParekh17} where it arises as the scaling limit of a corresponding ASEP. In the spirit of Proposition~\ref{MomentMatch}, \cite{BBCW17} computed a Laplace transform formula for the half-space SHE in terms of the (Pfaffian) GOE point process. We expect that using that as a starting point, our general methods will extend to yield the lower tail the half-space KPZ equation.

\subsection{Other integrable probabilistic systems}\label{secother}

Integrable probabilistic systems in the KPZ universality class \cite{Corwin12} fall into two classes -- determinantal (i.e., free Fermion) or non-determinantal. For determinantal models like the longest increasing subsequence, polynuclear growth model, directed last passage percolation with geometric (or, exponential) weights and the totally asymmetric simple exclusion process (TASEP) there have been a number of works, such as \cite{BDJ,BDMMZ01}, which have obtained optimal lower tail\footnote{For TASEP, the lower tail corresponds to the upper tail for the current of particles to pass the origin.} estimates via analysis of $2\times 2$ Riemann-Hilbert problems (often related to Painlev\'e equations). So far, our present work on the KPZ equation provides the only lower tail bounds for non-determinantal models\footnote{Of course, our analysis ultimately reduces to studying determinantal point processes.}.

Besides studying one-point lower tail decay and large deviations, there is much interest in understanding the large deviations of the entire space-time trajectory. For TASEP a recent attempt at this has been made in \cite{tsai}. The rate is still $N^2$, though the rate function is only bounded above and below in \cite{tsai}.
The stochastic six vertex model \cite{BCG2016} is a discrete time analogs of (T)ASEP. There has been significant efforts (summarized, for instance, in \cite{resh}) to study large deviations and surface tensions for the six vertex model. So far, the only rigorous results (i.e., large deviations for limit shapes) are for determinantal models such as uniform Aztec diamond or rhombus tilings (see, for example, \cite{kos,ko}).

Using the methods considered in this paper, we should be able to access tail / large deviation type results for a few other non-determinantal models. The starting point for our result is the identity in  Proposition~\ref{MomentMatch} which matches the SHE Laplace transform with a multiplicative function of the Airy point process. Similar formulas exist for the asymmetric simple exclusion process (ASEP) \cite[Theorem~1.1]{BO16}, stochastic six vertex model \cite[Corollary~4.4]{BorodinMoments}, $q$-TASEP \cite[Proposition~6.1]{DP17}. The methods of this paper should extend to these other models though will likely involve some new analysis (such as of $q$-Laplace transforms and the associated variants of Painlev\'e which arise for these different models).

 \section{Proof of the main result}\label{sec:Proof}

Recall $\fluc_T$ from \eqref{eq:eqnchit}. For large enough $s$, $\exp(-\exp(T^{\frac{1}{3}}(\fluc_T+s)))$ is approximately equal to $\mathbbm{1}(\fluc_T\leq -s)$. Motivated by this heuristic, we prove Theorem~\ref{MainResult} by estimating the Laplace transform formula $\mathbb{E}[\exp(-\exp(T^{\frac{1}{3}}(\fluc_T+s)))]$. We first state in Proposition~\ref{MostImportantTheorem} matching upper and lower bounds on the Laplace transform formula. Then, using Proposition~\ref{MostImportantTheorem}, we finish the proof of Theorem~\ref{MainResult} in Section~\ref{sec:MainResult}.

\bp\label{MostImportantTheorem}
  Fix $\epsilon, \delta\in (0,\frac{1}{3})$ and $T_0>0$. Then, there exist $s_0=s_0(\epsilon,\delta, T_0)$, $C=C(T_0)>0$, $K_1=K_1(\epsilon,\delta, T_0)>0$ and $K_2=K_2(T_0)>0$ such that for all $s\geq s_0$, one has
  \begin{align}\label{eq:UpperBoundOnInfiniteProduct}
 \mathbb{E}\Big[\exp\Big(-\exp\big(T^{\frac{1}{3}}(\fluc_{T}+s)\big)\Big)\Big] \leq e^{-\frac{4(1-C\epsilon) }{15\pi}T^{\frac{1}{3}} s^{5/2}} + e^{-K_1s^{3-\delta}-\epsilon T^{1/3}s} + e^{-\frac{(1-C\epsilon)}{12}s^3}
  \end{align}
  and
  \begin{align}\label{eq:LowerBoundOnInfiniteProduct}
  \mathbb{E}\Big[\exp\Big(-\exp\big(T^{\frac{1}{3}}(\fluc_{T}+s)\big)\Big)\Big]\geq  e^{-\frac{4(1+C\epsilon)}{15\pi}T^{\frac{1}{3}}s^{5/2}}+  e^{-K_2 s^{3}}.
  \end{align}
\ep
We postpone the proof of Proposition~\ref{MostImportantTheorem} to Section~\ref{sec:MostImpTheo}.


 \subsection{Proof of Theorem~\ref{MainResult}}\label{sec:MainResult}
We show that \eqref{eq:UpperBoundOnInfiniteProduct} (resp. \eqref{eq:LowerBoundOnInfiniteProduct}) implies  \eqref{eq:MainUpperBound} (resp. \eqref{eq:MainLowerBound}) of Theorem~\ref{MainResult}.

Let us first show that \eqref{eq:UpperBoundOnInfiniteProduct}$\Rightarrow$\eqref{eq:MainUpperBound}. Observe that using Markov's inequality
\begin{align}
 \mathbb{P}(\fluc_T\leq -s) &= \mathbb{P} \left(\exp\left(-\exp\left(T^{\frac{1}{3}}(\fluc_{T}+s)\right)\right)\geq e^{-1}\right) \leq e\mathbb{E}\left[\exp\left(-\exp\left(T^{\frac{1}{3}}(\fluc_{T}+s)\right)\right)\right].
\end{align}
\eqref{eq:UpperBoundOnInfiniteProduct} bounds the right-hand side and yields \eqref{eq:MainUpperBound}.

Now we show that \eqref{eq:LowerBoundOnInfiniteProduct}$\Rightarrow$\eqref{eq:MainLowerBound}. Fix some $\zeta\in (0,\epsilon)$. 
Observe that
\begin{align}\label{eq:UpperBoundONInfniyteProductToGetLowerBound}
\mathfrak{R}&:=\mathbb{E}\Big[\exp\Big(-\exp\big(T^{\frac{1}{3}}(\fluc_{T}+\bar{s})\big)\Big)\Big]\\&\leq
\mathbb{E}\Big[\mathbbm{1}\{\fluc_{T}\leq -s\}+\mathbbm{1}\{\fluc_T> -s\}\exp\Big(-\exp\big(\delta \bar{s}T^{\frac{1}{3}}\big)\Big)\Big],\quad \bar{s}:= (1-\zeta)^{-1} s.
\end{align}
where $\mathbbm{1}\{A\}$ is an indicator function. The above inequality implies that 
\begin{align}\label{eq:LowerBoundOnProbability}
\mathbb{P}\left(\fluc_{T}\leq - s\right) \geq \mathfrak{R} -\exp\Big(-\exp\big(\zeta \bar{s}T^{\frac{1}{3}}\big)\Big).
\end{align}
It follows from \eqref{eq:LowerBoundOnInfiniteProduct} that
\begin{equation}\label{eq:BoundOnB}
\mathfrak{R}\geq \exp\Big(-(1+C\epsilon+ C^{\prime}\zeta)\frac{4}{15\pi}T^{\frac{1}{3}}s^{\frac{5}{2}}\Big)+ \exp(-K_2s^{3})
\end{equation}
for all $s\geq S=S(\epsilon, \delta)$. Here, the $C^{\prime}\zeta$ terms appears because $\bar{s}^{\frac{5}{2}}\leq s^{\frac{5}{2}}(1+C^{\prime}\zeta)$ for some $C'>0$. Recalling that $\zeta<\epsilon$ we can replace $C^{\prime}\zeta$ in \eqref{eq:BoundOnB} by $C^{\prime}\epsilon$.

Now, we notice that there exists $S^{\prime}= S^{\prime}(\epsilon, T_0)$ such that for all $s\geq S^{\prime}$,
\begin{equation}\label{eq:ExponentialTermBound}
\exp\big(\zeta \bar{s} T^{\frac{1}{3}}\big)\geq  T^{\frac{1}{3}}\frac{4s^{\frac{5}{2}}}{15\pi}-\log \epsilon, \qquad \textrm{and}\quad
\exp\Big(-\exp\big(\zeta \bar{s} T^{\frac{1}{3}}\big)\Big)\leq \epsilon \exp\Big(-\frac{4}{15\pi}T^{\frac{1}{3}}s^{\frac{5}{2}}\Big).
\end{equation}
Plugging the lower bound \eqref{eq:BoundOnB} on $\mathfrak{R}$ and the upper bound \eqref{eq:ExponentialTermBound} on $\exp(-\exp(\zeta \bar{s} T^{\frac{1}{3}}))$ into the right-hand side of  \eqref{eq:LowerBoundOnProbability} yields, for all $s\geq \max\{S, S^{\prime}\}$,
\[\mathbb{P}(\fluc_{T}\leq -s)\geq (1-\epsilon)\exp\Big(-(1+(C+C^{\prime})\epsilon)\frac{4}{15\pi}T^{\frac{1}{3}}s^{\frac{5}{2}}\Big)+ \exp(-K_2s^{3}).\]
The multiplicative factor $(1-\epsilon)$ can be absorbed into  the exponential factor $(1+(C+C^{\prime})\epsilon))$ on the right-hand side above; and rewriting it as $(1+C\epsilon)$ for a slightly modified constant $C$ yields the right side of \eqref{eq:MainLowerBound}, thus completing the proof of Theorem~\ref{MainResult}.
\qed

\section{Airy point process}\label{sec:Airy}
To prove Proposition~\ref{MostImportantTheorem}, we use Proposition~\ref{MomentMatch} which connects the SHE and the Airy point process. In this section we recall or prove various properties about the Airy point process. Section \ref{secdet} reviews its determinantal structure. Section~\ref{sec:MostImpTheo} contains a proof of Proposition~\ref{MostImportantTheorem}. Section \ref{SAO} relates the Airy point process to the stochastic Airy operator and derives properties about the typical point locations and deviations from there. Section~\ref{Heuristic} contains a heuristic explanation for certain terms in our tail bound. Finally, Sections \ref{sec:LRgidity}, \ref{sec:UpperTail}, and \ref{Sec:Sand} provide proofs of, respectively, Theorem~\ref{thm:LRigidityBound}, Theorem~\ref{UpperTailLemma} and Proposition~\ref{ppn:SandwichResult}.

\subsection{Determinantal point process definition}\label{secdet}

The \emph{Airy point process} (written here as $\chia$ or $\mathbf{a}_1>\mathbf{a}_2>\cdots$) is a \emph{simple determinantal point process} \cite[Section 4.2]{AGZ10}. Let us briefly review these terms. Denote the Borel $\sigma$-algebra of the real line $\RR$ by $\mathcal{B}(\RR)$ and let $\mu$ be a sigma finite measure over $\RR$. A \emph{point process} is a probability distribution on locally finite configurations of the real points, or in other words, a non-negative integer-valued random measure $\chi$ on the measure space $M=(\RR, \mathcal{B}(\RR), \mu)$. A point process $\chi$ is called \emph{simple} if $\mu\big(\{\exists x: \chi(x)\neq 0\}\big)=0$. For any $k\geq 1$, the $k$-point \emph{correlation function} of $\chi$ with respect to the measure $\mu$ is the locally integrable function $\rho_k:\RR^{k}\to [0, \infty)$ such that for any mutually disjoint families of the Borel sets $B_1, \ldots , B_k\in \mathcal{B}(\RR)$,
\begin{equation}\label{eq:CorrelationFunction}
\mathbb{E}_{\nu}\Big[\prod_{i=1}^{k} \chi(B_i)\Big]= \int_{\RR^k} \rho_k(x_1, \ldots , x_k) d\mu(x_1)\ldots d\mu(x_k).
\end{equation}
A simple point process $\chi$ is \emph{determinantal} if there exists $K^{\chi}:\RR^2\to \mathbb{C}$ such that for all $k\geq 1$, $\rho_k(x_1,\ldots ,x_k)= \det\big[K^{\chi}(x_i,x_j)\big]_{1\leq i,j\leq k}$. We refer to $K^{\chi}$ as the \emph{correlation kernel} of $\chi$.

The Airy point process correlation kernel $K^{\mathrm{Ai}}$ relative to Lebesgue measure $\mu$ on $\RR$ is\footnote{Recall the \emph{Airy function} $\mathrm{Ai}(x):= \frac{1}{\pi} \int^{\infty}_0 \cos(tx+t^3/3) dt$.}
\begin{align}\label{eq:AiryKernel}
K^{\mathrm{Ai}}(x,y)= \frac{\mathrm{Ai}(x)\mathrm{Ai}^{\prime}(y)- \mathrm{Ai}(y)\mathrm{Ai}^{\prime}(x)}{x-y} = \int_0^\infty \mathrm{Ai}(x+r)\mathrm{Ai}(y+r)dr.
\end{align}
We will write $\chia$ to denote the Airy point process (random) measure. We may also write\footnote{This follows from a calculation like in Proposition \ref{ppn:ExpVarOfLS} which shows that almost-surely there are infinitely many particles in $\chia$ but only finitely many to the right of any given point.} $\chia=\sum_{i=1}^{\infty} \delta_{\mathbf{a}_i}$ for random points $\mathbf{a}_1>\mathbf{a}_2>\cdots$. We will use both of these notations.

An integral operator $\mathfrak{K}: L^2(M)\to L^2(M)$ with kernel $K:\R^2\to \mathbb{C}$ written as
\begin{align}\label{eq:IntegralKernel}
(\mathfrak{K}f)(x)= \int K(x,y) f(y) d\mu(y), \quad \text{for } f\in L^2(M)
\end{align}
is \emph{locally admissible} if for any compact set $D\subset \RR$, the operator $\mathfrak{K}_{D}= \mathbbm{1}_{D} \mathfrak{K}\mathbbm{1}_{D}$, having kernel $K_D(x,y)= \mathbbm{1}_{D}(x) K(x,y) \mathbbm{1}_{D}(y)$, has the following representation:
\begin{align}\label{eq:LA}
(\mathfrak{K}_D f)(x)= \sum_{k=1}^n \lambda_{k} \phi_k(x) \langle \phi_k, f\rangle_{L^2(M)},\qquad K_D(x,y) = \sum_{k=1}^n \lambda_{k} \phi_{k}(x) \overline{\phi_k(y)}
\end{align}
where $n$ may be finite or infinite, $\{\phi_k\}_{k}\in L^2(M)$ are orthonormal eigenfunctions and the eigenvalues $(\lambda^{D}_k)^{n}_{k=1}$ of $K_D$ are positive and satisfy $\sum^n_{k=1}\lambda^{D}_n< \infty$. We call $\mathfrak{K}$ \emph{good} if for all compact $D$ and all $1\leq k\leq n$, $\lambda^{D}_k\in(0,1]$. For a determinantal point process with  locally admissible and good correlation kernel, for any compact set $D\subset \RR$, $\chi(D)$ equals in distribution the sum of $n$ (same $n$ as in \eqref{eq:LA}) independent Bernoulli random variables with the respective probabilities of equality to 1 given by the $\lambda^{D}_1,\ldots,\lambda^{D}_n$ -- see \cite[Section 4.2]{AGZ10}.
\begin{lemma}\label{lem.locadd}
The kernel \eqref{eq:AiryKernel} of the Airy point process $K^{\mathrm{Ai}}$ is locally admissible and good.
\end{lemma}
We use this result in proving Theorem \ref{UpperTailLemma} (see \cite[Proposition~4.2.30]{AGZ10} for a proof).

%

\subsection{Proof of Proposition~\ref{MostImportantTheorem}}\label{sec:MostImpTheo}
As above, let $\mathbf{a_1}>\mathbf{a}_2>\ldots $ denote the Airy point process.
Denote
\begin{align}\label{eq:NewNotation}
\mathcal{I}_s(x) := \frac{1}{1+\exp\big(T^{\frac{1}{3}}(s+x)\big)} \qquad \text{and} \qquad \mathcal{J}_s(x) := \log\Big(1+\exp\big(T^{\frac{1}{3}}(s+x)\big)\Big)
\end{align}
so that for any $x\in \RR$, we have $\mathcal{I}_s(x)= \exp\big(-\mathcal{J}_s(x)\big)$. Proposition~\ref{MomentMatch} connects $\mathbb{E}_{\mathrm{Airy}}\big[\prod_{k=1}^{\infty}\mathcal{I}_s(\mathbf{a}_k)\big]$ with the Laplace transform of the SHE.
We now state upper and lower bounds on this expectation and then subsequently complete the proof of Proposition~\ref{MostImportantTheorem}.
\bp\label{ppn:ImportantLemma}
Fix any $\epsilon,\delta \in (0,\frac{1}{3})$ and $T_0>0$. Then, there exist $s_0=s_0(\epsilon,\delta,T_0)$, an absolute constant $C>0$, $K_1=K_1(\epsilon,\delta, T_0)>0$ and $K_2=K_2(T_0)>0$ such that for all $s\geq s_0$ and $T\geq T_0$,
\begin{align}\label{eq:FinalUpperBound}
 \mathbb{E}_{\mathrm{Airy}}\left[\prod_{k=1}^{\infty} \mathcal{I}_s(\mathbf{a}_k)\right]\leq e^{-\frac{4(1-C\epsilon) }{15\pi}T^{\frac{1}{3}} s^{5/2}} + e^{-K_1s^{3-\delta}-\epsilon T^{1/3}s} + e^{-\frac{(1-C\epsilon)}{12}s^3}
\end{align}
 and
\begin{equation}\label{eq:FinalLowerBound}
\mathbb{E}_{\mathrm{Airy}}\left[\prod_{k=1}^{\infty} \mathcal{I}_s(\mathbf{a}_k)\right]\geq e^{-\frac{4(1+C\epsilon)}{15\pi}T^{\frac{1}{3}}s^{5/2}}+  e^{-K_2 s^{3}}.
\end{equation}
\ep

\begin{proof}[Proof of Proposition~\ref{MostImportantTheorem}]
Using \eqref{eq:SHEAiryeqn}, \eqref{eq:UpperBoundOnInfiniteProduct}-\eqref{eq:LowerBoundOnInfiniteProduct} follows from \eqref{eq:FinalUpperBound}-\eqref{eq:FinalLowerBound} of Proposition~\ref{ppn:ImportantLemma}.
\end{proof}



 \subsection{Stochastic Airy operator}\label{SAO}
As  observed in \cite{ES07} and proved in \cite{RRV11}, the Airy point process equals in distribution the negated spectrum of the \emph{stochastic Airy operator}. This yields a way to compute the typical locations of the points and establish a uniform bound (Proposition \ref{ppn:SandwichResult}) on the deviations from those locations. This bound is used in the proof of \eqref{eq:FinalLowerBound} of Proposition~\ref{ppn:ImportantLemma}. It is not, however, tight enough to suffice for all of our needs, hence our need for Theorems \ref{thm:LRigidityBound} and  \ref{UpperTailLemma}.

\bd[Stochastic Airy operator]
Let $D=D(\RR^{+})$ be the space of the generalized functions, i.e., the continuous dual of the space $C^{\infty}_0$ of all smooth compactly supported test functions endowed with the topology of compact convergence. For any function $f$, we denote its $k$-th derivative by the symbol $f^{(k)}$ and define its action on any test function $\phi\in C^{\infty}_0$ by
\begin{align}\label{eq:DefineAction}
\prec \phi, f^{(k)}(x) \succ := (-1)^k \int f(x) \phi^{(k)}(x) dx.
\end{align}
Define the space of functions $H^{1}_{\mathrm{loc}}= H^{1}_{\mathrm{loc}}(\RR)$, where for any $f\in H^{1}_{\mathrm{loc}}$ and any compact set $I\subset \RR$, we have $f^{(1)}\mathbbm{1}_{I}\in L^2(\RR)$. The $\beta>0$ stochastic Airy operator $\mathcal{H}_{\beta}$ is a linear map
\begin{align}\label{eq:SAODef}
\mathcal{H}_{\beta}:H^{1}_{\mathrm{loc}}\to D\qquad \textrm{with}\quad \mathcal{H}_{\beta} f = -f^{(2)}+xf+\frac{2}{\sqrt{\beta}} fB^{\prime}.
\end{align}
Here, $B$ is a standard Brownian motion and $B^{\prime}$ is its derivative which belongs to the space $D$\footnote{To see that $fB^{\prime}\in D$, observe that $\int_{0}^{y} fB^{\prime} dx= -\int^{y}_{0} Bf^{\prime} dx + f(y)B_y - f(0)B_0$ by integration by parts. One can now see that the latter is continuous function. Thus, its derivative $fB^{\prime}$ belongs to the space $D$.}.
The non-random part of $\mathcal{H}_{\beta}$ is the \emph{Airy operator} $\mathcal{A}=-\partial^2_x+ x$.
Define the Hilbert space
\begin{align}\label{eq:Sp-NormDef}
L^*:=\big\{f: f(0)=0, \|f\|_{*}<\infty\big\} \quad \text{where}\quad \|f\|^2_{*}= \int^{\infty}_0 \big((f^{\prime})^2 + (1+x)f^2\big) dx.
\end{align}
A pair $(f,\mathbf{\Lambda})\in L^{*}\times \RR$ is an eigenfunction/value pair for $\mathcal{H}_{\beta}$ if $\mathcal{H}_{\beta} f= \mathbf{\Lambda}f$ (likewise for $\mathcal{A}$).
\ed

\bp[\cite{RRV11}]\label{ppn:connection}
 Let $\mathbf{a}=(\mathbf{a}_1> \mathbf{a}_2> \ldots )$ denote the Airy point process and $\mathbf{\Lambda}=(\mathbf{\Lambda}_1<\mathbf{\Lambda}_2 <\ldots )$ denote the eigenvalues of $\mathcal{H}_2$. Then, $\mathbf{a}$ and $-\mathbf{\Lambda}$ are equal in distribution.
 \ep

Results obtained in \cite{RRV11,Balint14} show that the spectrum of $\mathcal{H}_{\beta}$ lies within a uniform random band around the spectrum of the Airy operator $\mathcal{A}$. The following is a strengthening of such a result wherein the tail decay of the band width (here $C_{\epsilon}$) is controlled.

\bp\label{ppn:SandwichResult}
Denote the eigenvalues of the Airy operator $\mathcal{A}$ by $(\mathbf{\lambda}_1< \mathbf{\lambda}_2<\ldots )$ and the eigenvalues of $\mathcal{H}_{\beta}$ by $(\mathbf{\Lambda}^{\beta}_1,\mathbf{\Lambda}^{\beta}_2,\ldots)$. For any $\epsilon\in(0,1)$, we define the random variable $C_{\epsilon}$ as the minimal real number such that for all $k\geq 1$,
\begin{equation}\label{eq:SandwithchProperty}
(1-\epsilon) \mathbf{\lambda}_k -C_\epsilon \leq \mathbf{\Lambda}^{\beta}_k \leq (1+\epsilon)\mathbf{\lambda}_k +C_\epsilon
\end{equation}
Then, for all $\epsilon,\delta\in (0,1)$ there exist $s_0= s_0(\epsilon, \delta)$, and $\kappa= \kappa(\epsilon,\delta)$ such that for~$s\geq s_0$,
\begin{equation}\label{eq:RandomConstantSanwitch}
\mathbb{P}\Big(C_{\epsilon}\geq \frac{s}{\sqrt{\beta}}\Big)\leq \kappa\exp\big(-\kappa s^{1-\delta}\big).
\end{equation}
\ep

Notice that \eqref{eq:RandomConstantSanwitch} demonstrates a concentration inequality for the supremum of the deviations of the eigenvalues of $\mathcal{H}_{\beta}$ around their typical locations. We defer the proof of this proposition until Section~\ref{Sec:Sand}.


Finally, we state a result on the position of the eigenvalues of the Airy operator $\mathcal{A}$. Classical works \cite{MacTit59,T58} have addressed this question for more general operators $-\partial^2_x+V(x)$ for $V(x)$ satisfying certain regularity conditions. For the Airy operator, $\lambda_k$ coincides with the $k$-th zero of the Airy function.

\bp[\cite{MacTit59}]\label{MeanPosition}
Denote the eigenvalues of the Airy operator $\mathcal{A}$ by $\mathbf{\lambda}=(\mathbf{\lambda}_1< \mathbf{\lambda}_2<\ldots )$. Then for any $n\geq 1$, $\lambda_n$ satisfies
\begin{align}\label{eq:lambda_n}
\frac{1}{\pi}\int^{\lambda_n}_0 \sqrt{\left(\lambda_n -x\right)} dx = n-\frac{1}{4} +\mathcal{R}\left(n\right),\quad \textrm{or} \quad \lambda_n =\Big(\tfrac{3\pi}{2}\big(n-\tfrac{1}{4}+\mathcal{R}(n)\big)\Big)^{\frac{2}{3}}.
\end{align}
where $|\mathcal{R}(n)|\leq K/n$ for some large constant $K$.
 \ep

\subsection{Heuristics for Proposition~\ref{ppn:ImportantLemma}}\label{Heuristic}
There are two main contributions to $\mathbb{E}_{\mathrm{Airy}}\big[\prod_{k=1}^{\infty}\mathcal{I}_s(\mathbf{a}_k)\big]$ -- typical and atypical values of $\mathbf{a}$. Owing to Proposition~\ref{ppn:SandwichResult}, the typical values of $\mathbf{a}$ are close to the negatives of the Airy operator eigenvalues, whose locations are estimated in Proposition \ref{MeanPosition}.

The asymptotic formula in \eqref{eq:lambda_n} leads (as we now show) to the $\exp\big(-\tfrac{4}{15\pi} T^{1/3} s^{5/2}\big)$ term in \eqref{eq:FinalUpperBound} and \eqref{eq:FinalLowerBound}\footnote{The $\epsilon$ error factor comes from various approximation errors and the fact that the replacement is only true with high probability.}.
Replacing $\mathbf{a}_k$ by $-\lambda_k$ yields
\begin{align}\label{eq:Approx1}
\log \bigg(\prod_{k=1}^{\infty} \mathcal{I}_s(\mathbf{a}_k)\bigg)\approx\sum_{k=1}^{\infty} \mathcal{J}_s(-\lambda_k) =-\sum_{k=1}^{\infty}\log\Big(1+\exp\big(T^{\frac{1}{3}}(s-\lambda_k)\big)\Big).
\end{align}
When $s\gg\lambda_k$ and $T$ is bounded away from $0$, $\log\big(1+\exp(T^{\frac{1}{3}}(s-\lambda_k))\big)\approx T^{\frac{1}{3}}(s-\lambda_k)$.  By Proposition~\ref{MeanPosition}, $\lambda_k\approx (3\pi k/2)^{\frac{2}{3}}$, hence
\begin{align}\label{eq:Approx2}
\sum_{\{k:\lambda_k<s\}} \log\Big(1+\exp\big(T^{\frac{1}{3}}(s-\lambda_k)\big)\Big)&\approx T^{\frac{1}{3}}\sum_{k< \frac{2}{3\pi}s^{\frac{3}{2}}}\bigg(s-\big(\frac{3\pi k}{2}\big)^{\frac{2}{3}}\bigg)\\&\approx T^{\frac{1}{3}}\bigg(\frac{2}{3\pi}s^{\frac{5}{2}}- \frac{3}{5}\cdot\big(\frac{3\pi}{2}\big)^{\frac{2}{3}}\cdot\big(\frac{2}{3\pi} s^{\frac{3}{2}}\big)^{\frac{5}{3}}\bigg)= \frac{4}{15\pi} T^{\frac{1}{3}}s^{\frac{5}{2}}.
\end{align}
To obtain the last approximation we replace the sum $\sum_{k< x}k^{\frac{2}{3}}$ by the integral $\int^{x}_{0} z^{\frac{2}{3}} dz$ which is equal to $\frac{3}{5}\cdot x^{\frac{5}{3}}$. Thus \eqref{eq:Approx2} accounts for the first term in \eqref{eq:FinalUpperBound} and \eqref{eq:FinalLowerBound}.

To complete the above heuristic we must show that the sum of $\mathcal{J}_s(-\lambda_k)$ over all $\lambda_k>s$ can be ignored. For all $\lambda_k>s$, one has $0\leq \mathcal{J}_s(-\lambda_k)\leq \exp\big(T^{\frac{1}{3}}(s-\lambda_k)\big)$. Using this,
\[0\leq \sum_{\{k:\lambda_k>s\}}\mathcal{J}_s(-\lambda_k)\leq \sum_{k\geq  \frac{2}{3\pi}s^{\frac{3}{2}}} \exp\Big(T^{\frac{1}{3}}\big(s-\big(\frac{3\pi k}{2}\big)^{\frac{2}{3}}\big)\Big) \leq \int^{\infty}_{\frac{2}{3\pi}s^{\frac{3}{2}}} \exp\Big(T^{\frac{1}{3}}\big(s-\big(\frac{3\pi z}{2}\big)^{\frac{2}{3}}\big)\Big) dz.\]
The final integrand is less than $1$ inside $[\frac{2}{3\pi}s^{\frac{3}{2}}, \infty]$ and thanks to the inequality (Lemma~\ref{ExponentLemma})
\[s- \left(\frac{3\pi z}{2}\right)^{\frac{2}{3}}\leq - \left(\frac{3\pi(z-\tfrac{2}{3\pi}s^{\frac{3}{2}})}{2}\right)^{\frac{1}{3}} \quad \text{for all  } z\geq \left(\frac{2}{3\pi}\right)s^{\frac{3}{2}}+ \sqrt{\frac{2}{3\pi}}s^{\frac{3}{4}}. \]
we obtain the following bound
\begin{align}
\int^{\infty}_{\frac{2}{3\pi}s^{\frac{3}{2}}} \exp\Big(T^{\frac{1}{3}}\big(s-\big(\frac{3\pi z}{2}\big)^{\frac{2}{3}}\big)\Big) dz \leq \sqrt{\frac{2}{3\pi}}s^{\frac{3}{4}}+\int^{\infty}_0\exp\Big(-T^{\frac{1}{3}}\big(\frac{3\pi z}{2}\big)^{1/3}\Big)
dz.
\end{align}
The final integral evaluates to a constant times $(T/2)^{-\frac{1}{3}}\int^{\infty}_0 z^2\exp(-z)dz= (T/2)^{-\frac{1}{3}}\Gamma(3)$. Thus, when $T$ is bounded away from $0$, the contribution of the eigenvalues which are greater than $s$ is of the order $\mathcal{O}(s^{\frac{3}{4}})$ which is certainly less than $s^{\frac{5}{2}}$ for enough large $s$. 

The other terms in the bounds \eqref{eq:FinalUpperBound} and \eqref{eq:FinalLowerBound} come from the atypical deviations of the Airy points from their typical locations. For instance, if $\mathbf{a}_1$ is very negative, this will clearly effect the validity of the above heuristic. The proof of Proposition~\ref{ppn:ImportantLemma} boils down to controlling these atypical deviations and measuring their effect on the multiplicative functional in question.

Before we prove Proposition~\ref{ppn:ImportantLemma}, we give proofs of Theorems \ref{thm:LRigidityBound}, \ref{UpperTailLemma} and Proposition \ref{ppn:SandwichResult} which provide important control over the atypical deviations of the Airy point process.

\subsection{Proof of Theorem~\ref{thm:LRigidityBound}}\label{sec:LRgidity}
Let us denote
$A := \Big\{\chia([-s, \infty))- \mathbb{E}\big[\chia([-s, \infty))\big]\leq  - cs^{\frac{3}{2}}\Big\}.$
Using Markov's inequality, we find that for any $\lambda>0$
\begin{align}\label{eq:Mkv}
\mathbb{P}(A)\leq \exp\Big(-\lambda cs^{\frac{3}{2}} + \lambda\mathbb{E}\big[\chia([-s, \infty))\big]\Big) \mathbb{E}\Big[\exp\big(-\lambda\chia([-s, \infty)) \big)\Big].
\end{align}
Set $\lambda=s^{\frac{3}{2}-\delta}$. Owing to Proposition \ref{ppn:ExpVarOfLS} and Theorem~\ref{thm:CGFexpansion},
\begin{align}
\mathbb{E}\big[\chia([-s, \infty))\big] & = \frac{2}{3\pi} s^{\frac{3}{2}}+\mathfrak{D}_1(s),\label{eq:ExpAiry}\\
\mathbb{E}\Big[\exp\big(-\lambda\chia([-s, \infty))]) \big)\Big] &= F(-s;\lambda)\leq\exp\Big(- \frac{2\lambda}{3\pi}s^{\frac{3}{2}}+ K s^{3-\frac{19\delta}{15}}\Big)\label{eq:MGFAiry}
\end{align}
where $K=K(\delta)$ is a large constant and $s$ is large enough. Thus
\[\mathbb{P}(A)\leq\exp\left(-cs^{3-\delta}+ Ks^{3-\frac{19\delta}{15}}+\mathfrak{D}_1(s)\right).\]
Recalling that $|\mathfrak{D}_1(s)|$ is uniformly bounded for all $s>0$, we find the desired bound.

\subsection{Proof of Theorem~\ref{UpperTailLemma}}\label{sec:UpperTail}
Fix any $k\in \ZZ_{\geq 0}$. By Lemma \ref{lem.locadd} the kernel of the Airy point process is locally admissible and good. Thus (as discussed before Lemma \ref{lem.locadd}) for any compact set $D$, $\chia(D) \stackrel{d}{=}\sum^\infty_{i} X_i$ where the $X_i$'s are independent Bernoulli random variables satisfying $\mathbb{P}\left(X_i=1\right)= 1- \mathbb{P}(X_i=0)=\lambda^{D}_i$. Here $\lambda^{D}_i$'s are the eigenvalues of the operator $\mathbbm{1}_{D}K^{\mathrm{Ai}}\mathbbm{1}_{D}$.
Choose a sequence of compact set $D_n$  increasing to the interval $\mathfrak{B}_k$.
By Bennett's concentration inequality \cite{Bennett},
\begin{align}\label{eq:ConcIneq}
\mathbb{P}\Big(\chia(D_n)- \mathbb{E}\big[\chia(D_n)\big]\geq c s^\frac{3}{2}\Big)\leq \exp\Big(-\sigma^2_nh\big(\frac{c s^{\frac{3}{2}}}{ \sigma^2_n}\big)\Big)
\end{align}
where $h(u):= (1+u)\log (1+u) - u$. By the dominated convergence theorem, as $n\to \infty$, $\mu_n:=\mathbb{E}\left[\chia(D_n)\right] \to \mathbb{E}\left[\chia(\mathfrak{B}_k)\right]  $ and $\sigma^2_n:=\mathrm{Var}(\chia(D_n))\to \mathrm{Var}(\chia(\mathfrak{B}_k))$. By Proposition~\ref{ppn:ExpVarOfLS},
 \[ \mathrm{Var}\big(\chia(\mathfrak{B}_k)\big)=\frac{11 \log s}{ 12\pi^2} + \mathfrak{D}^{(k)}_2(s)\]
where $\mathfrak{D}^{(k)}_2(\cdot)$ is bounded as $s\to \infty$. Therefore, for any given $\epsilon>0$, there exist $S_0=S_0(\epsilon)$ and $N_0=N_0(\epsilon)$ such that for all $s\geq S_0$ and $n\geq N_0$,
 \begin{equation}\label{eq:UpAndLowBound}
 \frac{11(1-\epsilon)\log s}{12\pi^2}\leq \sigma^2_n\leq \frac{11(1+\epsilon)\log s}{12\pi^2}.
 \end{equation}
Since $h(u)\geq u(\log u-1)$, we find $\sigma^2_nh(cs^{\frac{3}{2}}/\sigma^2_n)\geq cs^{\frac{3}{2}}\big(\log(cs^{\frac{3}{2}}) -\log\sigma^2_n-1\big)$. Plugging the upper bound \eqref{eq:UpAndLowBound} on $\sigma^2_n$ into this inequality and exponentiating yields
\begin{equation}\label{eq:ExplBounds}
\exp\big(-\sigma^2_nh(c s^{\frac{3}{2}}/\sigma^2_n)\big)\leq \exp\big(-cs^{\frac{3}{2}}(\log(cs^{\frac{3}{2}})-(1+\epsilon)\log\log s)\big)
\end{equation}
for all $n\geq N_0$ and $s$ sufficiently large. Now, Fatou's lemma shows
\begin{align}\label{eq:Fatou}
\mathbb{P}\Big(\chia(\mathfrak{B}_k)- \mathbb{E}\big[\chia(\mathfrak{B}_k)\big]\geq c s^3\Big) &\leq \liminf_{n\to \infty}\mathbb{P}\Big(\chia(D_n)- \mathbb{E}\big[\chia(D_n)\big]\geq cs^3\Big).
\end{align}
Owing to \eqref{eq:ConcIneq} and \eqref{eq:ExplBounds}, we find that
\begin{equation}
\text{r.h.s of \eqref{eq:Fatou}}\leq \limsup_{n\to \infty}\exp\big(-\sigma^2_nh(c s^{\frac{3}{2}}/\sigma^2_n)\big)\leq \exp\big(-cs^{\frac{3}{2}}(\log(cs^{\frac{3}{2}})-(1+\epsilon)\log\log s)\big).
\end{equation}

%
\qed

\subsection{Proof of Proposition~\ref{ppn:SandwichResult}}\label{Sec:Sand}

We start with a lemma about the tails of the distribution of Brownian motion oscillations.

\bl\label{lm:ControlRegularity}
Let $B_x$ be a Brownian motion on $[0,\infty)$ and define
\begin{equation}\label{eq:ZDef}
 Z:=\sup_{x>0}\sup_{y\in [0,1)} \frac{|B_{x+y}-B_{x}|}{6\sqrt{\log(3+x)}}.
\end{equation}
Then, letting
$
\bar{B}_x = \int_{x}^{x+1} B_y dy$ and $\bar{B}_x^\prime = \frac{d}{dx} \bar{B}_x \,\big(=B_{x+1}-B_x\big),$
we have that (1) $\max\{|\bar{B}^{\prime}_x|, |\bar{B}_{x}-B_x|\}\leq 6Z\sqrt{\log(3+x)}$, and (2) there exist $K_1,K_2,s_0>0$ such that for all $s>s_0$
\begin{equation}\label{eq:ZTai;}
\mathbb{P}\left(Z\geq s\right)\leq K_1 e^{-K_2s^{2}}.
\end{equation}
 \el
 \begin{proof}
 The proof of (1) follows from the following inequalities:
 \begin{align}
 |\bar{B}^{\prime}_x|&= |B_{x+1}-B_x|\leq 6\sqrt{\log(3+x)}\sup_{y\in [0,1)}\frac{|B_{x+y}-B_x|}{6\sqrt{\log(3+x)}}\leq 6Z\sqrt{\log(3+x)},\\
 |\bar{B}_x-B_x|& \leq \int^1_0|B_{x+y}-B_x|dy\leq \sup_{y\in [0,1)}|B_{x+y}-B_x|\leq 6Z\sqrt{\log(3+x)}.
\end{align}
 Turning to (2), for any $y\in [0,1)$,
 \begin{align}\label{eq:IncrementBd}
 |B_{x+y}-B_x|&\leq |B_{x+y}-B_{\lceil x \rceil}|+ |B_{\lceil x\rceil}- B_{\lfloor x\rfloor }|+|B_x- B_{\lfloor x\rfloor }|\\&\leq 2\sup_{y\in [0,1]}|B_{\lceil x\rceil +y}-B_{\lceil x\rceil}|+2 \sup_{y\in [0,1)}|B_{\lfloor x\rfloor +y}-B_{\lfloor x\rfloor}|.
 \end{align}
 Therefore
 \begin{equation}\label{eq:SupBd}
 \sup_{y\in [0,1)}\frac{|B_{x+y}-B_{x}|}{\sqrt{\log(3+x)}}\leq 2\sup_{y\in [0,1]}\frac{|B_{\lceil x\rceil +y}-B_{\lceil x\rceil}|}{\sqrt{\log(3+x)}}+ 2\sup_{y\in [0,1]}\frac{|B_{\lfloor x\rfloor +y}-B_{\lfloor x\rfloor}|}{\sqrt{\log(3+x)}}.
 \end{equation}
To study $Z$ we must take the sup over all positive real $x$ of the above bound. However, at the cost of replacing $3+x$ by $2+x$ in the denominator, using \eqref{eq:SupBd} we can bound $Z\leq 4W$ where
$$
W:=\max_{n\in \ZZ_{\geq 1}} \frac{W_n}{6\sqrt{\log(2+n)}},\qquad \textrm{where}\quad W_n :=\sum_{y\in [0,1)} |B_{n+y}-B_n|.
$$
The $\{W_n\}_{n\in \ZZ_{\geq 1}}$ are iid, and an application of the reflection principle shows that
\begin{align}\label{eq:RefPrinc}
\mathbb{P}(W_n\geq a)\leq 2\mathbb{P}\big( |B_{n+1}-B_{n})|\geq a/2\big)\leq \frac{2}{a}e^{-a^2/8}.
\end{align}
The union bound shows that
$$
\mathbb{P}(Z \geq s)\leq \mathbb{P}(4\mathcal{W}\geq s)= \mathbb{P}\bigg(\bigcup_{n=0}^{\infty}\frac{W_n}{6\sqrt{\log(2+n)}}\geq \frac{s}{4}\bigg)\leq \sum_{n=0}^{\infty}\mathbb{P}\Big(W_n\geq \frac{3}{2}s\sqrt{\log (2+n)} \Big).
$$
Combining this with \eqref{eq:RefPrinc} yields the desired decay bound as long as $s$ is large enough.
%
%
 \end{proof}

\begin{proof}[Proof of Proposition~\ref{ppn:SandwichResult}]
We will make use the following convention: For any two operator $A, B:H^{1}_{\mathrm{loc}}\to D$, we write $A\leq B$ if for all $f\in L^{*}$, $\prec f, Af\succ\leq \prec f, Bf\succ$. If $A\leq B$, then $\lambda^{A}_k\leq \lambda^{B}_k$ where  $\lambda^{A}_k$ and $ \lambda^{B}_k$ are $k$-th lowest eigenvalues of the operators $A$ and $B$ respectively.

In our proof we bound $\mathcal{H}_\beta$ above/below by the Airy operator plus/minus an error with well-controlled tails. This requires establishing a random operator bound on $B'$. Decomposing Brownian motion $B_x= \bar{B}_{x}+ (B-\bar{B}_{x})$ ($\bar{B}_x$ is defined in Lemma \ref{lm:ControlRegularity}) we find that for $f\in C^{\infty}_0$,
 \begin{align}\label{eq:InnProdExpn}
 \prec f, B^{\prime} f\succ = \int_{0}^{\infty}  f^2 \bar{B}^{\prime}_x dx + \int^{\infty}_{0}  f(x)f^{\prime}(x) (\bar{B}_x-B_x) dx.
 \end{align}

\medskip
\noindent\textbf{Claim:}
  Fix $\epsilon, \delta\in (0,1)$. Let $K_1=K_1(\delta)\geq 1$ (resp. $K_2=K_2(\delta)\geq 1$) be a constant such that $\sqrt{\log(3+x)}\leq x^{\delta}$ (resp. $\log(3+x)\leq x^{\delta}$) for all $x\geq K_1$ (resp. $x\geq K_2$). Define
 \begin{align}\label{eq:UDef}
 	\mathcal{U}_{\epsilon}:= \max\Bigg\{Z\bigg(\Big(\frac{Z}{\epsilon}\Big)^{\frac{\delta}{(1-\delta)}}+K^{\delta}_1\bigg), Z^{2}\bigg(\Big(\frac{Z}{\epsilon}\Big)^{\frac{2\delta}{(1-\delta)}}+K^{\delta}_2\bigg)\Bigg\}.
 	\end{align}
 Then
 \begin{align}\label{eq:BoundBrownian}
 -10\epsilon \mathcal{A} -6\big(1+\tfrac{1}{2}\epsilon^{-1}\big)\mathcal{U}_{\epsilon} \leq B^{\prime} \leq 10\epsilon \mathcal{A} +6\big(1+\tfrac{1}{2}\epsilon^{-1}\big)\mathcal{U}_{\epsilon}
 \end{align}

\noindent\textsc{Proof of Claim:}
Recall that $\bar{B}^{\prime}_x= B_{x+1}-B_x$.  From Lemma~\ref{lm:ControlRegularity}, $|\bar{B}^{\prime}_x|\leq 6Z\sqrt{\log(3+x)}$ (see \eqref{eq:ZDef} for $Z$). Thus, we will start by establishing the following bound, valid for all $x\geq 0$:
\begin{align}\label{eq:StackedInequality}
 Z\sqrt{\log(3+x)}\leq \max\bigg\{ Z \Big((Z/\epsilon)^{\frac{\delta}{(1-\delta)}}+K^{\delta}_1\Big) +\epsilon x,\, Z \sqrt{(Z/\epsilon)^{\frac{2\delta}{(1-\delta)}}+K^{\delta}_2 + \epsilon^2 x}\bigg\}.
\end{align}
We explain the derivation of the first bound by $Z \Big((Z/\epsilon)^{\frac{\delta}{(1-\delta)}}+K^{\delta}_1\Big) +\epsilon x$, as the second bound follows a similar type of argument. Let $z_0:= \max\{(Z/\epsilon)^{\frac{1}{(1-\delta)}}, K_1\}$. For $x<z_0$,
$$
Z\sqrt{\log(3+x)}\leq Z\sqrt{\log(3+z_0)} \leq Z\cdot z^{\delta}_0\leq Z\Big((Z/\epsilon)^{\frac{\delta}{(1-\delta)}}+K^{\delta}_1\Big)\leq Z\Big((Z/\epsilon)^{\frac{\delta}{(1-\delta)}}+K^{\delta}_1\Big)+\epsilon x.
$$
The second inequality uses $\sqrt{\log(3+z_0)}\leq z^{\delta}_0$, the third uses $\max\{a,b\}\leq a+b$ for $a,b\geq 0$. For $x\geq z_0$,
$$
Z\sqrt{\log(3+x)} \leq Z(1+\sqrt{\log(3+x)})\leq Z+\epsilon x^{1-\delta}\sqrt{\log(3+x)} \leq  Z+\epsilon x \leq Z\Big((Z/\epsilon)^{\frac{\delta}{(1-\delta)}}+K^{\delta}_1 \Big)+\epsilon x
$$
The second inequality uses $Z\leq \epsilon x^{1-\delta}$ (as $(Z/\epsilon)^{\frac{1}{1-\delta}}\leq x$), the third uses $\sqrt{\log(3+x)}\leq x^{\delta}$ (since $x\geq K_1$), the fourth uses $(Z/\epsilon)^{\frac{\delta}{(1-\delta)}}+K^{\delta}_1\geq 1$.

Combining \eqref{eq:StackedInequality} with the definition of $\mathcal{U}_{\epsilon}$, we see that for all $x\geq 0$,
$$
Z\sqrt{\log(3+x)}\leq \max\Big\{\mathcal{U}_{\epsilon}+ \epsilon x,\, \sqrt{\mathcal{U}_{\epsilon}+ \epsilon^2 x}\Big\}.
$$
This, along with Lemma~\ref{lm:ControlRegularity} establishes that for all $x\geq 0$,
\begin{equation}\label{eq:IncrementBound}
|\bar{B}^{\prime}_x|\leq 6Z\sqrt{\log(3+x)}\leq 6(\mathcal{U}_{\epsilon}+ \epsilon x), \quad |\bar{B}_x-B_x|\leq 6Z\sqrt{\log(3+x)}\leq 6\sqrt{\mathcal{U}_{\epsilon}+\epsilon^2 x}.
\end{equation}

Using the formula for $\prec f, B^{\prime} f\succ$ in \eqref{eq:InnProdExpn}, along with the inequality $|f^{\prime}(x)f(x)(\bar{B}_x-B_x)|\leq 3\epsilon (f^{\prime}(x))^2+ (12\epsilon)^{-1} f(x)^2|\bar{B}_x-B_x|^2$ (which follows by applying $ab\leq \tfrac{1}{2}(a^2+b^2)$) we have that
$$
\left|\prec f, B^{\prime} f\succ\right|\leq \int^{\infty}_0 f^{2}(x)(\epsilon+ |\bar{B}^{\prime}_x|)dx + 3\epsilon\int^{\infty}_0 (f^{\prime}(x))^2  + (12\epsilon)^{-1}\int^{\infty}_0 (f(x))^2|\bar{B}_x-B_x|^2 dx.
$$
Plugging the bounds from \eqref{eq:IncrementBound} into the above expression yields
\begin{align}
\left|\prec f, B^{\prime} f\succ\right| &\leq 6\mathcal{U}_\epsilon\|f\|^2 +7\epsilon \langle f, \mathcal{A}f\rangle+ 3\epsilon \int_{0}^{\infty} (f^{\prime}(x))^2dx +3\epsilon^{-1} \int^{\infty}_{0} f^{2}\left(\mathcal{U}_{\epsilon}+\epsilon^2 x\right)dx\\
   &\leq 6\big(1+(2\epsilon)^{-1}\big)\mathcal{U}_{\epsilon}\|f\|^2 + 10\epsilon \langle f, \mathcal{A} f\rangle.
   \end{align}
which implies \eqref{eq:BoundBrownian} as claimed.
%
%
\medskip


Combining \eqref{eq:BoundBrownian} with the definition that $\mathcal{H}_{\beta}= \mathcal{A}+\frac{2}{\sqrt{\beta}} B^{\prime}$ yields
$$
\mathcal{A}\Big(1-\frac{20}{\sqrt{\beta}}\epsilon\Big) - \frac{12}{\sqrt{\beta}}\Big(1+\frac{1}{2\epsilon}\Big)\mathcal{U}_{\epsilon} \leq \mathcal{H}_{\beta} \leq \mathcal{A}\Big(1+\frac{20}{\sqrt{\beta}}\epsilon\Big) + \frac{12}{\sqrt{\beta}}\Big(1+\frac{1}{2\epsilon}\Big)\mathcal{U}_{\epsilon}.
$$
Replacing $\epsilon\mapsto \frac{\sqrt{\beta}}{20}\epsilon$ and using the tail bound \eqref{eq:ZTai;} on $\mathcal{U}_{\epsilon}$ yields Proposition~\ref{ppn:SandwichResult}.%
%
\end{proof}

\section{Proof of Proposition~\ref{ppn:ImportantLemma}}\label{sec:propproof}

We prove the upper bound \eqref{eq:FinalUpperBound} in Section~\ref{FinalUpperBound} and the lower bound \eqref{eq:FinalLowerBound} in Section~\ref{FinalLowerBound}. Before giving these proofs, we recall the behavior of the tail of $\mathbf{a}_1$ (the GUE Tracy-Widom distribution). There have been numerous works \cite{BoCe12,DuBa,BBD08,RRV11} to find the exact tails of $\mathbf{a}_1$  and the below proposition follows from these (e.g. \cite[Theorem 1.3]{RRV11}).

\bp\label{TracyWidom}
Let $\mathbf{a}_1$ denote the top particle in the Airy point process (which follows the Tracy-Widom GUE distribution). Then ($o(1)$ goes to zero as $s$ goes to infinity)
 \begin{align}
 \mathbb{P}\big(\mathbf{a}_1<-s\big)&= \exp\Big(-\frac{1}{12}\big(s^3+o(1)\big)\Big).\label{eq:TracyWidomLowerTail}
 \end{align}
\ep

\subsection{Proof of the upper bound (\ref{eq:FinalUpperBound})}\label{FinalUpperBound}
Recall $\mathcal{I}_s(\cdot)$ and $\mathcal{J}_s(\cdot)$ from \eqref{eq:NewNotation}, related by $\mathcal{I}_s(\cdot) = \exp\big(\mathcal{J}_s(\cdot)\big)$.
Thus, in order to obtain an upper bound on $\mathbb{E}\big[\prod_{k=1}^{\infty} \mathcal{I}_s(\mathbf{a}_k)\big]$, we derive a lower bound on $\sum_{k=1}^{\infty} \mathcal{J}_s(\mathbf{a}_k)$ by comparing the Airy point process with the corresponding eigenvalues $\lambda_k$ of the Airy operator (Section \ref{SAO}). Let us denote
 $
 \mathcal{D}_k:= (-\lambda_k -\mathbf{a}_k)_{+}= \max\{-\lambda_k-\mathbf{a}_k, 0\}
 $.

\bl\label{LowerBoundProp}
Fix some $\epsilon\in (0,1/3)$. Denote $\theta_0= \lceil  2s^{\frac{3}{2}}/3\pi\rceil$. There exists $S_0=S_0(\epsilon)>0$ and a constant $R>0$ such that for all $s\geq S_0$,
\begin{align}\label{eq:RevLB}
 \sum_{k=1}^{\infty} \mathcal{J}_s(\mathbf{a}_k)\geq T^{\frac{1}{3}}\bigg(\frac{4s^{\frac{5}{2}}}{15\pi}\big(1-8\epsilon\big)-\sum_{k=1}^{\theta_0} \mathcal{D}_k- R\bigg).
\end{align}
 \el

 \begin{proof}
 Using monotonicity of $\mathcal{J}_s(\cdot)$ and the inequality \eqref{eq:A-ASandwitch}, we obtain the following
\begin{align}\label{eq:LowerBound}
\sum_{k=1}^{\infty} \mathcal{J}_s(\mathbf{a}_k)= \sum_{k=1}^{\infty} \mathcal{J}_s\big(-\lambda_k-(-\lambda_k-\mathbf{a}_k)_{+}+(-\lambda_k-\mathbf{a}_k)_{-}\big)\geq \sum_{k=1}^{\infty} \mathcal{J}_s(-\lambda_k -\mathcal{D}_k).
\end{align}

We divide the sum on the right side of \eqref{eq:LowerBound} into three ranges: $[1,\theta_1]$, $(\theta_1,\theta_2)$ and $[\theta_2, \infty)$ where $\theta_1$ and $\theta_2$ are defined as (recall $\mathcal{R}(n)$ from Proposition~\ref{MeanPosition})
\[\mathcal{K}:= \sup_{n\geq 1}\{|n\mathcal{R}(n)|\},\qquad\theta_1:= \lceil 4\mathcal{K}\rceil, \qquad  \theta_2:= \Big\lceil\frac{ 2 s^{3/2}}{3\pi} + \frac{1}{2} \Big\rceil.   \]
 Note that as $\theta_1$ does not depend on $s$, but $\theta_2$ does, we choose $s$ large enough so $\theta_1<\theta_2$.

\bigskip
\noindent\textbf{Claim:}
\vskip-.56in
\begin{equation}\label{eq:RearLowerBd}
\sum_{k=1}^{\theta_1} \mathcal{J}_s(-\lambda_k -\mathcal{D}_k) \geq T^{\frac{1}{3}}\Bigg(\theta_1\bigg(s-\Big(\frac{3\pi(4\mathcal{K}+1)}{2}\Big)^{\frac{2}{3}}\bigg)-\sum_{k=1}^{\theta_1} \mathcal{D}_k\Bigg).
\end{equation}

\noindent\textsc{Proof of Claim:}
Since $\log(1+\exp(a))\geq a$ for any $a\in \RR$, $\mathcal{J}_s(\cdot)\geq T^{\frac{1}{3}}(s+\cdot)$. Using this and monotonicity of $\lambda_k$ in $k$, we find that
$$
 \sum_{k=1}^{\theta_1} \mathcal{J}_s(-\lambda_k -\mathcal{D}_k)\geq \sum_{k=1}^{\theta_1} T^{\frac{1}{3}}(s- \lambda_k-\mathcal{D}_k)\geq T^{\frac{1}{3}}\Big(\theta_1(s- \lambda_{\theta_1})- \sum_{k=1}^{\theta_1} \mathcal{D}_k\Big).
$$
From \eqref{eq:lambda_n}, $\lambda_{\theta_1}\leq  (3\pi(\theta_1-\frac{1}{4}+\mathcal{K}/\theta_1)/2)^{\frac{2}{3}}=(3\pi(4\mathcal{K}+1)/2)^{\frac{2}{3}}$; hence \eqref{eq:RearLowerBd} follows immediately.

\bigskip
\noindent\textbf{Claim:}
\vskip-.56in
\begin{equation}\label{eq:MidSecLowerBound}
 \sum_{k=\theta_1+1}^{\theta_2-1} \mathcal{J}_s(-\lambda_k-\mathcal{D}_k) \geq T^{\frac{1}{3}}\bigg(\frac{4s^{\frac{5}{2}}}{15\pi}\big(1-3\epsilon\big)- (\theta_1+1)s -\sum_{k=\theta_1+1}^{\theta_2-1} \mathcal{D}_k\bigg).
\end{equation}

\noindent\textsc{Proof of Claim:} We assume that $s\geq (3\pi\epsilon^{-1}/4)^{\frac{2}{3}}(1+\epsilon)$. Observe that
  \begin{align}\label{eq:b_LowerBound}
  \sum_{k=\theta_1+1}^{\theta_2-1} \mathcal{J}_s(-\lambda_k-\mathcal{D}_k)\geq   T^{\frac{1}{3}} \sum_{k=\theta_1+1}^{\theta_2-1}\Bigg(\bigg(s -\Big(\frac{3\pi k}{2}\Big)^{\frac{2}{3}}\bigg)-\sum_{k=\theta_1+1}^{\theta_2-1} \mathcal{D}_k\Bigg).
  \end{align}
This uses $\log(1+\exp(a))\geq a$ for all $a\in \RR$ and $\lambda_k\leq (3\pi k/2)^{\frac{2}{3}}$ for all $k>\theta_1$. Now we bound
  \begin{align}\label{eq:Integral1Bound}
&\sum_{k=\theta_1+1}^{\theta_2-1} \bigg(s-\Big(\frac{3\pi k}{2}\Big)^{\frac{2}{3}}\bigg) \geq
\sum_{k=\theta_1+1}^{\theta_2-1}\bigg(s-\Big(\frac{3\pi k}{2}\Big)^{\frac{2}{3}}\bigg)\geq  \int_{\theta_1+1}^{\theta_2-1} \bigg(s-\Big(\frac{3\pi z}{2}\Big)^{\frac{2}{3}}\bigg)dz\\
&\geq  \int^{\theta_2-1}_{0}\bigg(s-\Big(\frac{3\pi z}{2}\Big)^{\frac{2}{3}}\bigg) dz - (\theta_1+1)s =(\theta_2-1)\Big(s-\frac{3\cdot(3\pi)^{\frac{2}{3}}}{5\cdot 2^{\frac{2}{3}}}(\theta_2-1)^{\frac{2}{3}}\Big) - (\theta_1+1)s
  \end{align}
Noting that $(1-\epsilon)\frac{ 2s^{\frac{3}{2}}}{3\pi}\leq \theta_2-1\leq \frac{ 2 s^{\frac{3}{2}}}{3\pi}+1$ we may bound the above expression such that combining with \eqref{eq:b_LowerBound} we arrive at the claimed inequality \eqref{eq:MidSecLowerBound}.

\medskip
Plugging into \eqref{eq:LowerBound} the bounds \eqref{eq:RearLowerBd}, \eqref{eq:MidSecLowerBound}, and $\sum_{k=\theta_2}^{\infty}\mathcal{J}_s(-\lambda_k-\mathcal{D}_k)\geq 0$ yields
\begin{equation}\label{eq:LowerBound_A}
\sum_{k=1}^{\infty} \mathcal{J}_s(\mathbf{a}_k)\geq \frac{T^{\frac{1}{3}}}{2^{\frac{1}{3}}}\bigg(\frac{4s^{\frac{5}{2}}}{15\pi}\big(1-3\epsilon\big) -s - \sum_{k=1}^{\theta_2-1} \mathcal{D}_k- \theta_1\Big(\frac{3\pi(\mathcal{K}+1)}{2}\Big)^{3/2}\bigg)
\end{equation}
To finally arrive at the desired inequality in \eqref{eq:RevLB}, we use two more bounds.
Since we may assume $s\leq \frac{4\epsilon s^{5/2}}{3\pi}$ for all $s\geq S_0$, we can replace $-s$ by $- \frac{4\epsilon s^{5/2}}{3\pi}$ in the right side of \eqref{eq:LowerBound_A}. Finally, for all $\epsilon<1$, $\theta_1(3\pi(K+1)/2)^{3/2}$ can be bounded above by a large constant $R$ (independent of $s$ and $\epsilon$).
Incorporating these bounds into \eqref{eq:LowerBound_A} yields \eqref{eq:RevLB}.
\end{proof}

\begin{proof}[Proof of (\ref{eq:FinalUpperBound}) in Proposition~\ref{ppn:ImportantLemma}]
Multiplying \eqref{eq:RevLB} by $-1$ and exponentiating yields
\[\prod_{k=1}^{\infty}\mathcal{I}_s(\mathbf{a}_k)\leq  \exp\bigg(-T^{\frac{1}{3}}\Big(\frac{4s^{\frac{5}{2}}}{15\pi}\big(1-8\epsilon\big)- \sum_{k=1}^{\theta_0}\mathcal{D}_k-R\Big)\bigg).\]
Recalling $\theta_0= \lceil  2s^{\frac{3}{2}}/3\pi\rceil$ and defining $\mathcal{S}_{\theta_0} := \sum_{k=1}^{\theta_0} \mathcal{D}_k$ we have that
\begin{equation}\label{eq:RgRegime}
\mathbbm{1}\big\{\mathcal{S}_{\theta_0}<\epsilon s\theta_0\big\}\prod_{k=1}^{\infty} \mathcal{I}_s(\mathbf{a}_k) \leq \exp\Big(-T^{\frac{1}{3}}\frac{4s^{\frac{5}{2}}}{15\pi}(1-11\epsilon)\Big).
\end{equation}
If $S_{\theta_{0}}\geq \epsilon s\theta_0$, then there exists at least one $k\in [1,\theta_0]\cap \ZZ$ such that $\mathcal{D}_k$ is greater than $\epsilon s$. Thus,
$\big\{\mathcal{S}_{\theta_0}\geq \epsilon s \theta_0\big\}\subset \bigcup_{k=1}^{\theta_0}\big\{\mathcal{D}_k\geq \epsilon s\big\}$.
%
Summarizing the discussion above, we have that
\begin{align}\label{eq:FinalStepsOfUB}
 \mathbb{E}\bigg[\prod_{k=1}^{\infty} \mathcal{I}_s(\mathbf{a}_k)\bigg]&= \mathbb{E}\bigg[\mathbbm{1}\big\{\mathcal{S}_{\theta_0}<\epsilon  s\theta_0\big\}\prod_{k=1}^{\infty} \mathcal{I}_s(\mathbf{a}_k)\bigg] +  \mathbb{E}\bigg[\mathbbm{1}\big\{\mathcal{S}_{\theta_0}\geq \epsilon s\theta_0\big\}\prod_{k=1}^{\infty} \mathcal{I}_s(\mathbf{a}_k)\bigg]\\
 &\leq  \exp\Big(-T^{\frac{1}{3}}\frac{4s^{\frac{5}{2}}}{15\pi}\big(1-11\epsilon\big)\Big) +  \mathbb{E}\bigg[\mathbbm{1}\Big\{\bigcup_{k=1}^{\theta_0}\big\{\mathcal{D}_k\geq \epsilon s\big\}\Big\}\prod_{k=1}^{\infty} \mathcal{I}_s(\mathbf{a}_k)\bigg].
 \end{align}
We may bound indicator functions
$$\mathbbm{1}\Big\{\bigcup_{k=1}^{\theta_0}\big\{\mathcal{D}_k\geq \epsilon s\big\}\Big\} \leq \mathbbm{1}\Big\{\bigcup_{k=1}^{\theta_0}\big\{\mathcal{D}_k\geq \epsilon s\big\}\cap \big\{\mathbf{a}_1\geq -(1-\epsilon)s\big\}\Big\} + \mathbbm{1}\big\{\mathbf{a}_1\leq -(1-\epsilon)s\big\}.$$
Since $\mathcal{I}_s(\mathbf{a}_k)\leq 1$ for all $k\in \ZZ_{>0}$, when $\mathbf{a}_1\geq -(1-\epsilon)s$,
$$\prod_{k=1}^{\infty} \mathcal{I}_s(\mathbf{a}_k)\leq \frac{1}{1+\exp\big(T^{\frac{1}{3}}(s+\mathbf{a}_1)\big)}\leq \exp\big(-\epsilon s T^{\frac{1}{3}}\big).$$
Combining these observations and taking expectations implies
\begin{equation}\label{eq:SplitBd2}
\mathbb{E}\bigg[\mathbbm{1}\Big\{\bigcup_{k=1}^{\theta_0}\big\{\mathcal{D}_k\geq \epsilon s\big\}\Big\} \,\prod_{k=1}^{\infty} \mathcal{I}_s(\mathbf{a}_k)\bigg] \leq  \exp\big(-\epsilon sT^{\frac{1}{3}}\big) \mathbb{P}\bigg(\bigcup_{k=1}^{\theta_0} \{\mathcal{D}_k\geq \epsilon s\}\bigg) + \mathbb{P}\big(\mathbf{a}_1\leq -(1-\epsilon)s\big).
\end{equation}
By Proposition~\ref{TracyWidom}, there exists $C>0$ such that for $s$ large enough $\mathbb{P}\big(\mathbf{a}_1\leq-(1-\epsilon)s\big)\leq\exp\big(-\frac{s^3}{12}(1-C\epsilon)\big)$. Combining \eqref{eq:FinalStepsOfUB}, \eqref{eq:SplitBd2} and \eqref{eq:MassGap} in Lemma \ref{lemma:Qint}, we find \eqref{eq:FinalUpperBound}.
\end{proof}

 \bl\label{lemma:Qint}
   Fix $\epsilon, \delta\in (0,1/3)$. There exist $S_0=S_0(\eta,\delta)>0$ and $K_1=K_1(\eta,\delta)>0$ such that the following holds for all $s\geq S_0$. Divide the interval $[-s,0]$ into $\lceil 2\epsilon^{-1}\rceil$ segments $\mathcal{Q}_i:=[-j\epsilon s/2, -(j-1)\epsilon s/2)$ for $j=1,\ldots ,\lceil 2\epsilon^{-1}\rceil$. Denote the right and left end points of $\mathcal{Q}_j$ by $q_j$ and $p_j$. Define
   $k_j:= \inf\{k: -\lambda_k\geq q_j\}$ ($\lambda_1<\lambda_2<\ldots $ are the Airy operator eigenvalues). Then (recalling $\theta_0=\lceil  2s^{\frac{3}{2}}/3\pi\rceil$),
\begin{align}\label{eq:MassGap}
\mathbb{P}\big( \mathbf{a}_{k_j}\leq p_j\big) &\leq \exp(- K_1 s^{3-\delta})\qquad \forall j\in \big\{1, \ldots , \lceil 2\eta^{-1}\rceil\big\}\\
\mathbb{P}\Big(\bigcup_{k=1}^{\theta_0}\big\{\mathcal{D}_k\geq \epsilon s\big\}\Big) & \leq \exp(-K_1 s^{3-\delta}).
\end{align}
\el
 \begin{proof}
We prove the first line of \eqref{eq:MassGap}. For $1\leq j\leq \lceil2\epsilon^{-1}\rceil$, when $\mathbf{a}_{k_j}\leq p_j=-2^{-1}(j\epsilon s)$,
\begin{align}\label{eq:Impl}
\chia\big([-2^{-1}(j\epsilon s),\infty)\big)\leq k_j\leq \#\big\{k:-\lambda_k\geq -2^{-1}(j-1)\epsilon s\big\}.
\end{align}
Owing to Propositions \ref{ppn:ExpVarOfLS} and \ref{MeanPosition}, we have
\begin{align}\label{eq:EmpExpRel}
\#\big\{k:-\lambda_k\leq -x\big\}=:\frac{2}{3\pi}x^3+C_1(x), \quad\textrm{and}\quad \mathbb{E}\big[\chia([-x,\infty))\big]=:\frac{2}{3\pi}x^3+C_2(x)
 \end{align}
where $\sup_{x\geq 0}\{|C_1(x)|,|C_2(x)|\}<\infty$. Combining \eqref{eq:Impl} and \eqref{eq:EmpExpRel} shows that when $\mathbf{a}_{k_j}\leq p_j$,
\begin{align}
\chia&([-2^{-1}(j\epsilon s),\infty)) - \mathbb{E}\left[\chia([-2^{-1}(j\epsilon s),\infty))\right]\\&\leq \#\{k:-\lambda_k\geq -2^{-1}(j-1)\epsilon s\}-\#\{k:-\lambda_k\geq -2^{-1}j\epsilon s\}+C_1(2^{-1}j\epsilon s)-C_2(2^{-1}j\epsilon s)\\
&\leq \frac{(\epsilon s)^{\frac{3}{2}}}{3\sqrt{2}\pi}\big((j-1)^{\frac{3}{2}}-j^{\frac{3}{2}}\big)+C_1(2^{-1}j\epsilon s)-C_2(2^{-1}j\epsilon s)\\&\leq -M\sqrt{j}(\epsilon s)^{\frac{3}{2}} +C_1(2^{-1}j\epsilon s)-C_2(2^{-1}j\epsilon s).
\end{align}
for some $M>0$. Therefore,
$$
\mathbb{P}\big(\mathbf{a}_{k_j}\leq p_j\big)\leq \mathbb{P}\Big(\chia([p_j, \infty))- \mathbb{E}[\chia([p_j, \infty))]\leq  - M\sqrt{j}(\epsilon s)^{\frac{3}{2}}+2\sup_{x\geq 0}\{|C_1(x)|,|C_2(x)|\}\Big).
$$
For large enough $s$, $- M\sqrt{j}(\epsilon s)^{\frac{3}{2}}+2\sup_{x\geq 0}\{|C_1(x)|,|C_2(x)|\} \leq -\frac{M}{2}\sqrt{j}(\epsilon s)^{\frac{3}{2}}$ for all $j\in \big\{1, \ldots , \lceil 2\epsilon^{-1}\rceil\big\}$. The first line of \eqref{eq:MassGap} follows by applying \eqref{eq:ChiL} of Theorem~\ref{thm:LRigidityBound} which shows that there exist $S_0(\epsilon,\delta)$ and $K_1=K_1(\epsilon,\delta)$ such that for all $s\geq S_0$,
$$
\mathbb{P}\Big(\chia([-(j\epsilon s), \infty))- \mathbb{E}\big[\chia([-(j\epsilon s), \infty))\big]\leq  - \frac{M}{2}\sqrt{j}(\epsilon s)^{\frac{3}{2}}\Big)\leq \exp\big(- K_1 s^{3-\delta}\big).
$$

Turning to the second line of \eqref{eq:MassGap}, we assume (as allowed by \eqref{eq:lambda_n}) that $s$ is large enough so $\lambda_{\theta_0}<s$. We claim then that
\begin{equation}\label{eq:inclusion}
\bigcup_{k=1}^{\theta_0}\left\{\mathcal{D}_k\geq \epsilon s \right\} \subset \bigcup_{j=1}^{\lceil 2 \epsilon^{-1}\rceil} \{\mathbf{a}_{k_j}\leq p_j\}.
\end{equation}
To see this, consider any integer $1\leq k\leq  \theta_0$ and assume that $\mathcal{D}_k\geq \epsilon s$. Let $j$ be such that $-\lambda_k\in \mathcal{Q}_{j-1}$. Since $\mathcal{Q}_{j-1}$ is to the right of $\mathcal{Q}_j=[p_j,q_j]$, it follows that $\mathbf{a}_k\leq -\lambda_k-\epsilon s$. Moreover, $\mathbf{a}_{k_j}\leq \mathbf{a}_k$ because $-\lambda_{k_j}<-\lambda_k$. Combining these yields
$$\mathbf{a}_{k_j}<\mathbf{a}_{k}\leq -\lambda_{k}- \epsilon s= (\lambda_{k_j}- \lambda_{k})-\lambda_{k_j}-\epsilon s\leq -\lambda_{k_j}-\frac{\epsilon s}{2},$$
where the last inequality uses $0\leq (\lambda_{k_j}-\lambda_{k})\leq \frac{\epsilon s}{2}$ (as $\lambda_{k_j}, \lambda_{k} \in\mathcal{Q}_{j-1}$). Hence, the distance between $\mathbf{a}_{k_j}$ and $\lambda_{k_j}$ is greater than  or equal to $\epsilon s/2$. This shows $\mathbf{a}_{k_j} \leq p_j$, and hence \eqref{eq:inclusion}.

The first line of \eqref{eq:MassGap} along with \eqref{eq:inclusion} implies that
\begin{equation}\label{eq:ClaimResult}
\mathbb{P}\Big(\bigcup_{k=1}^{\theta_0}\big\{\mathcal{D}_k\geq \epsilon s \big\}\Big)\leq \sum_{i=1}^{\theta_0} \mathbb{P}\big(\mathbf{a}_{k_i}\leq p_i\big)\leq \lceil 2\epsilon^{-1}\rceil \exp\left(- K_1 s^{3-\delta}\right).
\end{equation}
As long as $s$ is sufficiently large, the $\lceil 2\epsilon^{-1}\rceil$ prefactor can be absorbed into the exponent at the cost of slightly modifying $K_1$.
\end{proof}

\subsection{Proof of the lower bound (\ref{eq:FinalLowerBound})}\label{FinalLowerBound}
In order to obtain a lower bound on $\mathbb{E}\big[\prod_{k=1}^{\infty} \mathcal{I}_s(\mathbf{a}_k)\big]$, we derive an upper bound on $\sum_{k=1}^{\infty}\mathcal{J}_s(\mathbf{a}_k)$.
\bl\label{UpperBoundProp}
There exists $B>0$ and $S_0$ such that for all $\epsilon\in (0,1/3)$ and all $s\geq S_0$,
\begin{align}
\sum_{k=1}^{\infty}\mathcal{J}_s(\mathbf{a}_k)\leq \mathcal{L}_{T,\epsilon}(s+C^{\mathrm{Ai}}_\epsilon) \label{eq:FinalUpBound}
  \end{align}
where
\begin{align} \label{eq:DefL}
\mathcal{L}_{T,\epsilon}(x) := T^{\frac{1}{3}}\Big(\frac{4 x^{\frac{5}{2}}}{15\pi}(1+3\epsilon)+2x+B\Big)+\frac{x^{\frac{3}{2}}}{3(1-\epsilon)^{\frac{3}{2}}}+ \sqrt{\frac{2}{\pi}}\frac{x^{\frac{3}{4}}}{(1-\epsilon)^{\frac{3}{4}}}+\frac{4}{T\pi(1-\epsilon)^3}.
\end{align}
\el
\begin{proof}
Using the monotonicity of $\mathcal{J}_s(\cdot)$ and the inequality \eqref{eq:A-ASandwitch}, we obtain
\begin{align}\label{eq:UpperBound}
\sum_{k=1}^{\infty}\mathcal{J}_s(\mathbf{a}_k) \leq \sum_{k=1}^{\infty} \mathcal{J}_s\big(-(1-\epsilon)\lambda_k+C^{\mathrm{Ai}}_\epsilon\big) = (\widetilde{\mathbf{I}})+(\widetilde{\mathbf{II}})+(\widetilde{\mathbf{III}}),
\end{align}
where $(\widetilde{\mathbf{I}})$, $(\widetilde{\mathbf{II}})$ and $(\widetilde{\mathbf{III}})$ equal the sum of $\mathcal{J}_s\big(-(1-\epsilon)\lambda_k+C^{\mathrm{Ai}}_\epsilon\big)$ over all integers $k$ in the intervals $[1,\theta^{\prime}_1]$, $(\theta^{\prime}_1, \theta^{\prime}_2)$ and $[\theta^{\prime}_2,\infty)$ respectively, and (similar to Section~\ref{FinalUpperBound}) $\theta^{\prime}_1$ and $\theta^{\prime}_2$ are
\[\theta^{\prime}_1 := \theta_1=\Big \lceil 4\sup_{n\in \ZZ_{>0}}n|\mathcal{R}(n)|\Big\rceil, \qquad \theta^{\prime}_2 := \Big\lceil \frac{2(s+C^{\mathrm{Ai}}_\epsilon)^{\frac{3}{2}}}{3\pi (1-\epsilon)^{\frac{3}{2}}}+\frac{1}{2}\Big\rceil.\]

For any integer $1\leq k\leq \theta^{\prime}_1$, $\mathcal{J}_s\big(-(1-\epsilon)\lambda_k+C^{\mathrm{Ai}}_\epsilon\big)\leq \mathcal{J}_s\big(-(1-\epsilon)\lambda_1+C^{\mathrm{Ai}}_\epsilon\big)$.
%
Using this upper bound and the inequality $\log(1+\exp(a))\leq a+\pi/2$ for $a>0$, we obtain
\begin{equation}\label{eq:RearSecUpperBd}
(\widetilde{\mathbf{I}}) \leq \theta^{\prime}_1 \mathcal{J}_s\big(-(1-\epsilon)\lambda_1+C^{\mathrm{Ai}}_\epsilon\big)\leq \theta^{\prime}_1  T^{\frac{1}{3}}\big(s-(1-\epsilon)\lambda_1 +C^{\mathrm{Ai}}_\epsilon\big)+\frac{\pi\theta^{\prime}}{2}.
\end{equation}

\bigskip
\noindent\textbf{Claim:}
\vskip-.56in
\begin{align}\label{eq:MidSecUpperBound}
\hskip.1in  (\widetilde{\mathbf{II}})\leq T^{\frac{1}{3}}\Big(\tfrac{4(s+C^{\mathrm{Ai}}_\epsilon)^{\frac{5}{2}}}{15\pi}\left(1+3\epsilon\right)+\big(2-\theta^{\prime}_1\big)(s+C^{\mathrm{Ai}}_\epsilon) - \tfrac{3(3\pi)^{2/3}(\theta^{\prime}_1)^{5/3}}{5\cdot 2^{2/3}}\Big)+\tfrac{\pi(\theta^{\prime}_2-\theta^{\prime}_1)}{2}.
\end{align}

\noindent\textsc{Proof of Claim:}
  For integer $k\in (\theta^{\prime}_1, \infty)$, it follows from the definition of $\theta^{\prime}_1$ that
  \begin{equation}\label{eq:lambdaLB}
  \lambda_{k}\geq \Big(\frac{3\pi(k-\frac{1}{4}-|\mathcal{R}(k)|)}{2}\Big)^{\frac{2}{3}}\geq \Big(\frac{3\pi(k-\frac{1}{2})}{2}\Big)^{\frac{2}{3}}.
  \end{equation}
  This and monotonicity of $\mathcal{J}_s(\cdot)$ implies that
  $$
  \mathcal{J}_s\big(-(1-\epsilon)\lambda_k +C^{\mathrm{Ai}}_\epsilon\big)\leq \mathcal{J}_s\bigg(-(1-\epsilon)\Big(\frac{3\pi(k-\frac{1}{2})}{2}\Big)^{\frac{2}{3}} +C_\epsilon\bigg).
  $$
  Leveraging this and using the inequality $\mathcal{J}_s(a)\leq a+\pi/2$ for any $a>0$, we obtain
  \begin{equation}\label{eq:PrelimBoundOnII}
  (\widetilde{\mathbf{II}})\leq \sum^{\theta^{\prime}_2-1}_{k=\theta^{\prime}_1+1}\Big(T^{\frac{1}{3}} f_s(k) + \frac{\pi}{2}\Big),\qquad \textrm{where} \quad f_s(z):=s+C^{\mathrm{Ai}}_\epsilon-(1-\epsilon)\Big(\frac{3\pi(z-\frac{1}{2})}{2}\Big)^{\frac{2}{3}}.
  \end{equation}%
  Bounding the sum in \eqref{eq:PrelimBoundOnII} by the corresponding the integral we find
  \begin{align}\label{eq:IntegralBound}
  (\widetilde{\mathbf{II}})\leq T^{\frac{1}{3}}\int^{\theta^{\prime}_2}_{\theta^{\prime}_1} f_s(z) dz+\frac{\pi(\theta^{\prime}_2-\theta^{\prime}_1)}{2}.
  \end{align}
  To bound $\int^{\theta^{\prime}_2}_{\theta^{\prime}_1} f_s(z)dz$, we observe that
\begin{align}
 \int^{\theta^{\prime}_2}_{\frac{1}{2}} f_s(z)dz &\leq (s+C^{\mathrm{Ai}}_\epsilon)\left(\frac{2(s+C^{\mathrm{Ai}}_\epsilon)^{\frac{3}{2}}}{3\pi(1-\epsilon)^{\frac{3}{2}}} +\frac{3}{2}\right)- (1-\epsilon)\frac{3}{5}\cdot \left(\frac{3\pi}{2}\right)^{\frac{2}{3}}\left(\frac{2(s+C^{\mathrm{Ai}}_\epsilon)^{\frac{3}{2}}}{3\pi(1-\epsilon)^{\frac{3}{2}}}\right)^{\frac{5}{3}}  \\
 &=  \frac{4(s+C^{\mathrm{Ai}}_\epsilon)^{\frac{5}{2}}}{15\pi(1-\epsilon)^{\frac{3}{2}}}+\frac{3}{2}(s+C^{\mathrm{Ai}}_\epsilon) \leq \frac{4(s+C^{\mathrm{Ai}}_\epsilon)^{\frac{5}{2}}}{15\pi}\left(1+3\epsilon\right)+\frac{3}{2}(s+C^{\mathrm{Ai}}_\epsilon),\label{eq:Int1Bound}\\
   \int^{\theta^{\prime}_1}_{\frac{1}{2}} f(z)dz &\geq (s+C^{\mathrm{Ai}}_\epsilon)\left(\theta^{\prime}_1-\frac{1}{2}\right)-  \int^{\theta^{\prime}_1}_{\frac{1}{2}} \left(\frac{3\pi (z-1)}{2}\right)^{\frac{2}{3}} dz \\&= (s+C^{\mathrm{Ai}}_\epsilon)\left(\theta^{\prime}_1-\frac{1}{2}\right)- \frac{3}{5}\cdot\left(\frac{3\pi}{2}\right)^{\frac{2}{3}}\cdot\left(\theta^{\prime}_1-\frac{1}{2}\right)^{\frac{5}{3}} \label{eq:Int2Bound}
  \end{align}
Combining these bounds with \eqref{eq:IntegralBound} yields the upper bound on $(\widetilde{\mathbf{II}})$ in \eqref{eq:MidSecUpperBound}.

\bigskip
\noindent\textbf{Claim:}
\vskip-.56in
  \begin{align}\label{eq:TailSecUpperBound}
  (\widetilde{\mathbf{III}})\leq \sqrt{\frac{2}{\pi}}\frac{(s+C^{\mathrm{Ai}}_\epsilon)^{\frac{3}{4}}}{(1-\epsilon)^{\frac{3}{4}}}+\frac{4}{T\pi(1-\epsilon)^3}.
  \end{align}

\noindent\textsc{Proof of Claim:}
  Using the inequality $\log(1+z)\leq z $ for all $z\geq 0$, we find
  \begin{equation}\label{eq:TrivUpBound}
  \mathcal{J}_s\big(-(1-\epsilon)\lambda_k+C^{\mathrm{Ai}}_\epsilon\big)\leq \exp\Big(T^{\frac{1}{3}}\big(s-(1-\epsilon)\lambda_k+C^{\mathrm{Ai}}_\epsilon\big)\Big).
\end{equation}
Plugging the lower bound on $\lambda_k$ from \eqref{eq:lambdaLB} into \eqref{eq:TrivUpBound}, we find (recalling $f_s(z)$ from \eqref{eq:PrelimBoundOnII})
  \begin{align}\label{eq:PrelimBoindOnIII}
  (\widetilde{\mathbf{III}})\leq \sum^{\infty}_{k=\theta^{\prime}_2}\exp\big(T^{\frac{1}{3}}f_s(k)\big).
  \end{align}
Noting that $f_s(k)\leq  f_s(\theta^{\prime}_2)< 0$ for all $k>\theta^{\prime}_2$, we find that for all $k>\theta^{\prime}_2+\sqrt{3\theta^{\prime}_2}$,
$$
f_s(k)<(1-\epsilon)\Big(\frac{3\pi(\theta^{\prime}_2-\frac{1}{2})}{2}\Big)^{\frac{2}{3}}-(1-\epsilon)\Big(\frac{3\pi(k-\frac{1}{2})}{2}\Big)^{\frac{2}{3}}\leq -(1-\epsilon)\Big(\frac{3\pi (k-\theta^{\prime}_2)}{2}\Big)^{\frac{1}{3}}.
$$
The first inequality uses $f_s(\theta^{\prime}_2)<0$ and the second follows from Lemma~\ref{ExponentLemma} (we assume $s$ is large enough so $\theta^{\prime}_2-\frac{1}{2}>27$). Utilizing this estimate yields
  \begin{align}\label{eq:UpperBound_a}
  (\widetilde{\mathbf{III}}) &\leq \sum^{k=\theta^{\prime}_2+\sqrt{3\theta^{\prime}_2}}_{k=\theta^{\prime}_2}\exp\big(T^{\frac{1}{3}}f_s(k)\big)+\sum_{k>\theta^{\prime}_2+\sqrt{3\theta^{\prime}_2}}\exp\big(T^{\frac{1}{3}}f_s(k)\big)\\
  &\leq \sqrt{3\theta^{\prime}_2} + \sum_{k=\theta^{\prime}_2+\sqrt{3\theta^{\prime}_2}}^{\infty} \exp\bigg(-(1-\epsilon)T^{\frac{1}{3}}\Big(\frac{3\pi (k-\theta^{\prime}_2)}{2}\Big)^{\frac{1}{2}}\bigg)\\
  &\leq \sqrt{3\theta^{\prime}_2} + \int^{\infty}_0\exp\bigg(-(1-\epsilon)T^{\frac{1}{3}}\Big(\frac{3\pi z}{2}\Big)^{\frac{1}{3}}\bigg)dz\\
  &= \sqrt{3\theta^{\prime}_2}+ \frac{4}{T\pi(1-\epsilon)^3}\leq \sqrt{\frac{2}{\pi}}\frac{(s+C^{\mathrm{Ai}}_\epsilon)^{\frac{3}{4}}}{(1-\epsilon)^{\frac{3}{4}}}+\frac{4}{T\pi(1-\epsilon)^3}.
  \end{align}
 The first inequality follows from \eqref{eq:PrelimBoindOnIII}; the second follows from the bound
\begin{align}\label{eq:BranchBound}
 \exp\big(T^{\frac{1}{3}}f_s(k)\big)\leq \left\{\begin{matrix}
  1, & \quad  k\in [\theta^{\prime}_2, \theta^{\prime}_2+\sqrt{3\theta^{\prime}_2}],\\
 \exp\left(-(1-\epsilon)T^{\frac{1}{3}}\left(\frac{3\pi (k-\theta^{\prime}_2)}{2}\right)^{\frac{1}{3}}\right),& \quad k\in [\theta^{\prime}_2+ \sqrt{3\theta^{\prime}_2},\infty);
\end{matrix}
 \right.
\end{align}
the third uses that the sum is bounded by the integral; and the last uses that for $s$ large enough, $\sqrt{3\theta^{\prime}_2}\leq~\sqrt{\frac{2}{\pi}}\frac{(s+C^{\mathrm{Ai}}_\epsilon)^{\frac{3}{4}}}{(1-\epsilon)^{\frac{3}{4}}}$. This completes the proof of \eqref{eq:TailSecUpperBound}.

\medskip
Plugging the upper bounds of $(\widetilde{\mathbf{I}})$, $(\widetilde{\mathbf{II}})$ and $(\widetilde{\mathbf{III}})$ obtained in \eqref{eq:RearSecUpperBd}, \eqref{eq:MidSecUpperBound} and \eqref{eq:TailSecUpperBound} respectively into \eqref{eq:UpperBound}, we arrive at \eqref{eq:FinalUpBound}.
  \end{proof}

\begin{proof}[Proof of (\ref{eq:FinalLowerBound})]
\mbox{}

\noindent\textbf{Claim:}
Fix any $\epsilon,\delta\in (0,1/3)$ and $T_0>0$. Then, there exists $\kappa=\kappa(\epsilon,\delta,T_0)>0$ and $S_0=S_0(\epsilon,\delta,T_0)>0$ such that for all $s\geq S_0$ and $T>T_0$
\begin{equation}\label{eq:LowerBound1stpartF}
 \mathbb{E}_{\mathrm{Airy}}\Big[\mathbbm{1}(\mathbf{a}_1\geq -s) \prod_{k=1}^{\infty} \mathcal{I}(\mathbf{a}_k) \Big]\geq \Big(1-2\kappa\exp(-\kappa s^{1-2\delta})\Big) \exp\Big(-\frac{4s^{\frac{5}{2}}}{15\pi}(1+9\epsilon)\Big).
\end{equation}

\noindent\textsc{Proof of Claim:}
Negating both sides of \eqref{eq:FinalUpBound} and exponentiating yields $\prod_{k=1}^{\infty} \mathcal{I}(\mathbf{a}_k)\geq  \exp\big(-\mathcal{L}_{T,\epsilon}(s+C^{\mathrm{Ai}}_{\epsilon}) \big)$.
Along with the monotonicity of $\mathcal{L}_{T,\epsilon}(\cdot)$, this implies
\begin{equation}\label{eq:LowerBound1stPart}
\mathbb{E}_{\mathrm{Airy}}\Big[\mathbbm{1}(\mathbf{a}_1\geq -s) \prod_{k=1}^{\infty} \mathcal{I}(\mathbf{a}_k) \Big]\geq \mathbb{P}\big(\mathbf{a}_1\geq -s, C^{\mathrm{Ai}}_\epsilon<s^{1-\delta}\big) \exp\big(-\mathcal{L}_{T,\epsilon}(s+s^{1-\delta}) \big).
 \end{equation}
Taking $s$ large enough we have the bounds
 \begin{align}
T^{\frac{1}{3}} \frac{4(s+s^{1-\delta})^{\frac{5}{2}}}{15\pi} \leq  T^{\frac{1}{3}} \frac{4s^{\frac{5}{2}}}{15\pi}(1+5\epsilon),\qquad
T^{\frac{1}{3}}  \big(2(s+s^{1-\delta})+B\big) \leq T^{\frac{1}{3}} \frac{4s^{\frac{5}{2}}}{15\pi}\epsilon,\\
\frac{(s+s^{1-\delta})^{\frac{3}{2}}}{3(1-\epsilon)^{\frac{3}{2}}} \leq T^{\frac{1}{3}}\frac{4s^{\frac{5}{2}}}{15\pi}\epsilon, \qquad
 \sqrt{\frac{2}{\pi}}\frac{(s+s^{1-\delta})^{\frac{3}{4}}}{(1-\epsilon)^{\frac{3}{4}}}\leq T^{\frac{1}{3}}\frac{4s^{\frac{5}{2}}}{15\pi}\epsilon,\qquad
 \frac{4}{T\pi(1-\epsilon)^3} \leq T^{\frac{1}{3}}\frac{4s^{\frac{5}{2}}}{15\pi}\epsilon.
 \end{align}
Using these bounds we find that
 \begin{equation}\label{eq:LLb}
 \mathcal{L}_{T,\epsilon}\big(s+s^{1-\delta}\big) \leq T^{\frac{1}{3}}\cdot\frac{4s^{\frac{5}{2}}}{15\pi}(1+9\epsilon).
\end{equation}
Thanks to \eqref{eq:SupDevTail} of Theorem~\ref{cor:AiryTails}, there exists $\kappa=\kappa(\epsilon,\delta)$ and $S^{\prime}_0= S^{\prime}_0(\epsilon,\delta)$ such that for all $s\geq S^{\prime}_0$, $\mathbb{P}\big(C^{\mathrm{Ai}}_{\epsilon}<s^{1-\delta}\big)>1- \kappa \exp(-\kappa s^{1-2\delta})$. Moreover, using \eqref{eq:TracyWidomLowerTail}, we find that for large enough $s$, $\mathbb{P}\big(\mathbf{a}_1 \leq -s \big)\leq \kappa\exp\big(-\kappa s^{1-2\delta}\big)$. This implies that for large enough $s$,
\begin{align}\label{eq:ProbLB}
\mathbb{P}\big(\mathbf{a}_1\geq -s, C^{\mathrm{Ai}}_\epsilon<s^{1-\delta}\big)\geq \mathbb{P}\big(\mathbf{a}_1\geq -s\big) + \mathbb{P}\big( C^{\mathrm{Ai}}_\epsilon<s^{1-\delta}\big)-1\geq 1-2\kappa\exp\big(-\kappa s^{1-2\delta}\big).&&&
\end{align}
Plugging this and  \eqref{eq:LLb} into \eqref{eq:LowerBound1stPart} yields \eqref{eq:LowerBound1stpartF}.

\medskip
\noindent \textbf{Claim:} Fix $\epsilon\in (0,1/3)$ and $T_0>0$. Then, there exists $K=K(\epsilon, T_0)>0$ and $S_0= S_0(\epsilon, T_0)>0$   such that for all $s\geq S_0$,
 \begin{align}\label{eq:LowerBound2stpartF}
\mathbb{E}_{\mathrm{Airy}}\Big[\mathbbm{1}(\mathbf{a}_1< -s) \prod_{k=1}^{\infty} \mathcal{I}(\mathbf{a}_k) \Big] \geq  \exp\big(-Ks^{3}\big).
 \end{align}

\noindent \textsc{Proof of Claim:}
We begin with a brief description of our proof technique. Let us denote $\theta^{\prime}_0:=\lceil s^{1+\delta}\rceil$. We consider a finite of sequence of intervals
$$
\mathfrak{I}_1:=[-s^2,-s), \mathfrak{I}_2:=[-2s^2, -s^2), \ldots , \mathfrak{I}_{\theta^{\prime}_0}:=[- \theta^{\prime}_0 s^2, - (\theta^{\prime}_0 -1)s^2).
$$
The length of each of the interval is $s^2$ and there are $\theta^{\prime}_0$ intervals. For any integer $\ell\in (1,\theta^{\prime}_0]\cap \ZZ$ (resp. $=1$), note that $\sum_{\mathbf{a}_k\in \mathfrak{I}_{\ell}} \mathcal{J}_s(\mathbf{a}_k)$ is less than or equal to $\sum_{\mathbf{a}_k\in \mathfrak{I}_{\ell}}\mathcal{J}_s(-(\ell-1) s^2)$ (resp. $\sum_{\mathbf{a}_k\in \mathfrak{I}_1}\mathcal{J}_s(-s)$) with equality when all the $\mathbf{a}_k$ in the interval $\mathfrak{I}_{\ell}$ coincide with the right end point $-(\ell-1)s^2$ (resp. $-s$). We show that with high probability the number of Airy points inside the interval $\mathfrak{I}_{\ell}$ cannot differ considerably from its expected value. Based on this, we argue that the probability of an abundant accumulation of the Airy points inside any of the intervals $\mathfrak{I}_1, \ldots , \mathfrak{I}_{\theta^{\prime}_0}$ is small in comparison to $\mathbb{P}(\mathbf{a}_1\leq -s)$. Moreover, the contributions of those Airy points which fall into any of those intervals are bounded from above by the result of moving the points to the right endpoint of the interval. Finally, using the upper tail estimate of $C^{\mathrm{Ai}}_\epsilon$ (see \eqref{eq:SupDevTail} of Theorem~\ref{cor:AiryTails}), we show that the $\mathbf{a}_k$'s which fall in the region $(-\infty, -\theta^{\prime}_0 s^2)$ hardly contribute to the product $\sum_{k=1}^{\infty}\mathcal{I}_s(\mathbf{a}_k)$.

Now, we provide  the details of the above sketch.
First, we find an upper bound on $\sum_{\mathbf{a}_k\in \widetilde{\mathfrak{I}}} \mathcal{J}_s(\mathbf{a}_k)$ where $\widetilde{\mathfrak{I}}:= \cup_{\ell=1}^{\theta^{\prime}_0} \mathfrak{I}_\ell$. Recall that the number of $\mathbf{a}_k$'s in a Borel set $D$ is given by $\chia(D)$. By replacing all the $\mathbf{a}_k$'s inside the interval $\mathfrak{I}_k$ by the right endpoint of the interval, we obtain
\begin{align}\label{eq:FreezeBound}
\sum_{\mathbf{a}_k\in \mathfrak{I}_\ell} \mathcal{J}(\mathbf{a}_k) \leq \left\{\begin{matrix}
\chia(\mathfrak{I}_\ell)\log\Big(1+\exp\big(T^{\frac{1}{3}}\big(s-(\ell-1)s^2\big)\big)\Big) & \quad \text{when }\ell>1,\\
\chia(\mathfrak{I}_1)\log(2) & \quad \text{when }\ell=1.
\end{matrix}\right.
\end{align}
Next, using Theorem~\ref{UpperTailLemma}, we observe that for large enough $s$, $\chia(\mathfrak{I}_l)$ is bounded above by $\mathbb{E}[\chia(\mathfrak{I}_\ell)]+\epsilon s^3$ with probability greater than $1-K_1\exp(-K_2s^3\log s)$. Owing to Proposition \ref{ppn:ExpVarOfLS}, there exists a constant $M$ such that for large enough $s$,
$$\mathbb{E}[\chia(\mathfrak{I}_\ell)] = \frac{2}{3\pi}[\ell^{\frac{3}{2}}- (\ell-1)^{\frac{3}{2}}]s^{3} +\mathfrak{D}_1(\ell s^2)-\mathfrak{D}_1((\ell-1)s^2)\leq \frac{M\sqrt{\ell}s^3}{\pi}.$$
Consequently, with probability exceeding $1-\theta^{\prime}_0 K_1\exp\big(-K_2s^3 \log s\big)$
\begin{align}\label{eq:FreezingBound}
\sum_{\mathbf{a}_k\in \widetilde{\mathfrak{I}}} \mathcal{J}_s(\mathbf{a}_k) &\leq \Big(\frac{M s^3}{\pi} + \epsilon s^3\Big)\,\bigg(\log 2+\sum_{\ell=2}^{\theta^{\prime}_0} \sqrt{\ell}\log \Big(1+ \exp\big(T^{\frac{1}{3}}(s-(l-1)s^2)\big)\Big)\bigg).
\end{align}
Since $\log(1+x)\leq x$ for all $x>0$ and $(\ell -1)s^2-s\geq (l-1)s^{2}(1-\epsilon)$ for all $s\geq \epsilon^{-1}$, we conclude there there exists a constant $C$ such that for large enough $s$, with probability exceeding $1-\theta^{\prime}_0 K_1\exp\big(-K_2s^3 \log s\big)$
\begin{equation}
\sum_{\mathbf{a}_k\in \widetilde{\mathfrak{I}}} \mathcal{J}_s(\mathbf{a}_k) \leq s^3\Big(\frac{M}{\pi} + \epsilon \Big)\,\bigg(\log 2+\sum_{\ell=2}^{\theta^{\prime}_0} \sqrt{\ell}\exp\big(-(\ell-1)(1-\epsilon)T^{\frac{1}{3}}s^2\big)\bigg)\leq Cs^{3}.\label{eq:FreezingFinalBound}
\end{equation}

We now turn to bound the remaining sum $\sum_{\mathbf{a}_k <-\theta^{\prime}_0 s^2} \mathcal{J}_s(\mathbf{a}_k)$. For this, we consider the following decomposition
     \begin{align}\label{eq:DeComSum}
     \sum_{k:\mathbf{a}_k< -\theta^{\prime}_0 s^2} \mathcal{J}_s(\mathbf{a}_k) &= (\mathbf{A}) + (\mathbf{B}), \\
    (\mathbf{A}):= \sum_{k:\mathbf{a}_k<- \theta^{\prime}_0 s^2, \lambda_k \leq   \theta^{\prime}_0 s^2} \mathcal{J}_s(\mathbf{a}_k) &, \qquad (\mathbf{B}):= \sum_{k:\mathbf{a}_k<- \theta^{\prime}_0 s^2, \lambda_k >   \theta^{\prime}_0 s^2} \mathcal{J}_s(\mathbf{a}_k).
\end{align}
Proposition~\ref{MeanPosition} shows that $\#\big\{\lambda_k\leq \theta^{\prime}_0 s^2\big\}\leq C s^{\frac{9}{2}+\frac{3\delta}{2}}$ for large enough $s$ and some constant $C>0$. This, along with the bound $\log(1+a)\leq a$ for all $a>0$ implies
$$
\mathcal{J}_s(\mathbf{a}_k)\leq \exp\big(T^{\frac{1}{3}}(s-\theta^{\prime}_0 s^2)\big)\leq \exp\big(-(1-\epsilon)T^{\frac{1}{3}}s^3\big)
$$
 when $\mathbf{a}_k\leq -\theta^{\prime}_0s^2$ and $s\geq \epsilon^{-\frac{1}{2}}$. Thus, for large enough $s$,
\begin{equation}\label{eq:Aboundss}
(\mathbf{A}) \leq C s^{\frac{9}{2}+\frac{3\delta}{2}}\exp\big(-(1-\epsilon)T^{\frac{1}{3}}s^{3}\big)\leq s^3.
\end{equation}

Now, we turn to bound $(\mathbf{B})$. Recall the inequality $\mathcal{J}_s(\mathbf{a}_k)\leq \mathcal{J}_s\big(-(1-\epsilon)\lambda_k+C^{\mathrm{Ai}}_\epsilon\big)$ which we obtain by using monotonicity of $\mathcal{J}_s$ and the inequality \eqref{eq:A-ASandwitch}. We will now employ Theorem~\ref{cor:AiryTails}, but to avoid confusion in notion let us temporarily rename the variables $s$ and $\delta$ in the statement of Theorem~\ref{cor:AiryTails} by $\tilde{s}$ and $\tilde{\delta}$. Then, taking $\tilde{s} = s^{3+\frac{\delta}{2}}$ and $\tilde{\delta}=\frac{\delta}{2(3+\delta/2)}$, the corollary implies there exists $\kappa_1= \kappa_1(\epsilon,\delta)>0$, $\kappa_2= \kappa_2(\epsilon,\delta)>0$ and $S_0= S_0(\epsilon,\delta)>0$ such that for all $s\geq S_0$, $\mathbb{P}\big(C^{\mathrm{Ai}}_\epsilon<s^{3+\frac{\delta}{2}}\big)\geq 1-\kappa_1\exp\big(-\kappa_2s^{3+\frac{\delta}{4}}\big)$. Since $\theta^{\prime}_0 s^2\approx s^{3+\delta}$, we have $s+s^{3+\frac{\delta}{2}}\leq (1-\epsilon)\theta^{\prime}_0s^2$ for large enough $s$. Consequently, for large enough $s$
 \begin{align}\label{eq:Bbound}
 \mathbb{P}\Big((\mathbf{B})\leq \sum_{\lambda_k> \theta^{\prime}_0 s^2} \mathcal{J}_s\big((1-\epsilon)( \theta^{\prime}_0 s^2- \lambda_k)-s\big)\Big) \geq 1-\kappa_1\exp\big(-\kappa_2s^{3+\frac{\delta}{4}}\big).
 \end{align}

Plugging the inequality \eqref{eq:JsBound} in Lemma~\ref{lm:Intermidiary} into \eqref{eq:Bbound} and using \eqref{eq:Aboundss} along with the fact that $(\theta^{\prime}_0s^2)^{\frac{3}{4}}\leq Cs^3$ for some constant $C$, we find that for large enough $s$
\begin{equation}\label{eq:ABbound}
\mathbb{P}\big((\mathbf{A})+(\mathbf{B})\leq Cs^3\big)\geq 1- \kappa_1\exp\big(-\kappa_2 s^{3+\frac{\delta}{4}}\big).
\end{equation}
Combining this with the probability bound computed on the event in \eqref{eq:FreezingFinalBound} implies that there exists a constant $C=C(\epsilon,\delta,T_0)>0$ such that for $s$ large enough
   \begin{align}\label{eq:LastClaimIneq}
\mathbb{P}\big(\mathcal{A}\big)\geq 1-\theta^{\prime}_0 K_1\exp(-K_2s^3\log s)-\kappa_1\exp(-\kappa_2s^{3+\frac{\delta}{4}}),
   \end{align}
where $\mathcal{A} :=\Big\{\sum_{k=1}^{\infty} \mathcal{J}_s(\mathbf{a}_k)\leq Cs^{3}\Big\}$. Negating both sides above, exponentiating, then multiplying by $\mathbbm{1}(\mathbf{a}_1\leq -s)$ and taking expectation, we obtain
\begin{align}\label{eq:FinalStep}
\mathbb{E}_{\mathrm{Airy}}\Big[\mathbbm{1}(\mathbf{a}_1\leq -s)\prod_{k=1}^{\infty} \mathcal{I}_s(\mathbf{a}_k)\Big]\geq \mathbb{P}\big(\{\mathbf{a}_1\leq -s\}\cap \mathcal{A}\big) \exp(-Cs^3).
\end{align}
It thus remains to estimate
\begin{align}\label{eq:FinalProBound}
\mathbb{P}\big(\{\mathbf{a}_1\leq -s\}\cap \mathcal{A}\big)&\geq \mathbb{P}(\mathbf{a}_1\leq -s)+\mathbb{P}(\mathcal{A})-1\\&\geq \exp(- s^3)-\theta^{\prime}_0 K_1\exp(-K_2s^3\log s)-\kappa_1\exp(-\kappa_2s^{3+\frac{\delta}{4}}).
\end{align}
The first inequality  uses $\mathbb{P}(A\cap B)\geq \mathbb{P}(A)+\mathbb{P}(B)-1$ for any events $A$ and $B$. The second uses the lower bound on $\mathbb{P}(\mathbf{a}_1\leq -s)$ in \eqref{eq:TracyWidomLowerTail} and the lower bound in \eqref{eq:LastClaimIneq}. Combining \eqref{eq:FinalProBound} with \eqref{eq:FinalStep} readily yields the claimed inequality \eqref{eq:LowerBound2stpartF} for some $K$ and $s$ large enough.
%
%

\medskip
Now we may complete the proof of \eqref{eq:FinalLowerBound} by combining \eqref{eq:LowerBound1stpartF} and \eqref{eq:LowerBound2stpartF} with
\begin{equation}\label{eq:LBSplit}
\mathbb{E}\Big[\prod_{k=1}^{\infty} \mathcal{I}_s(\mathbf{a}_k)\Big]= \mathbb{E}\Big[\mathbbm{1}(\mathbf{a}_1\geq -s)\prod_{k=1}^{\infty} \mathcal{I}_s(\mathbf{a}_k)\Big] + \mathbb{E}\Big[\mathbbm{1}(\mathbf{a}_1<-s)\prod_{k=1}^{\infty} \mathcal{I}_s(\mathbf{a}_k)\Big].
\end{equation}
\end{proof}


\bl\label{lm:Intermidiary}
As above, set $\theta^{\prime}_0= \lceil s^{1+\delta}\rceil$. Then, for all $s$ such that $\theta^{\prime}_0s^2>27$,
 \begin{equation}\label{eq:JsBound}
  \sum_{\lambda_k> \theta^{\prime}_0 s^2}\mathcal{J}_s\big((1-\epsilon)(\theta^{\prime}_0 s^2- \lambda_k)-s\big)\leq \sqrt{\frac{2}{\pi}}(\theta^{\prime}_0s^2)^{\frac{3}{4}}\log 2 + \frac{4}{T\pi(1-\epsilon)^3}.
\end{equation}
 \el
 \begin{proof}
For $s$ large enough, \eqref{eq:lambda_n} implies that
\begin{align}\label{eq:subseteq}
\big\{k:\lambda_k>\theta^{\prime}_0 s^2\big\}\subseteq \Big\{k: k\geq \frac{2}{3\pi}(\theta^{\prime}_0 s^2)^{\frac{3}{2}}\Big\}.
\end{align}
This implies that following (first) inequality
    \begin{align}
  &\sum_{\lambda_k> \theta^{\prime}_0 s^2} \mathcal{J}_s\big((1-\epsilon)(\theta^{\prime}_0 s^2- \lambda_k)-s\big)\leq\sum_{k\geq \frac{2}{3\pi}(\theta^{\prime}_0s^2)^{\frac{3}{2}}} \mathcal{J}_s\big((1-\epsilon)(\theta^{\prime}_0 s^2- \lambda_k)-s\big)\\&\quad\leq \sqrt{\frac{2}{\pi}}(\theta^{\prime}_0s^2)^{\frac{3}{4}}\log 2+ \sum_{k^{\prime}=\sqrt{\frac{2}{\pi}}(\theta^{\prime}_0s^2)^{\frac{3}{4}}}^{\infty} \exp\bigg(-(1-\epsilon)T^{\frac{1}{3}}\Big(\frac{3\pi (k^{\prime}-\frac{1}{2})}{2}\Big)^{1/3}\bigg). \label{eq:UseExpLem}
    \end{align}
To show the second inequality, let $\theta^{\prime\prime}_0:= \frac{2}{3\pi}(\theta^{\prime}_0 s^2)^{\frac{3}{2}}$ and $\theta^{\prime\prime\prime}_0:=\frac{2}{3\pi}(\theta^{\prime}_0 s^2)^{\frac{3}{4}}+ \sqrt{\frac{2}{\pi}}(\theta^{\prime}_0 s)^{\frac{3}{4}}$.
Using $\mathcal{J}_s(x)\leq \log 2$ and $\log(1+x)\leq x$ for $x\leq 0$, along with Lemma~\ref{ExponentLemma} (similarly to \eqref{eq:UpperBound_a}),
\begin{align}
\mathcal{J}_s\big((1-\epsilon)(\theta^{\prime}_0 s^2-\lambda_k)-s\big)\leq \left\{
  \begin{matrix}
  \log 2 &  k\in [\theta^{\prime\prime}_0, \theta^{\prime\prime\prime}_0]\cap \ZZ\\
   \exp\bigg(-(1-\epsilon)T^{\frac{1}{3}}\Big(\frac{3\pi(k-\theta^{\prime\prime}_0-\frac{1}{2})}{2}\Big)^{\frac{1}{3}}\bigg), &   k\in (\theta^{\prime\prime\prime}_0,\infty)\cap \ZZ.
   \end{matrix}
   \right.
   \end{align}
Using this bound and substituting $k^{\prime}= k-\theta^{\prime\prime}_0$, we obtain
$$
\sum_{k\geq \frac{2}{3\pi}(\theta^{\prime}_0s^2)^{\frac{3}{2}}} \mathcal{J}_s\big((1-\epsilon)(\theta^{\prime}_0 s^2- \lambda_k)-s\big)\leq \big(\theta^{\prime\prime\prime}_0-\theta^{\prime\prime}_0\big)\log2+\sum_{k^{\prime}>\theta^{\prime\prime\prime}_0-\theta^{\prime\prime}_0}\exp\bigg(-(1-\epsilon)T^{\frac{1}{3}}\Big(\frac{3\pi (k^{\prime}-\frac{1}{2})}{2}\Big)^{\frac{1}{3}}\bigg),
$$
which implies the second inequality in \eqref{eq:UseExpLem}. Bounding the sum by a corresponding integral  and evaluating yields the bound in \eqref{eq:JsBound}.
     \end{proof}

\bl\label{ExponentLemma}
  Fix $a>27$. Then, we have $(a+x)^{\frac{2}{3}}\geq a^{\frac{2}{3}}+x^{\frac{1}{3}}$ for all $x\geq \sqrt{3a}$.
  \el
  \begin{proof}
  Observe that for all $x\geq \sqrt{3a}$ and $a>27$, one can write $x<x^2$, and using $3a^{\frac{2}{3}}\leq a$, one has $3a^{\frac{4}{3}}x^{\frac{1}{3}}\leq ax$ and $3a^{\frac{2}{3}}x^{\frac{2}{3}}\leq ax^{\frac{2}{3}}\leq ax$.
 Combining these inequalities yields
$$
(a+x)^2=a^2 + x^2 + 2ax \geq a^2 +x + 3a^{\frac{4}{3}}x^{\frac{1}{3}}+3a^{\frac{2}{3}}x^{\frac{2}{3}}=(a^{\frac{2}{3}}+x^{\frac{1}{3}})^3.
$$
  \end{proof}

\section{Ablowitz-Segur solution of Painlev\'{e} II}\label{App:ASSol}

Recall (Section~\ref{secablowitz}) the Ablowitz-Segur solution $u_{\mathrm{AS}}(\cdot;\gamma)$ of  Painlev\'e II. We restate \cite[Theorem~1.10]{Bothner15} which provides the asymptotic form of $u_{\mathrm{AS}}(x;\gamma)$ as $x\to-\infty$. Lemmas~\ref{lem:cosform} and ~\ref{lem:FinalLemma} result from analyzing this and combine (in Section \ref{sec:CGF}) to yield a proof of Theorem~\ref{thm:CGFexpansion}.

 \bp[Theorem~1.10 of \cite{Bothner15}]\label{ppn:bothner}
 Let $u_{\mathrm{AS}}(x;\gamma)$ be the Ablowitz-Segur solution of Painlev\'e-II (Section \ref{secablowitz}) where $v:=-\log(1-\gamma)$. Denote $\tau = \frac{v}{(-x)^{3/2}} \in (0,\infty)$ and define $\kappa=\kappa(\tau)\in (0,1)$ implicitly as follows
 \begin{align}
 \tau = \frac{2}{3}\sqrt{\frac{2}{1+\kappa^2}} \left[E(\kappa^{\prime})- \frac{2\kappa^2}{1+\kappa^2}K(\kappa^{\prime})\right]
 \end{align}
 where
 \begin{align}
 \kappa^{\prime} = \sqrt{1-\kappa^2},\qquad
 K(\kappa)= \int^{1}_0 \frac{d\xi}{\sqrt{(1-\xi^2)(1-\kappa^2 \xi^2)}}, \qquad
 E(\kappa)= \int^{1}_0 \sqrt{\frac{1-\xi^2 \kappa^2}{1-\xi^2}} d\xi.
 \end{align}
 Further, define
 \begin{align}
V(\tau) := -\frac{2}{3\pi} \sqrt{\frac{2}{1+\kappa^2}} \Big(E(\kappa) -  \frac{1-\kappa^2}{1+\kappa^2} K(\kappa)\Big), \quad \chi = \chi(\tau) := 2\mathbf{i} \frac{K(\kappa)}{K(\kappa^{\prime})},
 \end{align}
and define the Jacobi theta and elliptic functions (with $q=e^{\mathbf{i}\pi\chi}$ and $z\in \CC$)
 \begin{align}
 \theta_2(z,q) = 2\sum_{m=0}^{\infty} q^{\left(m+\frac{1}{2}\right)^2} \cos((2m+1)\pi z),\qquad \theta_3(z,q) = 1+ 2\sum_{m=1}^{\infty} q^{m^2} \cos(2\pi mz),\\
 \mathrm{cd}\left(2zK\left(\frac{1-\kappa}{1+\kappa}\right), \frac{1-\kappa}{1+\kappa}\right) = \frac{\theta_3(0,q)\theta_2(z,q)}{\theta_2(0,q)\theta_3(z,q)}, \quad z\in \mathbb{C}\backslash \bigcup_{m,n\in \ZZ}\left\{\frac{1}{2}+ \frac{\chi}{2} + m+ \chi n\right\}.
  \end{align}
Then, for any fixed $\zeta\in (0, \frac{2\sqrt{2}}{3})$, there exist $x_0=x_0(\zeta)>0$, $x_1=x_1(\zeta)>0$, $c_0=c_0(\zeta)>0$, $c_1=c_1(\zeta)>0$, and $v_1=v_1(\zeta)>0$ such that (using $\kappa=\kappa(\tau)$ and $V=V(\tau)$ for short)
\begin{align}\label{eq:ASSol}
u_{\mathrm{AS}}(x;\gamma)= -\sqrt{-\frac{x}{2}}\frac{1-\kappa}{\sqrt{1+\kappa^2}}\,\mathrm{cd}\left(2(-x)^{3/2} VK\left(\frac{1-\kappa}{1+\kappa}\right), \frac{1-\kappa}{1+\kappa}\right) +J_1(x;\gamma),
\end{align}
with
\begin{align}\label{eq:J1Bound1}
|J_1(x;\gamma)|&\leq c_0(-x)^{-\frac{1}{10}},  &&\forall (-x)\geq x_0, 0< v\leq (-x)^{\frac{3}{2}}\left(\frac{2\sqrt{2}}{3}-\zeta\right),&&\\
|J_1(x;\gamma)|&\leq \frac{c_1}{\log(-x)},  &&\forall (-x)\geq x_1, v\geq v_1, \frac{2}{3}\sqrt{2}(-x)^{\frac{3}{2}}- \zeta \leq v < \frac{2\sqrt{2}}{3}(-x)^{\frac{3}{2}}.&&
\end{align}
\ep

The following result continue with the notation of Proposition~\ref{ppn:bothner}.

\bp[Proposition 3.2 and Corollary 3.3 of \cite{Bothner15}] \label{lem:bothner15}
There exist $\tau_0>0$ and $C= C(\tau)$  such that for all $\tau\leq \tau_0$,
\begin{align}\label{eq:KappaLimit}
\kappa(\tau) &= 1- 2\sqrt{\frac{\tau}{\pi}} + \frac{2\tau}{\pi} - \frac{29}{8}\left(\frac{\tau}{\pi}\right)^{3/2} + \mathfrak{Q}_1(\tau),\\
V(\tau) &= - \frac{2}{3\pi} - \frac{\tau}{2\pi^2} \log \tau + \frac{\tau}{2\pi^2}(1+\log 16\pi) + \mathfrak{Q}_2(\tau).\label{eq:VLimit}
\end{align}
where $|\mathfrak{Q}_1(\tau)|\leq C\tau^2$ and $|\mathfrak{Q}_2(\tau)|\leq C\tau^2$ for all $\tau\leq \tau_0$.
\ep

\bl\label{lem:cosform}
Fix $\eta_0\in (0,2/5)$ and let $v= (-x)^{3/2-\eta}$ for any $\eta\in (\eta_0, 2/5)$. Define
\begin{equation}\label{eq:phi}
\phi(x) = \pi(-x)^{\frac{3}{2}} V(\tau) + \frac{2}{3}(-x)^{\frac{3}{2}} - \frac{v}{2\pi}\log(8(-x)^{\frac{3}{2}}).
\end{equation}
Then, there exists $x_0=x_0(\eta_0)>0$, $C=C(\eta_0)>0$ and $C^{\prime}=C^{\prime}(\eta_0)>0$  such that 
\begin{align}\label{eq:LimitingValue}
u_{\mathrm{AS}}(x;\gamma) &=  (-x)^{-\frac{1}{4}}\sqrt{\frac{v}{\pi}} \cos\left(\pi(-x)^{\frac{3}{2}}V(\tau)\right) + J_2(x),\\
\phi(x) &=  \frac{v}{2\pi} \big((1+2\pi) - \log(v) + J_3(x)\big)\label{eq:FormOfPhi}
\end{align}
where 
 $|J_2(x)|\leq C(-x)^{\frac{1}{2}-\frac{3\eta}{2}}$ and $|J_3(x)|\leq C^{\prime}(-x)^{-2\eta}$ for all $x\geq x_0$.
\el																																																
\begin{proof}

 Using \cite[(22.11.4)]{NIST}, we get
\begin{align}\label{eq:FormOfCd}
\mathrm{cd}(z,\kappa) = \frac{2\pi}{K(k)k} \sum_{n=1}^{\infty}(-1)^n \frac{q^{n+\frac{1}{2}}}{1-q^{2n+1}} \cos\big((2n+1)\zeta\big)
\end{align}
where $\zeta = \frac{\pi z}{2K(\kappa)} $ and $q = e^{\mathbf{i}\pi\chi} = \exp(- \pi K(\kappa^{\prime})/K(\kappa))$.
\medskip

\noindent\textbf{Claim:} There exists $0\leq\kappa_0<1$ and some constant $C_1=C_1(\kappa_0)>0$ such that for all $\kappa\leq \kappa_0$,
\begin{align}\label{eq:CDSausage}
\cos\big(\pi z/2K(\kappa)\big)-C_1\kappa^2\leq \mathrm{cd}(z,\kappa)\leq \cos\big(\pi z/2K(\kappa)\big)+ C_1\kappa^2.
\end{align}

\noindent\textsc{Proof of Claim:}
 Thanks to \cite{Fettis} and \cite[(19.5.5), (19.5.8)]{NIST}, there exist $0\leq \kappa_0<1$ and $0<C_2=C_2(\kappa_0)<C_3=C_3(\kappa_0)$ such that for all $\kappa\leq \kappa_0$
\begin{align}\label{eq:LimofqandK}
C_2\kappa^4+ \frac{\kappa^2}{16}\leq  q\leq \frac{\kappa^2}{16}+C_3 \kappa^4,\qquad \frac{\pi}{2} +C_2\kappa^2\leq K(\kappa)\leq \frac{\pi}{2}+C_3\kappa^2.
\end{align}
When $\kappa\leq \kappa_0$, plugging \eqref{eq:LimofqandK} into \eqref{eq:FormOfCd} yields
\[\Big|\mathrm{cd}(z,k)- \cos\big(\pi z/2K(\kappa)\big)\Big|\leq \sum_{n=1}^{\infty}\frac{(C_3\kappa)^{2n}}{1-\frac{\kappa^2}{16}-C_2\kappa^4}+C_2\kappa^2.\]
For small enough $\kappa$, the right hand side of the above inequality is bounded by $C_1\kappa^2$ (for some constant $C_1$) which proves \eqref{eq:CDSausage}.
\medskip

Owing to \eqref{eq:KappaLimit}, one has $\frac{1-\kappa}{1+\kappa}=: \sqrt{\tau/\pi}+ \mathfrak{H}_1(\tau)$ such that $|\mathfrak{H}_1(\tau)|\leq C_5\tau$ for small enough $\tau$ where $C_5>0$ is a constant. It is worth noting that $\tau=(-x)^{-\eta}\leq (-x)^{-\eta_0}$ and $\tau$ is converging to $0$ as $x\to -\infty$. By the virtue of \eqref{eq:KappaLimit}, $(1-\kappa)/(1+\kappa)\to 0$ as $x\to -\infty$. Therefore, applying \eqref{eq:CDSausage} and letting $\tilde{\kappa}:=\frac{1-\kappa}{1+\kappa}$, there exist $x_1=x_1(\eta_0)>0$ and $C_6=C_6(\eta_0)>0$ such that
\begin{align}\label{eq:CdToCos}
\mathrm{cd}\left(2(-x)^{3/2} VK(\tilde{\kappa}), \tilde{\kappa}\right)= \cos\big(\pi (-x)^{\frac{3}{2}} V\big)+\mathfrak{H}_2(x;\eta)
\end{align}
where $|\mathfrak{H}_2(x;\eta)|\leq C_6(-x)^{-\eta}$ for all $(-x)\geq x_1$. Combining this and \eqref{eq:KappaLimit} yields
\begin{align}
\sqrt{-\frac{x}{2}}& \frac{1-\kappa}{\sqrt{1+\kappa^2}} \mathrm{cd}\left(2(-x)^{\frac{3}{2}} VK\left(\tilde{\kappa}\right), \tilde{\kappa}\right)= \sqrt{\frac{v}{\pi(-x)^{\frac{1}{2}}}} \cos(\pi (-x)^{\frac{3}{2}} V(\tau))+ J_2(x)\label{eq:cd}
\end{align}
where $|J_2(x)|\leq C_6(-x)^{\frac{1}{2}-\frac{3\eta}{2}}$ for all $(-x)\geq x_1$. Plugging  \eqref{eq:cd} into \eqref{eq:ASSol} along with the inequality $|J_1(x;\gamma)| \leq C_8(-x)^{\frac{1}{2}-\frac{3\eta}{2}}$ (thanks to \eqref{eq:J1Bound1} and $\frac{1}{2}-\frac{3\eta}{2}>-\frac{1}{10}$) yields \eqref{eq:LimitingValue}. By \eqref{eq:VLimit}, there exist $x_2=x_2(\eta_0)>0$ and $C_7=C_7(\eta_0)>0$ such that
\begin{align}\label{eq:tV}
\pi (-x)^{\frac{3}{2}} V(\tau) = - \frac{2}{3}(-x)^{\frac{3}{2}} + \frac{v}{2\pi} \log(8(-x)^{\frac{3}{2}}) - \frac{v}{\pi} \log (v/2\pi) + \frac{v}{2\pi}(1+2\pi) +J_3(x)
\end{align}
where $|J_3(x)|\leq C_7(-x)^{-2\eta}$. By using \eqref{eq:tV} in \eqref{eq:phi}, we get \eqref{eq:FormOfPhi}.
\end{proof}																							

For the next lemma, we will treat $v=s^{\frac{3}{2}-\delta}$ as constant. This may seem to be contrary to the formulation of Lemma~\ref{lem:cosform}. However, as we explain in the beginning of the proof of Lemma \ref{lem:FinalLemma}, we may set $s^{\frac{3}{2}-\delta}= (-x)^{\frac{3}{2}-\eta}$ where $\eta=\eta(x)$ is chosen so as to match both sides. Since the result of Lemma \ref{lem:cosform} is stated as true uniformly over varying $\eta$, we it remains valid.

\bl\label{lem:FinalLemma}
 Fix $\delta\in(0,2/3)$. Set $v=s^{\frac{3}{2}-\delta}$ (see the discussion at the beginning of the proof). Recall $\phi(\cdot)$ from \eqref{eq:phi} and define
\begin{equation}\label{eq:psi}
\psi(x)=-\frac{4}{3}(-x)^{\frac{3}{2}}+ \frac{v}{\pi} \log(8(-x)^{\frac{3}{2}})+ 2\phi(x).
\end{equation}
Choose $\theta\in (0,\delta)$ such that $(\delta-\theta)\in (0,2/5)$.
 Then, there exist $s_0=s_0(\theta)>0$ and $C=C(\theta)>0$ such that for all $s\geq s_0$
\begin{align}\label{eq:OscillatoryTerm}
\Big|\int^{-s^{1-\frac{2}{3}\theta}}_{-s} \frac{(x+s)}{(-x)^{\frac{1}{2}}} \cos\left(\psi(x)\right) dx\Big| \leq Cs^{\frac{3}{2}}\max\{s^{-\delta+\theta},s^{-\frac{3}{2}+\theta)},s^{-2(\delta-\frac{2}{3}\theta)(1-\theta)}, s^{-\frac{5}{2}(1-\frac{2}{3}\theta)-(\delta-\theta)}\}.&&&&
\end{align}
\el
\begin{proof}
First, observe that for $x\in [-s, -s^{1-\frac{2}{3}\theta})$, there exists $\eta=\eta(x)\in (\delta-\theta,2/5)$ such that $s^{\frac{3}{2}-\delta}= (-x)^{\frac{3}{2}-\eta(x)}$. Letting $\eta_0=\delta-\theta$ (which is $\leq 2/5$) we have that $\eta=\eta(x)\geq \eta_0$ for $x\in [-s, -s^{1-\frac{2}{3}\theta})$. Thus, we can apply Lemma~\ref{lem:cosform} to conclude that there exist $s_0=s_0(\delta,\theta)>0$ and $C^{\prime}=C^{\prime}(\delta,\theta)>0$ such that for all $s\geq s_0$
\begin{equation}\label{eq:PhiBound}
|J_3(x)|= \left|\frac{2\pi}{v}\phi(x)- (1+2\pi)+\log(v)\right|\leq C^{\prime}(-x)^{-(\delta-\theta)}
\end{equation}
for any $x\in (-s, -s^{1-\theta})$ where $J_3(\cdot)$ is the same as in \eqref{eq:phi}.

Now, in order to show \eqref{eq:OscillatoryTerm} we will divide the interval of integration $[-s,-s^{1-\frac{2}{3}\theta}]$ into the unique disjoint union of consecutive closed intervals $\mathcal{I}_1,\mathcal{I}_2, \ldots , \mathcal{I}_{k}$ (here $k$ is some non-negative integer) such that these intervals satisfy the following two conditions: (1) The right end point of $\mathcal{I}_1$ is $-s^{1-\frac{2}{3}\theta}$ and left end point of $\mathcal{I}_k$ is $-s$; (2) For any $1\leq j<k$, if $\mathcal{I}_j=[a_j,b_j]$ then,
\begin{equation}\label{eq:pqrel}
a_j= b_j-\pi(-b_j)^{-\frac{1}{2}}.
\end{equation}

Fix any $1\leq j<k$ and write $\mathcal{I}_j=[a,b]$ where $a=a_j$ and $b_j$ satisfy \eqref{eq:pqrel}. Then, for any $t\in [0,1]$, using Taylor's expansion and \eqref{eq:PhiBound} (to get the second equality in \eqref{eq:Label2}) we find that there exist $s_1=s_1(\delta,\theta)>0$ and $C=C(\delta,\theta)>0$ such that for all $s\geq s_1$,
\begin{align}
&\frac{4}{3}\big(-b+(-b)^{-\frac{1}{2}}t\pi\big)^{\frac{3}{2}}= \frac{4}{3}(-b)^{\frac{3}{2}}+ 2\pi t + J_4(b),\qquad\qquad
\phi\big(-b+ (-b)^{-\frac{1}{2}}t\pi\big) = \phi(-b)+ J_5(b) \\
&\qquad\frac{v}{\pi} \log \big(8\big(-b+ (-b)^{-\frac{1}{2}}t\pi\big)^{\frac{3}{2}}\big) = \frac{v}{\pi} \log (8 (-b)^{\frac{3}{2}}) + \frac{3vt}{2} (-b)^{-\frac{3}{2}}+J_6(b)\label{eq:Label2}
\end{align}
where
\begin{align}
|J_4(b)|\leq C(-b)^{-\frac{3}{2}}, \quad |J_5(b)|\leq C(-b)^{-2(\delta-\theta)}, \quad |J_6(b)|\leq C(-b)^{-\frac{5}{2}-(\delta-\theta)/(1-\frac{2}{3}\theta)}\label{eq:3Bound}.
\end{align}
Therefore, for any $t\in [0,1]$, combining  \eqref{eq:Label2} and \eqref{eq:3Bound}, we arrive at
\begin{equation}\label{eq:psidecomp}
\psi(b-(-b)^{-\frac{1}{2}}t\pi)= \psi(b)+2\pi t+\frac{3vt}{2}(-b)^{-\frac{3}{2}} +J_7(b)
\end{equation}
where for some constant $C=C(\delta,\theta)>0$ and all large enough $s$
\begin{equation}\label{eq:J8Bound}
|J_7(b)|\leq C\max\big\{(-b)^{-\frac{3}{2}},(-b)^{-2(\delta-\theta)}, (-b)^{-\frac{5}{2}-(\delta-\theta)/(1-\frac{2}{3}\theta)}\big\}.
\end{equation}
\medskip
\noindent\textbf{Claim:}
There exist $s_3=s_3(\delta,\theta)>0$ and $C=C(\delta,\theta)>0$ such that for all $s\geq s_3$ and all intervals $[a,b]$ with $a=b-(-b)^{-\frac{1}{2}}\pi$ and $-s<a<b<-s^{1-\frac{2}{3}\theta}$,
\begin{align}\label{eq:EachIntCon}
\int_{a}^{b}\frac{(x+s)}{(-x)^{\frac{1}{2}}}\cos\big(\psi(x)\big)dx = \frac{\pi(b+s)}{(-b)}\Big(\sin(\psi(b))-\sin\big(\psi(b)+\frac{3v}{2}(-b)^{-\frac{3}{2}}\big)+ J_8(b)\Big) &&
\end{align}
where $$|J_8(b)|\leq C\max\big\{(-b)^{-\frac{3}{2}},(-b)^{-2(\delta-\theta)}, (-b)^{-\frac{5}{2}-(\delta-\theta)/(1-\frac{2}{3}\theta)}\big\}.$$

\noindent\textsc{Proof of Claim:}
Note that any point in the interval $[a,b]$ can be written as $b-(-b)^{-\frac{1}{2}}t\pi$ for some $t\in [0,1]$. We use a shorthand notation $\psi_t(b)$ for $\psi(b)-2\pi t+ (-b)^{-\frac{3}{2}}vt$. Owing to \eqref{eq:psidecomp}, we have that
$\psi(b-(-b)^{-\frac{1}{2}} t\pi)= \psi_t(b)+J_7(b)$ where $J_7(b)$ satisfies \eqref{eq:J8Bound}.
Applying the formula $\cos(x+y)= \cos(x)\cos(y) - \sin(x) \sin(y)$,
\begin{align}\label{eq:cosexpansion}
\cos\big(\psi(b-t\pi(-b)^{-\frac{1}{2}})\big)= \cos\big(\psi_t(b)\big)\cos\big(J_7(b)\big)-\sin\big(\psi_t(b)\big)\sin\big(J_7(b)\big).
\end{align}
Appealing to \eqref{eq:J8Bound} and the mean value theorem shows that for some $C>0$,
\begin{equation}\label{eq:CosSinBound}
\max\big\{|\cos\big(J_7(b)\big)-1|,|\sin\big(J_7(b)\big)|\big\} \leq  C\max\big\{(-b)^{-\frac{3}{2}},(-b)^{-2(\delta-\theta)}, (-b)^{-\frac{5}{2}-(\delta-\theta)/(1-\frac{2}{3}\theta)}\big\}.
\end{equation}
By applying virtue of the change of variable $t= \frac{1}{\pi}(-b)^{\frac{1}{2}}(x-b)$, we can show that
\begin{align}\label{eq:integral}
\int^{b}_{a}\frac{(x+s)}{(-x)^{\frac{1}{2}}}\cos\big(\psi(x)\big)dx =\frac{\pi}{(-b)^{\frac{1}{2}}} \int^{1}_{0} \left(\frac{b +s}{(-b)^{\frac{1}{2}}}+\mathfrak{H}_3(b)\right)\cos\big(\psi(b-(-b)^{-\frac{1}{2}}t\pi)\big)dt &&&
\end{align}
where the term $\mathfrak{H}_3(b)$ (which depends on $t$ and is computed from the change of variables) can readily be bounded uniformly in $t$ as $|\mathfrak{H}_3(b)|\leq C(-s)^{-1+\frac{2}{3}\theta}$ for some constant $C>0$.
Plugging \eqref{eq:cosexpansion} and \eqref{eq:CosSinBound} into the right side of \eqref{eq:integral} and evaluating yields \eqref{eq:EachIntCon}.
\medskip

Now, we turn to the final step of the proof where we sum the contributions over all the intervals $\mathcal{I}_1,\ldots,\mathcal{I}_k$. Recall that $\mathcal{I}_j=[a_j,b_j]$ where $a_j=b_j-(-b_j)^{-\frac{1}{2}} \pi$
for $1\leq j<k$. Thus, \eqref{eq:EachIntCon} holds for $a=a_j$ and $b=b_j$ for all $1\leq j<k$. Summing \eqref{eq:EachIntCon} over $1\leq j<k$ and bounding $\big|\sin\big(\psi(b)\big)- \sin\big(\psi(b)+\tfrac{3}{2}v(-b)^{-\frac{3}{2}}\big)\big|$ by $\tfrac{3}{2}v(-b)^{-\frac{3}{2}}$, we get
\begin{equation}\label{eq:starstar}
\Big|\sum_{j=1}^{k-1} \int_{\mathcal{I}_j}(-x)^{-\frac{1}{2}}(x+s) \cos\big(\psi(x)\big)dx\Big| \leq \sum_{j=1}^{k-1} \frac{\pi(b_j+s)}{(-b_j)}\Big(\frac{3v}{2}(-b_j)^{-\frac{3}{2}}+ J_8(b_j) \Big).
\end{equation}
Using the bound $\sum_{j=1}^{k-1}\frac{\pi(b_j+s)}{(-b_j)} \leq 2 \int^{-s^{1-\theta}}_{-s}(-x)^{-\frac{1}{2}}(x+s) dx = \frac{8s^{\frac{3}{2}}}{3}$ yields the following
\begin{align}\label{eq:SumInt}
\text{r.h.s of \eqref{eq:starstar}}\leq \frac{8s^{\frac{3}{2}}}{3\pi}\Big(\frac{3}{2}\max_{1\leq j<k} v(-b_j)^{-\frac{3}{2}}+ \max_{1\leq j<k} |J_8(b_j)|\Big).
\end{align}
Since $v=s^{\frac{3}{2}-\delta}$ and $(-b_j)\geq s^{1-\frac{2}{3}\theta}$, we find that $v(-b_j)^{-\frac{3}{2}}\leq s^{-\delta+\theta}$ for all $1\leq j<k$. Likewise,
\begin{align}\label{eq:J9Max}
\max_{1\leq j<k} |J_8(b_j)|\leq C\max\big\{(s^{1-\frac{2}{3}\theta})^{-\frac{3}{2}},(s^{1-\frac{2}{3}\theta})^{-2(\delta-\theta)}, (s^{1-\frac{2}{3}\theta})^{-\frac{5}{2}-(\delta-\theta)/(1-\frac{2}{3}\theta)}\big\}.
\end{align}
 Combining \eqref{eq:SumInt}, \eqref{eq:J9Max} and observing that $\big|\int_{\mathcal{I}_k} (-x)^{-\frac{1}{2}}(x+s) \cos\big(\psi(x)\big)dx\big|\leq 1$ (this follows since $0\leq x+s\leq (-b_k)^{-\frac{1}{2}}\pi$, $b_k-a_k\leq  (-b_k)^{-\frac{1}{2}}\pi$, and $-b_k<-s^{1-\theta}$), we arrive at \eqref{eq:OscillatoryTerm}.
\end{proof}

\subsection{Proof of Theorem~\ref{thm:CGFexpansion}}\label{sec:CGF}
We divide the integral $\int_{-s}^{\infty}(x+s)u^2_{\mathrm{AS}}(x;\gamma) dx$  into two parts $(\mathbf{a})$ and $(\mathbf{b})$
\begin{align}\label{eq:SepInt}
(\mathbf{a}) := \int^{\infty}_{0}(x+s) u^{2}_{\mathrm{AS}}(x;\gamma) dx \quad \text{and, }\quad (\mathbf{b}):= \int^{0}_{-s} (x+s) u^2_{\mathrm{AS}} (x;\gamma) dx.
\end{align}
Owing to the exponential decay (see \eqref{eq:ASBoundary}) of $u_{\mathrm{AS}}(x;\gamma)$ as $x\to \infty$, $(\mathbf{a}) = \mathcal{O}(s)$ (i.e., $(\mathbf{a})$ is bounded above by $Cs$ as $s$ grows to $\infty$). To estimate the value of $(\mathbf{b})$, fix some $\theta \in (0,\delta)$ such that $(\delta-\theta)\in (0,2/5)$ and divide $(\mathbf{b})$ into
\begin{align}\label{eq:DivInt}
(\mathbf{b1}):= \int^{0}_{-s^{1-\frac{2}{3}\theta}} (x+s)& u^2_{\mathrm{AS}}(x;\gamma)dx, \quad
&(\mathbf{b2}):= \int^{-s^{1-\frac{2}{3}\theta}}_{-s}(x+s)u^2_{\mathrm{AS}}(x;\gamma) dx.
\end{align}
\medskip
\noindent\textbf{Claim:}
There exist $s_0=s_0(\delta,\theta)>0$ and $C=C(\delta,\theta)>0$ such that for all $s\geq s_0$,
\begin{align}\label{eq:II3Bound}
(\mathbf{b2}) = \frac{2v}{3\pi}s^{\frac{3}{2}}\big(1+J_{9}(s)\big) + J_{10}(s)
\end{align}
where $|J_{10}(s)|\leq Cs^{3-\frac{5(\delta-\theta)}{2}}$ and
\begin{equation}\label{eq:J10Bound}
|J_{9}(s)|\leq C\max\{s^{-\frac{2}{3}\theta},s^{-(\delta-\theta)},s^{-\frac{3}{2}+\theta},s^{-2(\delta-\theta)(1-\frac{2}{3}\theta)}, s^{-\frac{5}{2}(1-\frac{2}{3}\theta)-(\delta-\theta)}\}.
\end{equation}

\noindent\textsc{Proof of Claim}:
As explained in the proof of Lemma \ref{lem:FinalLemma}, for any $x\in (-s, -s^{1-\frac{2}{3}\theta})$, there exists $\eta=\eta(x)\in (\delta-\theta,2/5)$ such that $v=s^{\frac{3}{2}-\delta}=(-x)^{\frac{3}{2}-\eta(x)}$. Squaring both sides of \eqref{eq:LimitingValue} shows that there exist $x_0=x_0(\delta,\theta)>0$ and $C=C(\delta,\theta)>0$ such that for all $x\geq x_0$
\begin{equation}\label{eq:usquare}
u^2_{\mathrm{AS}}(x;\gamma) = \frac{1}{(-x)^{\frac{1}{2}}}\frac{v}{\pi}\cos^2\left(-\frac{2}{3}(-x)^{\frac{3}{2}}+ \frac{v}{2\pi} \log(8(-x)^{\frac{3}{2}})+ \phi(x)\right) + J_{11}(x)
\end{equation}
where $|J_{11}(x)|\leq C(-x)^{1-\frac{5(\delta-\theta)}{2}}$.
Now, we plug \eqref{eq:usquare} inside the integral of $(\mathbf{b2})$. Recalling $\psi(\cdot)$ from \eqref{eq:psi} and using the identity $\cos^2(z)= \frac{1}{2}(\cos 2z+ 1)$, we arrive at
   \begin{align}\label{eq:CiNTEGRAL}
   (\mathbf{b2})
    =\frac{v}{2\pi}\int^{-s^{1-\theta}}_{-s} \frac{(x+s)}{(-x)^{\frac{1}{2}}}\Big(1+ \cos\big(\psi(x)\big)\Big) dx + \int^{-s^{1-\theta}}_{-s} (x+s)J_{11}(x)dx.
\end{align}
By a direct computation, it follows
\begin{equation}\label{eq:LeadTerm}
\frac{v}{2\pi}\int^{-s^{1-\theta}}_{-s} \frac{1}{(-x)^{\frac{1}{2}}}(x+s)dx=\frac{2v}{3\pi}s^{\frac{3}{2}}(1+\mathfrak{H}_4(s)).
\end{equation}
where $|\mathfrak{H}_4(s)|\leq C s^{-\frac{2}{3}\theta}$ for some $C>0$.
Owing to the upper bound on $J_{11}(x)$, we get
\begin{equation}\label{eq:BoundOnJ12}
\Big|\int^{-s^{1-\theta}}_{-s}(x+s)J_{11}(x) dx\Big|\leq C s^{3- \frac{5(\delta-\theta)}{2}}
\end{equation}
Plugging \eqref{eq:LeadTerm} and \eqref{eq:BoundOnJ12} into \eqref{eq:CiNTEGRAL} and invoking \eqref{eq:OscillatoryTerm}, we find \eqref{eq:II3Bound}.

Now, we are ready to complete the proof of Theorem~\ref{thm:CGFexpansion}. Note that  $(\mathbf{a})$ and $(\mathbf{b1})$ are both positive real numbers. Owing to this and \eqref{eq:II3Bound}, we observe
	\begin{align}\label{eq:Finaleq}
\int_{-s}^{\infty}(x+s)& u^2_{\mathrm{AS}}(x;\gamma) dx = (\mathbf{a})+(\mathbf{b})=(\mathbf{a}) +(\mathbf{b1}) + (\mathbf{b2})\geq  \frac{2v}{3}s^{\frac{3}{2}} +\mathfrak{R}.
\end{align}		
where
$
\mathfrak{R}:= \frac{2v}{3}s^{\frac{3}{2}}J_{9}(x)+J_{10}(x).
$
Finally, we set $\theta= 2\delta/5$. As $\delta \in (0,2/3)$, necessarily $\delta-\theta=3\delta/5\in (0,2/5)$. Plugging $\theta$ into \eqref{eq:J10Bound}, we get $|J_{9}(s)|\leq Cs^{-\frac{4\delta}{15}}$ and $|J_{10}(s)|\leq Cs^{3-\frac{3\delta}{2}}$.
Combining these inequalities (and recalling $v=s^{\frac{3}{2}-\delta}$) yields $|\mathfrak{R}|\leq C^{\prime}s^{3-\frac{19\delta}{15}}$ for all large enough $s$ where $C^{\prime}>0$ is a constant. This completes the proof of Theorem~\ref{thm:CGFexpansion}.
\qed


\bibliographystyle{alpha}
\bibliography{Reference}

\end{document}